%% file: coresizes.tex
\documentclass[a4paper,reqno]{paper}
\usepackage[utf8]{inputenc}
\usepackage[T1]{fontenc}
\usepackage{amsfonts,enumerate}
\usepackage[usenames, dvipsnames, table ]{xcolor}
\usepackage{amsmath}

\usepackage{amssymb}
\usepackage{amsthm}
\usepackage{graphics}
 \usepackage{array,multirow,graphicx}
\usepackage[backref = page]{hyperref}
\hypersetup{
colorlinks =  true, 
linkcolor = purple!80, 
citecolor =teal
}
\usepackage[capitalize]{cleveref}
\usepackage{tikz}
\usepackage{geometry}
\usepackage{graphicx}
\usepackage{faktor}
\usepackage{genyoungtabtikz}
\usepackage{dsfont}
\usepackage{mathrsfs}
\usepackage{nicematrix}
\usepackage{stmaryrd}
\usepackage{mynewcommands}
\usepackage{booktabs}
\usepackage{ytableau}
\usepackage{pgf}
\usetikzlibrary{calc}
\usetikzlibrary{positioning}
\usepackage{float}
\usepackage{titlesec}
\usetikzlibrary{arrows, backgrounds}
\usepackage{pdflscape}
\usepackage{rotating}
\usepackage{physics}
\usepackage[new]{old-arrows}
\usepackage[all]{xy}
\usepackage{pgfplots}
\pgfplotsset{compat=1.18} 
\usepackage[justification=centering]{caption}
\usepackage{paralist}
\usepackage{dirtytalk}
\usepackage{url}
\usepackage{shuffle}
\usepackage{blkarray,bigstrut}
\usepackage{nicematrix}
\usepackage{tcolorbox}
\usepackage{makecell}
\usepackage{afterpage}

\usepackage{frcursive}

\usetikzlibrary{decorations.pathreplacing,decorations.markings}
\usetikzlibrary{arrows.meta}
\usepackage{longtable}
\Crefname{equation}{}{}
\usepackage{titletoc,tocloft}
\titleformat{\subsection}
{\normalfont\fontsize{11}{17}\sffamily\bfseries}
{\thesubsection}
{1em}
{}

\let\svthefootnote\thefootnote
\newcommand\freefootnote[1]{
  \let\thefootnote\relax
  \footnotetext{#1}
  \let\thefootnote\svthefootnote
}

\makeatletter
\DeclareRobustCommand*{\mfaktor}[3][]
{
   { \mathpalette{\mfaktor@impl@}{{#1}{#2}{#3}} }
}
\newcommand*{\mfaktor@impl@}[2]{\mfaktor@impl#1#2}
\newcommand*{\mfaktor@impl}[4]{
   \settoheight{\faktor@zaehlerhoehe}{\ensuremath{#1#2{#3}}}%
   \settoheight{\faktor@nennerhoehe}{\ensuremath{#1#2{#4}}}%
      \raisebox{-0.5\faktor@zaehlerhoehe}{\ensuremath{#1#2{#3}}}%
      \mkern-4mu\diagdown\mkern-5mu%
      \raisebox{0.5\faktor@nennerhoehe}{\ensuremath{#1#2{#4}}}%
}
\makeatother

\makeatletter
\renewcommand*\env@matrix[1][*\c@MaxMatrixCols c]{%
  \hskip -\arraycolsep
  \let\@ifnextchar\new@ifnextchar
  \array{#1}}
\makeatother

\renewcommand\arraystretch{1.7}

\newcommand{\SC}{\mathcal{SCC}}

\begin{document}

\sloppy

\title{Generalised core partitions and Diophantine equations}

\author{
Olivier Brunat\thanks{
Universit\'e Paris Cit\'e, Institut de math\'ematiques de
         Jussieu -- Paris Rive Gauche,  UFR de math\'e\-matiques, Case
7012, 75205 Paris Cedex 13, France. Email address: {\tt olivier.brunat@imj-prg.fr}
},
Nathan Chapelier-Laget\thanks{Université du Littoral Côte d'Opale, 62228 Calais, France. Email address:  {\tt nathan.chapelier@univ-littoral.fr}},
Thomas Gerber\thanks{Université Lyon 1, 69100 Villeurbanne, France. Email address: 
{\tt gerber@math.univ-lyon1.fr}
}
}
\maketitle

\begin{flushright}
\textit{To Jérémie}
\end{flushright}
\vspace{1cm}

\hrule 

\begin{abstract} 
We study generalised core partitions arising from 
affine Grassmannian elements in arbitrary Dynkin type.
The corresponding notion of size is given by the atomic length in the sense of \cite{CG2022}.
In this paper, we first develop the theory for extended affine Weyl groups.
In a series of applications, we give some remarkable parametrisations of the solutions of certain Diophantine equations
resembling Pell's equation, by refining the results of \cite{BrunatNath2022} and \cite{Alpoge2014}, and generalising them to
further  types.
\end{abstract}

\hrule

\setcounter{tocdepth}{2}
 \tableofcontents

\section*{Introduction}
\addcontentsline{toc}{section}{\protect\numberline{}Introduction}

\subsection*{The ubiquity of core partitions}

The study of core partitions, those integer partitions with no hooks of a given size, is a recurring theme in 
algebraic combinatorics. Over the years, it has provided
striking links between different fields of mathematics such as
\begin{itemize}
\item representation theory, originally via the Nakayama conjecture \cite{Brauer1947, Robinson1947}
which paved the way to investigating block theory for symmetric groups
\cite{JamesKerber1984}
and related algebraic structures 
such as Hecke algebras \cite{GeckPfeiffer2000, Fayers2006, GeckJacon2011}, finite groups of Lie type \cite{FS1982, Olsson1986},
\item enumerative combinatorics, for instance via the various generalisations and interpretations
of the Rogers-Ramanujan identities \cite{AndrewsEriksson2004, Gordon1961, Andrews1974, Bressoud1979, GKS1990, HirschhornSellers1996b, Chen2013};
or via some famous enumeration results 
\cite{GO1996, Anderson2002, Han2010, HanOno2011, Johnson2018},
\item number theory, via the study of certain modular forms
arising from the generating series of core partitions
\cite{GO1996, HirschhornSellers1996a, HirschhornSellers1999} 
and relationships with certain Diophantine equations \cite{OnoSze1997, HanusaNath2013, Alpoge2014}.
\end{itemize}

Recently, there has been some exciting developments 
in all of these directions, for instance
with the study of blocks for cyclotomic Hecke algebras \cite{JaconLecouvey2021, ChlouverakiJacon2023};
generalisations of classical enumeration results \cite{Fayers2019, STW2023, LW2024};
or the explicit parametrisation of the solutions of certain Pell-type equations \cite{BrunatNath2022, MS2023, HansonJameson2022}.
All of these results rely on a careful study of the combinatorics
of cores, either in their classical definition or as a natural generalisation.
In particular, the notion of (generalised) size of a 
(generalised) core plays an essential role.

\subsection*{The atomic length: a new approach to core partitions}

In previous work, the last two authors have introduced 
a statistic, called the \textit{atomic length}, defined on
the Weyl group attached to a Kac-Moody algebra of given finite or affine type \cite{CG2022}.
This approach, representation-theoretic by nature, enables us to 
study cores in greater generality and in a very natural fashion. 
Moreover, it was shown that appropriate specialisation of the parameters
defining the atomic length enables to recover previous constructions in the literature \cite{GKS1990, Lascoux2001, STW2023}.
For instance, one recovers the usual $(n+1)$-core partitions and their size simply by considering type $A_n^{(1)}$ 
and specialising at the fundamental weight $\La_0$.

\medskip

In arbitrary type, specialising at $\La_0$ yields a statistic on so-called
\textit{affine Grassmannian} elements, which we will therefore call \textit{generalised cores} in what follows. 
This is consistent with the philosophy of \cite{STW2023}.
In this paper, it will be relevant
to consider the more general $\La$-atomic length. 
For instance, in type $C_2^{(1)}$, specialising at the fundamental weight $\La=\La_1$ will be particularly fruitful.

\medskip

In \cite{CG2022}, much of the study focused on Kac-Moody algebras of \textit{finite} types. This culminated in an analogue of the
famous Granville--Ono theorem \cite{GO1996} for all finite Dynkin types.
The case of \textit{affine} Kac-Moody algebras was initiated, but not extensively investigated.

\medskip

The first objective of this paper is to define and study the \textit{extended affine Weyl group}
in the context of Kac-Moody algebras, see \Cref{sec_ext_AG}, which, as far as we know, had not been done in the literature before.
This leads to the definition of \textit{extended atomic length} in \Cref{sec_ext_AL}.
Using both the decomposition of the extended affine Weyl group as the semidirect product
of the finite Weyl group with the coweight lattice and as the semidirect product of the fundamental group with the affine Weyl group, we show that, with the exception of twisted type $A_{2n}^{(2)}$, the atomic length $\sL_{\La_0}$ admits the following remarkable property, see \Cref{Theorem uniform extended atomic length}.
\begin{center}
\begin{tcolorbox}[width=15cm, boxrule=0mm, colback=white]
\centering
For any $\sigma$ in the fundamental group,
$\sL_{\La_0}(\sigma w)= \sL_{\La_0}(w)$ for all Weyl group element $w$.
\end{tcolorbox}
\end{center}
This is a direct analogue of the well-known formula $\ell(\sigma w)=\ell(w)$, which characterises the fundamental group using the Coxeter length.

\subsection*{Labelling solutions of Pell-type equations by generalised cores}

The second goal of this paper is
to demonstrate a striking application of this theory: in several Dynkin types of small rank,
we prove bijectively that the generalised cores
parametrise the solutions of certain Pell-type Diophantine equations, namely
\begin{itemize}
    \item in type $A_2^{(1)}$, this was proved in \cite{BrunatNath2022}, which we recall in \Cref{thm_A2}, 
    \item in type $C_2^{(1)}$, see \Cref{thm_C2},
    \item in type $C_2^{(1)}$ for a second, complementary equation, see \Cref{thm_C2_1},
    \item in type $D_3^{(2)}$, see \Cref{thm_D3t},
    \item in type $D_4^{(3)}$, see \Cref{thm_D43}.
\end{itemize} 
These results, all in rank $2$, rely on the same following mechanism.
\begin{center}
\begin{tcolorbox}[width=15cm, boxrule=0mm, colback=white]
\begin{itemize}
    \item[\textbf{Step 1}] Consider the Diophantine equation of the form 
    \begin{equation}\label{Pell} x^2+dy^2=k \end{equation}
    where $d\in\{1,3\}$ and $k$ has a certain fixed residue class (these will depend on the Dynkin type).
    \item[\textbf{Step 2}]  Endow the solution set of \Cref{Pell} with a \textit{free action} of an appropriate finite group $G$.
    \item[\textbf{Step 3}]  Find an explicit map $\varphi$ that sends (extended) affine Grassmannian elements to solutions of \Cref{Pell}.
    \item[\textbf{Step 4}]  Show that $\Im(\varphi)$ is a \textit{complete set of representatives} of the $G$-orbits.
\end{itemize}
\end{tcolorbox}
\end{center}

\medskip

Let us focus more closely on Equation \Cref{Pell} in the case $d=1$. 
First, considering the affine Grassmannian elements in type $C_2^{(1)}$
of atomic length equal to $N$ (which are naturally in bijection with
the self-conjugate $4$-cores of size $N$), we obtain a description of the solutions of
the Diophantine equation $x^2+y^2=8N+5$ by the above process (see
\Cref{thm_C2}). In particular, this gives a simple bijective proof of an identity due to Alpoge \cite[Theorem 7]{Alpoge2014}.
Moreover, in contrast with Alpoge's approach, that requires involved number-theoretic methods, 
our proofs are quite elementary once the machinery of the atomic length
is properly established.
Secondly, using the atomic length associated with the
fundamental weight $\Lambda_1$ in type $C_2^{(1)}$, we obtain a
description of the solutions of the Diophantine equation $x^2+y^2=8N+1$
in terms of extended Grassmannian elements
(see \Cref{thm_C2_1}).
Combining these two results, we obtain the following.
\begin{center}
\begin{tcolorbox}[width=15cm, boxrule=0mm, colback=white]\centering
For any odd $k$, the solutions of $x^2+y^2=k$ are in correspondence with
the (extended) affine Grassmannian elements of type $C_2^{(1)}$.
\end{tcolorbox}
\end{center}
Note that, by \Cref{lembij}, this gives a way to
describe the integral solutions of $x^2+y^2=k$ for \textit{any} positive integer $k$, which
is of independent interest.

\medskip

In addition, we can use the combinatorial description of
the affine Grassmannian elements in terms of integer partitions of \cite{STW2023} and \cite{LW2024}
to obtain neat enumeration corollaries in further Dynkin types: \Cref{cor_D3t}, \Cref{cor_D43}.
For \Cref{cor_C2_1}, the question of finding a similar combinatorial model is open.

\medskip

Moreover, in type $A_n^{(1)}$ of small ranks,
we shed some light on the crank phenomenons revealed in \cite{BrunatNath2022} 
by using the \textit{extended} atomic length on 
the extended affine Weyl group (\Cref{sec_sols_rank2}).
This gives a completely natural refinement of the parametrisation of the solutions of the associated Diophantine equations
by using more involved group actions.
Type $A_3^{(1)}$ is addressed in \Cref{sec_pell_A3},  where the structure of the integral solutions of the equation 
\begin{equation}\label{eqpigA3}
x^2+2y^2+3z^2 =48N+30
\end{equation}
arising from extended affine Grassmannian elements is elucidated. 
More precisely, there is a dihedral group $G$ of order 12, acting on the integer solutions of \Cref{eqpigA3} 
such that, see \Cref{pig A3},
\begin{center}
\begin{tcolorbox}[width=15cm, boxrule=0mm, colback=white]\centering
Every solution of \Cref{eqpigA3} arises from a $4$-core modulo the action of $G$. 
\end{tcolorbox}
\end{center}
As a consequence, using results of Bhargava \cite{bhargava2000conway}, we obtain an alternative, constructive proof of
the Granville-Ono theorem \cite{GO1996} in rank $3$, claiming the existence of $4$-cores of any size, see \Cref{thm 4-cores}.

\subsection*{Structure of the paper}

In \Cref{sec_KM}, we describe the starting point of our theory,
namely the construction of Weyl groups associated to an affine
Kac-Moody algebra (or equivalently, a generalised Cartan matrix of affine type).
We take advantage of this section to explain in detail which conventions and constructions are necessary in our setting, by referring mainly to Carter's book \cite{Carter2005}.

\medskip

\Cref{sec_ext} is dedicated to the study of the \textit{extended} version
of the affine Weyl group and affine Grassmannian.
It contains extensions to all twisted Kac-Moody types of several classical results,
which did not seem to appear elsewehere in the literature.
In particular, we give some explicit bijections which will provide
us with different ways to think about (extended) affine
Grassmannian elements.

\medskip

\Cref{sec_AL} is concerned with the general study of the \textit{affine atomic length} \cite{CG2022} and its \textit{extended} version, which we introduce here.
This statistic depends on a parameter $\La$ which is a dominant weight, but we will focus mainly on the case $\La=\La_0$ (and sometimes on the fundamental weight $\La=\La_i$, $i\geq 1$.)
We establish new results that will be reinvested in subsequent sections, as the extended
atomic length will naturally appear in some of our main results.
In particular, 
we show in \Cref{decompo atomic D} that the extended affine atomic length is compatible
with the semidirect decomposition of the extended affine Weyl group;
and we derive \cref{Theorem uniform extended atomic length}, which is central.
Finally, in \Cref{section gaussian equation gen}, we explain how we can systematically attach solutions of
certain Diophantine equations to (extended) affine Grassmannian elements.

\medskip

In \Cref{sec_cores}, we take a breath and summarise well-known results on the combinatorics of $(n+1)$-cores using
their interpretation as affine Grassmannian elements of type $A_n^{(1)}$, essentially following Lascoux's approach \cite{Lascoux2001}.
Though short, this section is important as it legitimises the definition of the generalised cores as
the affine Grassmannian elements in arbitrary affine type.
Moreover, we give a closed formula for the Diophantine equation
naturally arising in this case,
and show that extended affine Grassmannian elements provide in fact \textit{integral} solutions.

\medskip

We quickly explain in \Cref{sec_comb_models}  how combinatorial models for generalised cores can be constructed by reviewing the recent results of \cite{STW2023} and \cite{LW2024}.
These constructions are reminiscent of the type $A_n^{(1)}$ combinatorics, but looking at certain self-conjugate core partitions instead.

\medskip

In \Cref{sec_pell}, we temporarily forget about 
generalised cores and focus solely on the Pell-type equations that will naturally arise in \Cref{sec_param} in the rank $2$ cases, namely
$$x^2+y^2 = k \mand x^2+3y^2=k,$$
where $k$ has a fixed congruence class modulo a certain integer (which depends on the type).
We first focus on the first equation and explain how this
can be reduced to the case where $k=1\mod 4$.
Then, we observe that the solution sets of these Diophantine equations are
endowed with certain finite group actions, which turn out to be free
under favourable arithmetic conditions.
This is similar to the philosophy of cranks \cite{GKS1990} used recently in 
\cite{BrunatNath2022}.

\medskip

\Cref{sec_param} is concerned with rank $2$ Dynkin types and
contains many new results. 
There, using the universal approach given by the atomic length on (extended) affine Weyl groups, we are able to label
the solutions of families of Pell-type Diophantine equations by
affine Grassmannians elements (of appropriate Dynkin type).
First, we focus on type $A_2^{(1)}$ and refine the results of \cite{BrunatNath2022}.
Then, we establish several analogous results in other types,
and we are even able to use the combinatorial models
of \Cref{sec_comb_models} to give elegant analogues of
\cite{BrunatNath2022, Alpoge2014}.
This is achieved in type $A_2^{(1)}$ (\Cref{thm_A2}), $C_2^{(1)}$ (for two
different Diophantine equations, \Cref{thm_C2} and \Cref{thm_C2_1}), 
$D_3^{(2)}$  (\Cref{thm_D3t}) and
$D_4^{(3)}$ (\Cref{thm_D43}).

\medskip

In \Cref{sec_pell_A3}, we focus on type $A_3^{(1)}$.
This turns out to be much more involved than its rank $2$ counterpart studied in \cite{BrunatNath2022},
but we are nevertheless able to establish a series of results in the same spirit.
More precisely, we introduce an appropriate explicit group $G$ (isomorphic to the dihedral group of order $12$)
acting on the solutions of the Diophantine equation $$x^2+2y^2+3z^2=48N+30$$ for any fixed $N\in\N$.
Unlike in type $A_2^{(1)}$, the $G$-action is not free, but almost: we are able to describe the $G$-orbits precisely (\Cref{prop orbites}, \Cref{thm size orbits}).
In the main result of that section, \Cref{pig A3}, we fully describe
the relationship between integral solutions and extended affine Grassmannian elements.
From this,  we are able to give a short and constructive proof of the fact that there are $4$-cores of any given size 
(the celebrated result of \cite{GO1996}), this is \Cref{thm 4-cores}.

\medskip

For completeness, we finish in \Cref{sec_Pell_hyperoct} 
with general considerations in cases where the underlying finite type is $B_n$ or $C_n$
(there are $5$ such families),
where the action of the hyperoctahedral group in dimension $n$ is considered.

\subsection*{Perspectives}

One can expect to get analogous parametrisations for
further Diophantine equations, by playing around with the other degrees of freedom, namely by
\begin{itemize}
\item investigating the \textit{other Dynkin types} of small rank, as initiated in \Cref{sec_pell_A3} and \Cref{sec_Pell_hyperoct},
\item considering \textit{more general weights} $\La$. 
Indeed, our first try (in type $C_2^{(1)}$ for $\La=\La_1$) was very successful, see \Cref{sec_Pell_C2},
\item using the \textit{extended} atomic length, which has already paid off in Sections \ref{sec_Pell_A2}, \ref{sec_Pell_C2}, \ref{sec_pell_A3}. 
\end{itemize}

\medskip

In turn, in order to get the most complete picture, it would be interesting to obtain combinatorial models (using partitions/abaci)
for (extended) affine Grassmannian elements in the favourable cases.
Moreover, one could imagine using the methods of \Cref{sec_pell_A3} (the case of $A_3^{(1)}$)
and invoke the universality of certain quadratic forms \cite{bhargava2000conway}, in order
to give an alternative, more constructive proof of the Granville-Ono theorem (in higher rank) and its analogues in other types.

\medskip

Finally, it would be natural to expect
results relating self-conjugate $6$-cores and $8$-cores (by considering types $C_3^{(1)}$ and $C_4^{(1)}$) and integral solutions
of appropriate quadratic equations (with $3$ and $4$ variables respectively).
In this context, some important results announced by Alpoge \cite[Theorem 4 and 6]{Alpoge2014} 
(and their consequences such as a proof of \cite[Conjecture 3.5]{HanusaNath2013}) 
contain some mistakes.
The case of $6$-cores was recently fixed by Hanson and Jameson \cite[Theorem 2]{HansonJameson2022}.
It would be interesting to fix the case of $8$-cores as well, 
and to see how this relates with our approach.
This is initiated in \Cref{sec_Pell_hyperoct}.

\renewcommand{\be}{q}

\section{Kac-Moody algebras, Weyl groups and other generalities}
\label{sec_KM}

\newcommand{\W}{W_{\mathrm{aff}}}

\subsection{Affine Kac-Moody algebras}
    
Let $A = (a_{ij})_{0\leq i, j\leq n}$
be a generalised (symmetrisable) Cartan matrix such that
the corresponding Kac-Moody algebra $\mathfrak{g}$ is of affine type \cite{Kac1984}. 
In particular,
the rank of $A$ is $n$.
Consider a realisation of $\mathfrak{g}$ (see \cite{Kac1984} for details),
in particular $\mathfrak{h}$ is a real vector space of dimension $n+2$ with dual space $\mathfrak{h}^\ast$
and we fix $\Delta=\{ \al_0,\ldots, \al_n\}\subseteq  \mathfrak{h}^\ast$ and
$\Delta^\vee=\{ \al_0^\vee,\ldots, \al_n^\vee\}\subseteq  \mathfrak{{h}}$ so that
$\langle \al_j, \al_i^\vee \rangle = a_{ij},$
where $\langle . , .\rangle$ denotes the natural pairing between $\mathfrak{h}^\ast$ and $\mathfrak{{h}}$. 
Moreover, to each type, one associates positive integers $a_0, \ldots, a_n$
and $a_0^\vee, \ldots, a_n^\vee$ 
given by the unique integer vector with no common factor in the kernel of $A$ and $A^\mathrm{t}$ respectively. 
Note that in untwisted types, 
$a_1, \ldots, a_n$ are the coefficients of the highest root in the basis of the simple roots $\al_1, \ldots , \al_n$ \cite[Definitions pp. 252, 485]{Carter2005}. 
In type $A_{2n}^{(2)}$, we have $a_0 = 2$, but otherwise we always have $a_0 = 1$. 
Note also that we always have $a_0^\vee=1$.

\medskip

There is a classification of these matrices into different Dynkin types,
which we record in \Cref{table}, alongside other important data whose definition we will recall below.
Note that transposing the Cartan matrix induces a natural duality at the level of Kac-Moody algebras,
and we denote $\mathfrak{g}^\vee$ the algebra associated to $A^\mathrm{t}$.
Its realisation is obtained by switching $\mathfrak{h}$ with $\mathfrak{h}^\ast$ and $\Delta$ with $\Delta^\vee$.
In \Cref{table}, the types with a superscript "{\footnotesize$(1)$}" are called \textit{untwisted},
and each of these has a dual \textit{twisted} type (with superscript "{\footnotesize$(2)$}" or  "{\footnotesize$(3)$}")
\footnote{Note that twisted type $A_{2n}^{(2)}$ of \Cref{table} is isomorphic to its dual.}. 
Untwisted types coincide with the classic geometrical construction of affine Weyl groups obtained for instance in \cite{BOURB}.

\input{table}

\subsection{Root space, dual root space and some important elements}

We introduce now two important spaces that we shall use all along the article
\begin{align}
\label{def:V0}
    V_0=\bigoplus_{i=1}^n \R\al_i \quad \text{and}\quad V = V_0\oplus \R\al_0. 
\end{align}

\begin{Rem}\label{remark conventions}
Beware that in \Cref{table},
we have used a renormalisation of Bourbaki's standard realisation
and of the corresponding lattices \cite{BOURB, Carter2005}.
This ensures that the norm of elements in $V_0$ (computed by the formula of \Cref{rk_norm_1}) 
coincides with the Euclidean norm.
\end{Rem}

We also introduce the following important elements.
The \textit{Coxeter number} is
\begin{align}\label{coxter number}
    h=\sum_{i=0}^n a_i.
\end{align}

We set
\begin{align}\label{elements delta and c}
    \delta = \sum_{i=0}^n a_i\al_i\in  V
\mand
c = \sum_{i=0}^n a_i^\vee \al_i^\vee \in \mathfrak{h}
\end{align}

and we define
\begin{align}\label{root theta}
    \theta = \delta - a_0\alpha_0 =  \sum_{i=1}^n a_i\alpha_i \in V_0.
\end{align}

\subsection{Weyl group}

For $0\leq i\leq n $, the reflections $s_i : \mathfrak{h}^\ast \to \mathfrak{h}^\ast$ given by the formula
\begin{align}\label{generators s_i}
    s_i(x) = x -\langle x, \al_i^\vee \rangle \al_i
\end{align}
generate a subgroup 
\begin{align}\label{def Weyl group from Kac}
    W = \langle s_0, s_1,\dots, s_n\rangle \leq \mathrm{GL}(\mathfrak{h}^\ast)
\end{align}
called the \emph{Weyl group}  of $\mathfrak{g}$. The subgroup of $W$ generated by $s_1,\dots, s_n$ is called the \emph{finite Weyl group} of $\mathfrak{g}$ and we denote it by $W_0$. In other words, $W_0$ is the Weyl group of the Cartan matrix $A_0$ obtained from $A$ by deleting the first row and the first column.
The groups $W$ and $W_0$ are Coxeter groups with sets of generators respectively $S = \{s_0, s_1, \ldots, s_n\}$ and $S_0 = \{ s_1, \ldots, s_n\}$.  
The group $W_0$ acts on $\Delta_0=\Delta\setminus\{\al_0\}$ 
and we denote $\Phi_0 = W_0(\Delta_0)$. The elements of $\Phi_0$ are called \textit{roots} and each root $\beta$ decomposes uniquely as $\beta = \sum_{i=1}^n b_i\alpha_i$ with $b_i \in \mathbb{Z}$. The \textit{height} of $\beta$ is the integer $\text{ht}(\beta) = \sum_{i=1}^n b_i$.

\medskip

Each root $\beta \in \Phi_0$ is such that there exists $\alpha_i \in \Delta_0$ and $w \in W_0$ with $w(\alpha_i) = \beta$. We define the dual root of $\beta$, denoted $\beta^{\vee}$, as $\beta^{\vee} := w(\al_i^{\vee})$. The dual root system of $\Phi_0$ is $\Phi_0^{\vee} := \{\beta^{\vee}~|~\beta \in \Phi_0\}$. It is known that if $\beta = \sum_{i=1}^n b_i\alpha_i \in \Phi_0$, then 
 $\beta^{\vee}$ decomposes as 
 $\beta^\vee = \sum_{i=1}^n \frac{b_i|\alpha_i|^2}{|\beta|^2}\alpha_i^{\vee}$ \cite[Section 3.14, ex. 3.2]{Kac1984}.
 
Formula \Cref{generators s_i} extends for any $\al \in \Phi_0$, and we have 
\begin{align}\label{reflexion s_alpha}
    s_{\alpha}(x) = x -\langle x, \al^\vee \rangle \al.
\end{align}

In the rest of the paper, the following classification of the different Kac-Moody types will be useful:
\begin{center}
\begin{tabular}{lll}
  \textbf{Type I}   & : & untwisted type or type $A_{2n}^{(2)}$,
  \\
 \textbf{Type II} & : &  twisted type other than $A_{2n}^{(2)}$.
\end{tabular}
\end{center}

The following is \cite[Proposition 17.18]{Carter2005}.

\begin{Prop}\label{prop_theta}
The element $\theta$ is
\begin{itemize} 
    \item the highest root of $\Phi_0$ in Type I,
    \item the highest short root of $\Phi_0$ in Type II.
\end{itemize}
\end{Prop}

Finally, we have the following important formula (\cite[$\S6.4$]{Kac1984})
\begin{align}\label{def alpha_0 vee}
    \al_0^{\vee} = c - a_0\theta^{\vee}.
\end{align}

\subsection{Weights}\label{section lattices}

\subsubsection{Affine dominant weights} 

Let us choose an element $\La_0^\vee\in\mathfrak{h}$ extending $\Delta^\vee$ into a basis of $\mathfrak{h}$
and such that $\langle \al_j, \Lambda_0^\vee\rangle=\delta_{0,j}$ for all $0\leq j \leq n$.
The \textit{fundamental weights} $\Lambda_i$, $0\leq i\leq n$ are the elements of $\mathfrak{h}^\ast$ determined by
\begin{align}
\left\langle \La_i, \al_j^\vee \right\rangle
=\delta_{ij} 
\text{ \quad for all } 0\leq i,j\leq n
\quad\mand\quad
\left\langle \La_i, \Lambda_0^\vee \right\rangle = 0
\text{ \quad for all } 0\leq i\leq n.
\end{align}
We have for $0\leq j\leq n$ 
\begin{align}
\left\langle \delta, \al_j^\vee \right\rangle 
=
\left\langle \sum_{k=0}^na_k\al_k , \al_j^\vee \right\rangle 
=
\sum_{k=0}^na_k \left\langle  \al_k , \al_j^\vee \right\rangle
=
\sum_{k=0}^n a_k a_{jk}
=
0,
\end{align}
where the last equality follows from the definition of $a_0,a_1,\ldots, a_n$,
and
\begin{align}
\left\langle \delta, \La_0^\vee \right\rangle 
=
\left\langle \sum_{k=0}^na_k\al_k , \La_0^\vee \right\rangle 
=\sum_{k=0}^na_k \left\langle  \al_k , \La_0^\vee \right\rangle 
= a_0,
\end{align}
therefore $\{  \La_0,\La_1,\ldots,\La_n, a_0^{-1}\delta \}$
is the dual basis of  $\{  \al_0^\vee,\al_1^\vee,\ldots,\al_n^\vee, \La_0^\vee \}$ in $\mathfrak{h}$.

\medskip

The \textit{affine weight lattice} is defined by
\begin{align}
 P  := \{\La \in \mathfrak{h}^*~|~\langle \La, \al_i^{\vee} \rangle \in \mathbb{Z}~\forall i=0,\dots, n ~ \text{ and } ~ \langle \La, \Lambda_0^{\vee} \rangle \in \mathbb{Z}\} 
 =\bigoplus\limits_{i=0}^{n} \mathbb{Z}\La_i \oplus \frac{1}{a_0}\mathbb{Z}\delta,
\end{align}
and the set of affine \textit{dominant} weights is given by 
\begin{align}
    P^+ = \bigoplus\limits_{i=0}^{n} \N\La_i \oplus \frac{1}{a_0}\mathbb{Z}\delta.
\end{align}

\medskip

We also have the set of \textit{affine fundamental coweights} $\{ \La_0^\vee, \La_1^\vee,\ldots, \La_n^\vee \}$
which verify the dual relations. 
Finally we let $\rho^\vee$ be
\begin{equation}\label{def rho vee}
    \rho^\vee = \sum_{i=0}^{n} \La_i^\vee.
\end{equation}

\subsubsection{Finite dominant weights}

We set, for all $1\leq i\leq n$,
\begin{align}
\label{omega_i_def}
\om_i = \La_i - \frac{a_i^\vee}{a_0^\vee} \La_0.
\end{align}
Then one can check that
\begin{align}\label{omega_rels}
\left\langle \om_i, \al_j^\vee \right\rangle
=\delta_{ij} \text{ \quad for all } 1\leq i,j\leq n
\mand 
\omega_i\in V_0 \quad\text{ for all } 1\leq i\leq n.
\end{align}
Therefore, these elements are called the \textit{finite fundamental weights},
as they are the fundamental weights of the corresponding finite type Kac-Moody algebra
with Cartan matrix $A_0=(a_{ij})_{1\leq i,j\leq n}$.

\medskip 

The (finite) \textit{weight lattice} is the dual lattice (for the pairing $\langle~,~\rangle$) of the lattice $\bigoplus_{i=1}^n\mathbb{Z}\alpha_i^{\vee}$, that is  
\begin{align}
    P_0 := \{x \in V_0~|~\langle f,x \rangle \in \mathbb  {Z} ~\text{~for any}~ f \in \bigoplus_{i=1}^n\mathbb{Z}\alpha_i^{\vee}\} = \bigoplus_{i=1}^n \mathbb{Z}\omega_i.
\end{align}
The set $P_0^+= \bigoplus_{i=1}^n \N \om_i$ is called the set of \textit{integral dominant weights} associated to $P_0$.
There is a similar notion of \textit{fundamental coweights}, defined by
\begin{align}
\om_i^\vee = \La_i^\vee - \frac{a_i}{a_0} \La_0^\vee
\end{align}
and which give the dual basis of $\Delta_0$.
Similarly, we  consider the \textit{coweight lattice}
  \begin{equation}\label{lattice P0vee}
 P_0^\vee := \bigoplus_{i=1}^n \mathbb{Z}\omega_i^{\vee}.
  \end{equation}

\subsection{Bilinear form on $\mathfrak{h}^*$}

The space $\mathfrak{h}^*$ decomposes as follows (\cite[Section 6.2.5]{Kac1984})
\begin{align}\label{decompo h etoile}
    \mathfrak{h}^* = V_0 \oplus \mathbb{R}\delta \oplus \mathbb{R}\Lambda_0, 
\end{align}

and is equipped with a non-degenerate symmetric bilinear form $(-~|~-)$ satisfying (\cite[(6.2.2), (6.2.4)]{Kac1984})

\begin{equation}\label{list formula inner product}
\left\{
\begin{aligned}
&\left( \alpha_i \mid \alpha_j  \right) = a_i^\vee a_i^{-1} a_{ij}, \quad \text{for} \quad i,j = 0,1,\dots, n  \\
&\left( \delta \mid \alpha_i  \right) = 0, \quad \text{for} \quad i = 0,1,\dots, n  \\
&\left( \delta \mid \delta  \right) = 0  \\
&(\alpha_i~|~\Lambda_0) = 0  \quad \text{for} \quad i = 1,\dots, n \quad \text{and}~(\alpha_0~|~\Lambda_0) = a_0^{-1}  \\
&(\delta~|~\Lambda_0) = 1 \\
&(\Lambda_0~|~\Lambda_0) = 0 .
\end{aligned}
\right.
\end{equation}
 
The bilinear form $(-~|~-)$ restricted to $V_0$ is an inner product.

\begin{Rem}\label{rk_norm_1}
Note that the formulas above completely determine the values of $(\al_i\mid\al_j)$.
In particular, the squared norm of a simple root is given by 
\begin{align}\label{inner product al_i with al_i}
    |\al_i|^2 =(\al_i~|~\al_i) = 2\frac{a_i^\vee}{a_i}.
\end{align}
\end{Rem}

The root $\theta$ defined in \Cref{root theta} satisfies the following  equality (\cite[Section 6.4]{Kac1984})
 \begin{align}\label{equality norm theta}
     (\theta~|~\theta) = 2a_0,
 \end{align}

and for any $\al \in \Phi_0$, one has
\begin{align}\label{inner product theta with simple root}
    (\theta~|~\al) \in \mathbb{N} ~\text{~if~}~\al \in \Phi_0^+ \quad \text{~and~} \quad (\theta~|~\al) \in -\mathbb{N} ~\text{~if~}~\al \in \Phi_0^-.
\end{align}
To see \Cref{inner product theta with simple root}, it suffices to do it for the simple roots $\al_j, j=1,\dots, n$:
\begin{align}
    (\theta~|~\al_j) = (\delta-a_0\al_0~|~\al_j) \overset{\Cref{list formula inner product}}{=} -a_0(\al_0~|~\al_j) \overset{\Cref{list formula inner product}}{=} -a_0^{\vee}a_{0j} \overset{a_0^\vee=1}{=} -a_{0j} \in \mathbb{N}.
\end{align}

It turns out that the form $(-~|~-)$ is $W$-invariant, that is
\begin{align}\label{invariance form ()}
    \big(w(x)~|~w(y)\big) = \big(x~|~y\big) ~\quad \forall w \in W,~\forall x,y \in \mathfrak{h}^*.
\end{align}

\subsection{The map $\zeta$ and some important formulas}

Let $\zeta: \mathfrak{h}^* \to \mathfrak{h}$ be the isomorphism defined by the condition 
$\langle x, \zeta(y) \rangle = (x~|~y)$ for any $x,y \in \mathfrak{h}^*$. 
For any $i=0,1,\dots, n$, we have (\cite[Formula (6.2.3)]{Kac1984})
\begin{align}\label{formula Kac 6.2.3}
    \zeta(\al_i) = \frac{a_i^{\vee}}{a_i}\al_i^{\vee}= \frac{(\al_i\mid\al_i)}{2}\al_i^{\vee} \quad \text{and} \quad  \zeta(\delta) = c \quad \text{and} \quad \zeta(\theta) = a_0\theta^{\vee}.
\end{align}

In the following, it will be convenient to 
replace $\Phi^\vee$ by $\Phi^\dagger := \zeta^{-1} (\Phi^\vee) \subset \mathfrak{h}^\ast$,
and to transport the properties of $\langle-,-\rangle$ to $(-\mid-)$.
We denote
\begin{equation}\label{def root dagger}
    \al_i^{\dagger} := \zeta^{-1}(\al_i^\vee) =  \frac{2}{(\al_i~|~\al_i)}\al_i.
\end{equation}

By \Cref{formula Kac 6.2.3}, it follows for any $y \in \mathfrak{h}^*$ and any $i=0,1,\dots, n$, that
\begin{align}\label{connection inner product et pairing}
    \langle y, \al_i^\vee \rangle = \frac{a_i}{a_i^\vee}(y~|~\al_i) = (y\mid \al_i^\dagger)
    \quad \text{and} \quad \langle \al_i, c \rangle = 0.
\end{align}
It also follows from \Cref{formula Kac 6.2.3}, for any $y \in \mathfrak{h}^*$, that
\begin{align}\label{formula theta vee duality}
    \langle y, \theta^{\vee} \rangle = a_0^{-1}\langle y, \zeta(\theta) \rangle = a_0^{-1} \left(y~|~\theta \right).
\end{align}

Since $a_0^\vee=1$, from \Cref{list formula inner product} and \Cref{connection inner product et pairing} we deduce that for any $y \in V$
\begin{align}\label{crochet c avec h etoile}
    \langle y, c \rangle = 0 \quad \text{and} \quad \langle \Lambda_0, c  \rangle = 1.
\end{align}

We also deduce that
\begin{align}\label{action W sur delta et action W_0 Lambda_0}
    w(\delta) = \delta ~~\forall w \in W \quad \text{~and~} \quad w(\Lambda_0) = \Lambda_0 ~~\forall w \in W_0 \quad \text{~and~} \quad w(\Lambda_0) \in \Lambda_0 + V~~\forall w \in W.
\end{align}
We obtain the following useful formula.
\begin{align}\label{formula utile crochet c}
    \langle w(y), c  \rangle = \langle y, c\rangle \quad \forall y \in \mathfrak{h}^*,~\forall w \in W.
\end{align}

Now, we have
$$
\delta_{ij}
= \langle \al_j, \omega_i^{\vee} \rangle 
= \langle \omega_i, \al_j^{\vee}  \rangle 
= \frac{a_j}{a_j^{\vee}}(\omega_i~|~\al_j)
= \frac{a_j}{a_j^{\vee}}(\alpha_j~|~\omega_i) 
= \frac{a_j}{a_j^{\vee}}\langle \al_j, \zeta(\omega_i) \rangle.
$$
This shows that for all $1\leq i\leq n$,
\begin{align}\label{zeta on omega_i}
    \zeta(\omega_i) = \frac{a_i^{\vee}}{a_i}\omega_i^{\vee},
\end{align}
which shall be compared to \Cref{formula Kac 6.2.3}.
Finally, for $i=1,\dots, n$, we define the element $\varpi_i \in V_0$ by
 \begin{align}\label{varpi in terms of omega}
  \varpi_i := \zeta^{-1}(\omega_i^{\vee}) \overset{\Cref{zeta on omega_i}}{=} \frac{a_i}{a_i^{\vee}}\omega_i.
 \end{align}

Therefore, for any $j = 0,1,\dots,n$, we have
\begin{align}
    \langle \varpi_i, \al_j^{\vee} \rangle = \frac{a_i}{a_i^{\vee}}\delta_{ij}.
\end{align}

Finally, we have two dual basis relationships with respect to $(-\mid -)$, namely
\begin{align}\label{dual_bases}
    (\varpi_i\mid \al_j) =\delta_{ij}
    =
    (\om_i\mid \al_j^\dagger).
\end{align}

\subsection{The lattice ${{M}}$ and the translations}

Let $\frac{1}{a_0}W_0 \cdot \theta$ be the $W_0$-orbit of the element $\frac{1}{a_0}\theta$. 
Following \cite[Section 17.3]{Carter2005}\footnote{In \cite{Carter2005}, the notation used is $M^\ast$. We choose ${{M}}$ instead.}, we define now the lattice ${{M}} \subset V_0$ by
\begin{align}\label{orbit theta}
{{M}} := \frac{1}{a_0}\mathbb{Z}W_0\cdot \theta.
\end{align} 

\begin{Rem}\label{rem_duality_trick}
Let us recall some key properties of the lattice introduced in~\Cref{orbit theta}, 
as they will play a central role in the sequel.

\begin{enumerate}[(i)]
    \item Suppose We are in untwisted type. Then $\theta$ is the highest root
    of $\Phi_0$, and $\theta^\dag = \theta$. Since $\theta^\dag$ is a
    short root of $\Phi_0^\dag$, we obtain
    $$
    M = \mathbb{Z} W_0 \cdot \theta = \mathbb{Z} W_0 \cdot \theta^\dag = \mathbb{Z} \Phi_0^\dag.
    $$
   Note that $\theta$ is the highest root of the dual root system of
   $\Phi_0^\dag$, since $(\Phi_0^\dag)^\dag = \Phi_0$.

    \item Suppose that $A$ is of twisted type $A_{2n}^{(2)}$. Then
    $\theta$ is again the highest root of $\Phi_0$, and, by our choice of
    normalization, we have
    $$
    \theta^\dag = \frac{1}{2} \theta.
    $$
    Thus,
    $$
    M = \mathbb{Z} W_0 \cdot \theta^\dag = \mathbb{Z} \Phi_0^\dag,
    $$
    since $\theta^\dag$ is a short root of $\Phi_0^\dag$. As in the
    untwisted case, $\theta$ is the highest root of the dual root system
    of $\Phi_0^\dag$.

    \item Suppose we are in twisted type. Then $\theta$ is the highest short root of $\Phi_0$, hence $M = \mathbb{Z} \Phi_0$, and $\theta^\dag$ is
    the highest root of $\Phi_0^\dag$. With our choice of
    normalization, we observe that
    $$
    \theta^\dag = \zeta^{-1}(\theta^\vee),
    $$
    so that the coefficients of $\theta^\dag$ in the basis $\Delta^\dag$
    are precisely the $a_i^\vee$ for $1 \leq i \leq n$.
\end{enumerate}
\end{Rem}

A concrete description of ${{M}}$ in terms of the simple roots can also
be found in \cite[Proposition 17.23 and Remark 17.34]{Carter2005},
where it is shown that

\begin{equation}
\label{lattice M dagger explicit}
\begin{array}{ccccc}
{{M}} 
& = 
& 
\left\{
\begin{array}{ll}\ds
\bigoplus_{i=1}^n\mathbb{Z}\frac{a_i}{a_i^{\vee}}\alpha_i 
&\text{\quad  in Type I}\\
\ds\bigoplus_{\alpha_i~\text{short}}\mathbb{Z}\frac{a_i}{a_i^{\vee}}\alpha_i \oplus \bigoplus_{\alpha_i~\text{long}} \mathbb{Z}\alpha_i 
&\text{\quad in Type II}
\end{array}
\right.
& = 
&
\left\{
\begin{array}{ll}
\ds
\bigoplus_{i=1}^n\mathbb{Z}\alpha_i^{\dagger} 
& \text{\quad  in Type I} \\
\ds \bigoplus_{i=1}^n\mathbb{Z}\alpha_i 
&\text{\quad in Type II}.
\end{array}
\right.
\end{array}
\end{equation}

Accordingly, we define ${{{M^\dagger}}}$ by

\begin{equation}
\label{lattice M new gen}
\begin{array}{ccc}
{{{M^\dagger}}} 
&
= 
&
\left\{
\begin{array}{ll}
\ds \bigoplus_{i=1}^n\mathbb{Z}\alpha_i 
&\text{\quad  in Type I} \\
 \ds \bigoplus_{i=1}^n\mathbb{Z}\alpha_i^{\dagger} 
 &\text{\quad in Type II}.
\end{array}
\right.
\end{array}
\end{equation}

\medskip

We refer to \Cref{table} for this explicit description of ${{M}}$ in terms of the simple roots,
and we also give its realisation in terms of the standard basis of the ambient space. 

\medskip

For any element $q \in V$ we define the linear map $t_q \in \text{Hom}(\mathfrak{h}^*)$,
called the \textit{translation} by $q$ (with the terminology justified in \Cref{group W_a}), by the formula 
\begin{equation}
\label{translation}
t_q(v) := v + \langle v,c\rangle q - \left( (v~|~q)+\frac{1}{2}|q|^2\langle v,c\rangle \right)\delta.
\end{equation}

An easy computation gives 
    $s_0 = t_{\frac{1}{a_0}\theta}s_{\theta}$,
and using \Cref{invariance form ()} and \Cref{formula utile crochet c}, for any $w \in W$ and any $q,r \in V$, we have that
\begin{align}\label{formula conjugate t_q}
    t_qt_r = t_{q+r} \quad \text{and} \quad wt_qw^{-1} = t_{w(q)}.
\end{align}

Let $T({{M}})$ be the group of translations $t_q$ with $q \in {{M}}$. We have that (see \cite[Prop. 6.5]{Kac1984})
\begin{align}\label{chara affine weyl group via semi direct}
    W = T({{M}})\rtimes W_0.
\end{align}

Therefore any element $w \in W$ decomposes uniquely as $w = t_q\overline{w}
$ where $q \in {{M}}$ and $\overline{w} \in W_0$.

\subsection{Alcove geometry for Weyl groups of affine Kac-Moody algebras}
\label{sec_alcoves}

In the non-twisted cases, the concept of alcoves is well-understood and follows easily from the definition of the affine Weyl group, in the sense of Bourbaki \cite{BOURB} or Humphreys \cite{humphreys1992reflection}. 
However, it is also possible to define alcoves for any Weyl groups associated to affine Kac-Moody algebras. 
In this section, we explain it by following \cite[Section 6.6]{Kac1984}, while providing additional precision.
The Weyl group $W$ is a subgroup of $\text{GL}(\mathfrak{h}^*)$. However, for our purpose, the space $\mathfrak{h}^*$ is unnecessarily large, and we can instead view $W$ as a subgroup of $\text{Aff}(V_0)$, the space of affine transformations of $V_0$. The alcove geometry of $W$ appears from this consideration in $V_0$.

\subsubsection{The group of affine transformations}\label{group W_a}

Let $k \in \mathbb{R}$. We define the following hyperplane of $\mathfrak{h}^*$
\begin{align}\label{hyperplane in V via c}
    \mathcal{H}_k := \{v \in \mathfrak{h}^*~|~\langle v, c\rangle = k\} 
     = \{v+k\Lambda_0,~v \in V\} = \{x + a\delta+k\Lambda_0, x \in V_0, a \in \mathbb{R}\}.
\end{align}

The second equality of \Cref{hyperplane in V via c} follows easily  from \Cref{decompo h etoile} and \Cref{crochet c avec h etoile}. These hyperplanes possess an important property: they are invariant under the action of $T(V_0)\rtimes W$. To see this, it suffices to consider the action of the simple generators $s_i$ for $i=0,1,\dots, n$ along with \cref{crochet c avec h etoile}, and to see that $T(V_0)$ acts on $\mathcal{H}_k$ as well.

\medskip

Now, from \Cref{generators s_i} we do not have $W$ that acts on $V_0$, because of $s_0$. Consequently, $T(V_0)\rtimes W$ does not act on $V_0$ either. It is however possible to find a space that is 
in bijection with $V_0$ and on which $T(V_0)\rtimes W$ has a well-defined action.
Since $T(V_0)\rtimes W$ acts on $\mathcal{H}_1$ and stabilises $\mathbb{R}\delta$ \Cref{action W sur delta et action W_0 Lambda_0}, it follows that $T(V_0)\rtimes W$ also acts on the affine space $\mathcal{H}_1/\mathbb{R}\delta = \{ [x+\La_0]  \,;\, x\in V_0 \}$
by $g \cdot [v] := [g(v)]$.
We consider now the structural map of this action, defined by
$\phi(g) : [v] \mapsto [g(v)] $.
It is easy to see that $\phi_{|W} : W\to \mathrm{Sym}(\cH_1/\R\delta)$ is an injective morphism. 
Now, the following map is clearly a bijection
\begin{equation}\label{iso invariant W-space}
\begin{array}{ccccc}
     \Theta & : & \cH_1/\R\delta & \longrightarrow & V_0  \\
            &   & \left[ x + a\delta + \Lambda_0 \right]& \longmapsto &  x.
\end{array}
\end{equation}

We define the morphism $\text{af} : T(V_0) \rtimes W \rightarrow \text{Aff}(V_0)$ 
making the following diagram commutes  $\forall g \in T(V_0) \rtimes W$.

\begin{center}
\begin{tikzpicture} \label{def Qrl}
\node at (0,0) {$\cH_1/\R\delta$} ;
\node at (3,0) {$\cH_1/\R\delta$} ;
\node at (3,-2) {$V_0$} ;
\node at (0,-2) {$V_0$} ;

\node at (0.5,0) (1)  {} ;
\node at (2.5, 0) (2) {} ; 
\draw [>=stealth,->] (1) to (2);
\node at (1.5,0.3) {$\phi(g)$};

\node at (0,-0.3) (3)  {} ;
\node at (0, -1.8) (4) {} ; 
\draw [>=stealth,->] (3) to (4);
\node at (-0.4,-1) {$\Theta$};

\node at (3,-0.3) (5)  {} ;
\node at (3, -1.8) (6) {} ; 
\draw [>=stealth,->] (5) to (6);
\node at (3.5,-1) {$\Theta$};

\node at (0.5,-2) (7)  {} ;
\node at (2.5, -2) (8) {} ; 
\draw [>=stealth,->] (7) to (8);
\node at (1.5,-2.3) {$\text{af}(g)$};

\end{tikzpicture}
\end{center}

For any $g \in W$ and any $v = x + a\delta + \Lambda_0 \in \mathcal{H}_1$, we have $\text{af}(g)\circ \Theta([v]) = \text{af}(g)(x) = \Theta \circ \phi(g)([v]) = \Theta([g(v)]) = \Theta([g(x) + ag(\delta) + g(\Lambda_0)]) = \Theta([g(x) + g(\Lambda_0)])$, where the last equality comes from \Cref{action W sur delta et action W_0 Lambda_0}. Therefore, $T(V_0) \rtimes W$ acts on $V_0$ by
\begin{equation}\label{action W sur V_0}
    \text{af}(g)(x) = \Theta([g(x) + g(\Lambda_0)]).
\end{equation}
We now analyse\Cref{action W sur V_0} in detail 
restricted to three important subsets of $W$, namely $W_0$, $T(V_0)$, and $\{s_0\}$. 
\begin{itemize}
    \item[$\bullet$] Let $w \in W_0$. Then, for all $x \in V_0$,
 we have 
    \begin{align*}
        \text{af}(w)(x) & = \Theta([w(x) + \Lambda_0]) & \text{by}~\Cref{action W sur delta et action W_0 Lambda_0}\\
        &= w(x). &\text{since} ~ w(V_0) \subset V_0
    \end{align*}

\item[$\bullet$]  Let $q\in V_0$. Then, for all $x \in V_0$,
 we have 
\begin{align*}
    \text{af}(t_q)(x)
    & = \Theta([t_q(x) + t_q(\Lambda_0)]) & \text{by}~\Cref{action W sur V_0}\\
                      & = \Theta([x + \langle x,c\rangle q  + \Lambda_0 + \langle \Lambda_0,c\rangle q]) & \text{by}~ \Cref{translation}~\text{and}~\Cref{iso invariant W-space}\\
                      & = \Theta([x + \Lambda_0 + q]) & \text{by}~\Cref{crochet c avec h etoile}  \\
                      & = x + q. & \text{by}~\Cref{iso invariant W-space}
\end{align*}
    
    \item[$\bullet$] For all $x \in V_0$,
 we have 
    \begin{align*}
    \text{af}(s_0)(x) 
                        & = \Theta([s_0(x) + s_0(\Lambda_0)]) & \text{by}~\Cref{action W sur V_0} \\
                        & = \Theta([x - \langle \al_0^{\vee},x\rangle \al_0 + \Lambda_0 - \langle \al_0^{\vee},\Lambda_0\rangle \al_0]) & \text{by}~\Cref{generators s_i}\\
                        & = \Theta([x - \Big(\langle \al_0^{\vee},x\rangle +1\Big)\al_0 + \Lambda_0]) & \text{by} ~ \Cref{crochet c avec h etoile}\\ 
                        & = \Theta([x - \Big(\langle c - a_0\theta^{\vee},x\rangle +1\Big)\al_0 + \Lambda_0]) & \text{by}~\Cref{def alpha_0 vee}\\
                        & = \Theta([x - \Big(-a_0\langle \theta^{\vee},x\rangle +1\Big)\al_0 + \Lambda_0]) & \text{by}~\Cref{crochet c avec h etoile}\\ 
                        & = \Theta([x  + \langle \theta^{\vee},x\rangle (\delta - \theta) - a_0^{-1}(\delta-\theta) + \Lambda_0]) & \text{by}~\Cref{root theta}\\
                        & =\Theta([x + \Big( a_0^{-1} -  \langle \theta^{\vee},x\rangle \Big) \theta  + \Big(\langle \theta^{\vee},x\rangle - a_0^{-1}\Big)\delta  + \Lambda_0]) & \\
                        & = x + a_0^{-1}\Big( 1 -  (\theta~|~x)\Big)\theta. & \text{by}~\Cref{iso invariant W-space}, \Cref{formula theta vee duality}
\end{align*}
\end{itemize} 

Since $\phi_{|W}$ is an injective morphism, $\text{af}_{|W}$ is also an injective morphism and we denote by $W_a$ its image, that is
\begin{align}\label{iso W Wa}
    W_a := \langle \text{af}(s_i)\, ; \,i =0,1,\dots,n \rangle \simeq W,
\end{align}

By \Cref{chara affine weyl group via semi direct},  we have 
$
W_a = \text{af}(T({{M}})) \rtimes \text{af}(W_0).
$
Note that in untwisted types, the group $W_a$ is usually given as the definition of the affine Weyl group (see for example \cite{BOURB, Shi87}).

\subsubsection{Alcoves}

Let $\alpha \in \Phi_0^+$ and $k \in \mathbb{Z}$. We define the affine hyperplane $H_{\al,k}$ of $V_0$ as follows

\begin{equation}\label{def hyperplans}
H_{\al,k} :=  \left\{
\begin{aligned}
& \left\{ x \in V_0~|~(\alpha~|~x) = k \right\} \text{ \quad  in Type I} \\
& \left\{ x \in V_0~|~(\alpha^{\dagger}~|~x) = k \right\} \quad \text{in Type II}.\\
\end{aligned}
\right.
\end{equation}

\begin{Rem}
The hyperplane $H_{\al,k}$ is defined by means of the scalar product instead of the duality, as was done in  \cite[Remark 17.34]{Carter2005}.
However, our collection of hyperplanes correspond exactly to that of \cite{Carter2005}, and therefore define the same alcove geometry (see below).
\end{Rem}

The connected components of 
\begin{equation}
 V_0 ~\backslash \bigcup\limits_{\begin{subarray}{c}
 ~ ~\alpha \in \Phi_0^{+} \\ 
  k \in \mathbb{Z}
\end{subarray}}
H_{\alpha,k} 
\end{equation}
are simplices called \emph{alcoves}. We denote $\mathcal{A}_e$ the alcove
\begin{equation}
  \mathcal{A}_e = \left\{ x \in V_0~|~ (\al_i~|~x) > 0~i=1,\dots,n ~\text{and}~(\theta~|~x) < 1 \right\}.
\end{equation}

The alcove $\mathcal{A}_e$ is called the \textit{fundamental alcove}.
The vertices of $\mathcal{A}_e$ are given by
\begin{equation}\label{vertices fundamental alcove}
\{0\} \sqcup \left\{\frac{\varpi_i}{a_i}\, ; \, i = 1,2,\dots, n\right\}.
\end{equation}

The group $W$ acts on the set of alcoves as follows
\begin{align}\label{def action W sur set alcoves}
    w\mathcal{A} := \text{af}(w)(\mathcal{A}).
\end{align}

It turns out that this action is regular, giving therefore a bijection between $W$ and the set of alcoves
$w \mapsto \cA_w,$
where $\cA_w:= w\cA_e$.

\begin{Rem}\label{rem_wall_crossing}
We can color the faces of the alcoves, called the \textit{walls}, as follows.
First, color each wall of the fundamental alcove by $i\in\{1,\ldots,n\}$
if it lies on $H_{\al_i,0}$ and by $0$ if it lies on $H_{\theta,1}$.
Then, color all remaining walls by recursively reflecting the colored walls across all adjacent hyperplanes.
It is clear that each alcove as a unique wall of color $i$, called an \textit{$i$-wall}.
Moreover,
for $w\in W$ and $0\leq i \leq n$,
$(ws_i)\cA_e$ is the alcove obtained by reflecting the alcove $w\cA_e$ along its $i$-wall.
\end{Rem}

Finally, we will need the \textit{fundamental chamber} $C_0$, which is defined by
\begin{equation}\label{fundamental chamber}
    C_0 = \{x \in V_0~|~(x~|~\alpha_i) > 0~\forall i=1,\dots,n\}.
\end{equation}

\subsection{Affine Grassmannian}
\label{sec_AG}

In this section we recall the notion of \textit{affine Grassmannian}, 
a well-studied subset of the affine Weyl group which turns out to be naturally in bijection with the $(n+1)$-cores in type $A_n^{(1)}$ (see \Cref{sec_cores}). Therefore, for an arbitrary affine Weyl group, affine Grassmannian elements
give a natural generalisation of cores, and we will therefore sometimes use the terminology \textit{generalised cores}.

\medskip

Let us start by some generalities on transverse sets.
Let $(W,S)$ be a Coxeter system with length function $\ell$. Let $\Phi = \Phi^+ \sqcup \Phi^-$ be the root system of $W$ with simple system $\Delta$. Let $I,J \subset S$ and let $(W_I,S_I)$ and $(W_J,S_J)$ be the corresponding parabolic subgroups of $W$. Let $T$ be the set of reflections of $W$, that is $T = \{wsw^{-1}~|~s \in S, w \in W\}$. 

\medskip

The (left-)inversion set an element $w \in W$ is defined  by 
\begin{align}\label{inversion sets coxeter}
    N(w) = \{t \in T ~|~\ell(tw) < \ell(w)\}.
\end{align}

The set $W^J$ is called the $J$-right transverse set of $W$ and is defined by
\begin{align}\label{def transverse set}
     W^J = \{w \in W~|~\ell(ws) > \ell(w)~\forall s \in J\}.
\end{align}

Each element $w \in W$ has unique decomposition  $w = w^Jw_J{}$ with  $w^J \in W^J$ and $w_J \in W_J$. This decomposition satisfies $\ell(w) = \ell(w^J) + \ell(w_J)$.  The element $w^J$ is the unique element of minimal length in $wW_J$. Therefore the following map is a natural bijection

\begin{equation}\label{eq:array2}
\begin{array}{ccccc}
\chi_{\raisebox{-0.5ex}{$\scriptstyle J$}}& : &\faktor{W}{W_J} & \longrightarrow &  W^J \\
     &   &\left[w\right]       & \longmapsto      & w^J.
\end{array}
\end{equation}

Let $W$ be an affine Weyl group with $S = \{s_0,s_1,\dots, s_n\}$ and let $W_0$ be its finite part with $I_0 =  \{s_1,\dots, s_n\}$ so that $W_0 = W_{I_0}$.
The \textit{(right) affine Grassmannian} is by definition the set 
$$
\faktor{W}{W_{0}}.
$$
We also have the analogous notion of \textit{left affine Grassmannian}, and a similar notion 
for any $I \subset S$ instead of $I_0$.
Using the map \Cref{eq:array2}, the affine Grassmannian is in bijection with 
 \begin{equation}
 \label{eq:w0}
W^0:= W^{I_0} = \{w \in W~|~\ell(ws_k) > \ell(w)~\forall k=1,\dots,n\}.
 \end{equation}
Note that $W^0$ is the set of elements $w\in W$ such that every reduced expression of $w$ ends with $s_0$.
Elements of $W^0$ are called \textit{affine Grassmanian elements}. We will use the map \Cref{eq:array2} in the situation $J = I_0$, and simply denote it by $\chi$, that is 
\begin{equation}\label{def chi}
\begin{array}{ccccc}
\chi & : &\faktor{W}{W_0} & \longrightarrow &  W^0 \\
     &   &\left[w\right]       & \longmapsto      & w^0.
\end{array}
\end{equation}

\medskip

Recall that we have defined the fundamental chamber $C_0$ in \Cref{sec_alcoves}.
It turns out that the bijection $w\mapsto \cA_w$ between $W$ and the alcoves restricts
to a bijection between $\{w\in W \mid w^{-1}\in W^0\}$ and the alcoves of the fundamental chamber, in other terms
\begin{align}\label{bij_alcoves_ag}
w^{-1}\in W^0  \quad \Longleftrightarrow\quad  \cA_w\subset C_0.
\end{align}

We now want to express the bijection \Cref{def chi} in terms of the lattice ${{M}}$. To streamline later arguments (see \Cref{rem_AL_AG}), we introduce the following straightforward bijection
\begin{equation}\label{map_f}
\begin{array}{ccccc}
f & : & \faktor{W}{W_0}  & \longrightarrow & {{M}}      \\
  &    & \left[t_q\overline{w}\right]       & \longmapsto       & q
\end{array}
\end{equation}
whose inverse map is given by $q \mapsto [t_q]$. In particular this shows that $ {{M}} \simeq W^0$ via $q \mapsto (t_q)^0$.

\begin{Rem}\label{rem_AG}
For an element $w = t_q\overline{w}\in W$ we set $\mathfrak{r}(w) := q$. Let $w = w^0w_0$ be the right $I_0$-decomposition of $w \in W$.
It is easy to see that for $q\in {{M}}$, the decomposition $t_q=(t_q)^0(t_q)_0$ implies that $\mathfrak{r}((t_q)^0) = q$. Therefore, for any $w = t_q\overline{w} \in W$ we also have $f([w]) = \mathfrak{r}(w^0)$.
\end{Rem}

\section{Extended affine Weyl group and extended affine Grassmannians}

\label{sec_ext}

\subsection{The lattice $L$ and the fundamental group $F$}
\label{sec_F} 

 We define the lattice $L$ as the dual lattice of ${{{M^\dagger}}}$ (defined in \Cref{lattice M new gen}), that is
\begin{equation}\label{lattice L}
      L = \{v \in V_0~|~(v~|~u) \in \mathbb{Z} ~\text{~for any~} ~ u \in {{{M^\dagger}}} \} = \left\{
\begin{aligned}
& \bigoplus_{i=1}^n\mathbb{Z}\varpi_i \text{ \quad  in Type I} \\
&  \bigoplus_{i=1}^n\mathbb{Z}\omega_i \quad \text{in Type II}.\\
\end{aligned}
\right.
 \end{equation}

\begin{Prop}\label{MvsL}
    We have ${{M}} \subset L$ and these two lattices are $W_0$-equivariant.
\end{Prop}

\begin{proof}
By \Cref{lattice M dagger explicit} we have two cases. 
\begin{itemize}
    \item In Type I, it suffices to show that $\al_i^{\dagger} \in L$ for any $i =1,\dots, n$, that is that $(\al_i^{\dagger}~|~\al_j) \in \mathbb{Z}$ for any $j = 1,\dots, n$. This follows from \Cref{inner product al_i with al_i}.
    \item In Type II, we show that $\al_i\in L$ for any $i =1,\dots, n$, that is that $(\al_i~|~\al_j^{\dagger}) \in \mathbb{Z}$ for any $j = 1,\dots, n$. As above, this follows from \Cref{inner product al_i with al_i}.
\end{itemize} 
The fact that $W_0({{M}}) = {{M}}$ follows directly from \Cref{formula conjugate t_q}, but can also be easily checked from \Cref{generators s_i} together with \Cref{lattice M dagger explicit}. Let us show now that $W_0(L) = L$. It suffices to show it for the simple reflections and for the generators of $L$. Once again we have two cases to consider. 
\begin{itemize}
    \item In Type I, we have
    \begin{align*}
        s_j(\varpi_i) = \varpi_i - \frac{a_i}{a_i^{\vee}}\delta_{ij}\al_j  = \left\{
    \begin{array}{ll}
    \varpi_j &\text{if} ~ j = i \\
    \varpi_j - \al_j^{\dagger}  & \text{otherwise}. 
    \end{array}
    \right.
    \end{align*}
    In this situation we have ${{M}} = \bigoplus_{i=1}^n\mathbb{Z}\alpha_i^{\dagger}$. Therefore, since ${{M}} \subset L$, we have that $\al_j^{\dagger} \in L$ and then in both cases  $s_j(\varpi_i) \in L$.
    
    \item In Type II, we have
    \begin{align*}
        s_j(\omega_i) = \omega_i - \delta_{ij}\al_j  = \left\{
    \begin{array}{ll}
    \varpi_j &\text{if} ~ j = i \\
    \varpi_j - \al_j  & \text{otherwise}. 
    \end{array}
    \right.
    \end{align*}
    In this situation we have ${{M}} = \bigoplus_{i=1}^n\mathbb{Z}\alpha_i$. Therefore, since ${{M}} \subset L$, we have that $\al_j \in L$ and then in both cases  $s_j(\omega_i) \in L$.
    \hfill
\end{itemize}
\end{proof}

The \textit{fundamental group}\footnote{The terminology is justified as follows: Let $G$ be a semi-simple simply connected algebraic group with Weyl group $W$. Then  $L/{{M}} \simeq \pi_1\big(G/Z(G)\big)$ where $Z(G)$ is the center of $G$.}, denoted $F$, is the group defined by
\begin{align}\label{fund_group}
    F := \faktor{L}{{{M}}}.
\end{align}

\subsection{Two semidirect decompositions of the extended affine Weyl group}
\label{sec_ext_AWG}

Let $T(L)=\{ t_q~|~q\in L\}$.
The \textit{extended (affine) Weyl group} associated to $W$ is 
the subgroup of $\mathrm{GL}(\mathfrak{h}^\ast)$
defined by 
\begin{equation}\label{def extended group}
\widehat{W} := T(L) \rtimes W_0.
\end{equation}

Commonly we will write $\widehat{w}$ for an element in $\widehat{W}$. Any element $\widehat{w} \in \widehat{W}$ decomposes uniquely as $\widehat{w}=t_rw$ with $r \in L$ and $w \in W_0$. 
The group $\widehat{W}$ also acts on the set of alcoves by the same formula as \Cref{def action W sur set alcoves}.
 We denote then
\begin{equation}
\Sigma := \text{Stab}_{\widehat{W}}(\mathcal{A}_e) = \{ f \in  \widehat{W}~|~f\mathcal{A}_e = \mathcal{A}_e\}.
\end{equation}

Let $I = \{1,\dots,n\}$.  We denote, for $j \in I$,
\begin{align}
\Delta_j = \{\alpha_k~|~k \in I~\text{and} ~k \neq j\} \quad \text{and} \quad W_j = \langle s_k~|~\al_k \in \Delta_j \rangle.
\end{align}
The group $W_j$ is a parabolic subgroup of $W_0$. Let $\Phi_j = \Phi_j^+ \sqcup \Phi_j^-$ be the corresponding root system. There exists a unique element $w_{0,j} \in W_j$ such that $w_{0,j}(\Phi_j^+) = \Phi_j^-$. In particular $w_{0,j}$ satisfies $w_{0,j}(\al_k) \in \Phi_0^-$ for any $k \in I \setminus \{j\}$.  This element is the maximal element of the group $W_j$, and then we have $w_{0,j}^{-1} = w_{0,j}$.
 Recall that $\theta$, along with the coefficients $a_i$, are introduced in \Cref{root theta}.

Let us now define $J$ as the following subset of $I$.
\begin{equation}\label{set J}
J = \{0\} \sqcup \left\{
\begin{aligned}
& \{ i \in I~|~ a_i = 1 \}   \text{ \quad  in Type I,}  \\
&  \{ i \in I~|~ a_i^\vee = 1 \}   \text{ \quad  in Type II.}
\end{aligned}
\right.
\end{equation}

Recall that $\varpi_i$ was introduced in \Cref{sec_F}. For $j \in J$, we denote by $\sigma_j \in \widehat{W} $ the element
\begin{equation}\label{def sigmaj_j}
\sigma_j := t_{\varpi_j}w_{0,j}w_0.
\end{equation}

We set the conventions $\varpi_0 := 0$, $\omega_0 := 0$, $w_{0,0} := w_0$ so that $\sigma_0:=1$.
The elements $w_{0,i}$, $\sigma_j$ and $\varpi_j$ verify some crucial properties which we list in the following lemma. Some of these properties are standard in the untwisted cases, but as far as we know, none of them appear uniformly in the literature across both untwisted and twisted types.

\begin{Prop}\label{lemma Garnier}\
    \begin{enumerate}
    \item For any $j \in J$ we have $\varpi_j = \omega_j$.
    \item We have  $\Sigma = \{\sigma_j\, ; \, j \in J\}$. 
        \item Let $i\neq j$. Then $w_{0,i}(\alpha_j) \in -\Delta_i$ and $w_{0,i}(\alpha_i) \in \Phi_0^+ \setminus \Phi_i^+$.  In particular $w_{0,i}\big(\Delta_0 \setminus \{\alpha_i\}\big) \subset - \Delta_0$.
        \item If $j \in J$ then $w_{0,j}(\alpha_j) = \theta$ in Type I, and $w_{0,j}(\alpha_j^{\dag})=\theta$ in Type II.
        \item If $j \in J\backslash \{0\}$ then $\varpi_j \notin {{M}}$.
    \end{enumerate}
\end{Prop}

\begin{proof}
\
\begin{enumerate}
\item This is a direct check using \Cref{varpi in terms of omega} and \Cref{table}.
\item
We apply \cite[Chapter VI, Section  2, no 3, Proposition 6]{BOURB}. In
more precise terms, for cases (i) and (ii) of \Cref{rem_duality_trick} we obtain
that $\Sigma$ consists of the elements $\sigma_i$ defined by
formula \Cref{def sigmaj_j}, where the indices $i$ correspond to the
coefficients equal to $1$ in the decomposition of the highest root of the
dual system to $\Phi_0^\dag$ (which is $\theta$) in terms of the basis $\Delta$.
By definition of $\theta$, these indices are exactly those such that $a_i
= 1$.
In case (iii), \cite[Chapter VI, Section  2, no 3, Proposition 6]{BOURB}
states that the elements of $\Sigma$ are $\sigma_i' = t_{\omega_i} w_{0,i}
w_0$, where $i$ runs over the indices corresponding to the coefficients equal
to $1$ in the decomposition of the highest root of the dual system to
$\Phi_0$, namely $\theta^\dag$, expressed in the basis $\Delta^\dag$.
These are precisely the indices such that $a_i^\vee = 1$. On the other hand, by Point (1), we have
in this case $\omega_i = \varpi_i$, so that $\sigma_i' = \sigma_i$. This
completes the proof of the statement.
\item  This has been proved in \cite[Remark 2.2.2 and Lemma 2.2.5]{Garnier}.
\item 
In cases (i) and (ii) of \Cref{rem_duality_trick}, the result follows directly from \cite[Lemma 2.2.5]{Garnier}, since $\theta$ is the highest root of $\Phi_0$. For case
(iii), we apply \cite[Lemma 2.2.5]{Garnier} to the root system $\Phi_0^\dag$, which shows that
$w_{0,j}(\alpha_j^\dag) = \theta^\dag$, since $\theta^\dag$ is the highest
root of $\Phi_0^\dag$. The result then follows by noting that $\theta =
\theta^\dag$ in this case.
 \item
 Let $j \in J$ and assume that $\varpi_j \in {{M}}$. Since $\sigma_j\mathcal{A}_e = \mathcal{A}_e$, the alcove corresponding to $w_{0,j}w_0 \in W_0$ is translated by $\varpi_j$ to the fundamental alcove. Therefore, a non trivial element of ${{M}}$ translates an element of $W_0$ to another element of $W_0$. This is impossible since $T({{M}}) \cap W_0 = \{1\}$.
\end{enumerate}
\end{proof}

\begin{Exa} Let us illustrate \Cref{lemma Garnier} (1) on all rank $2$ types, depicted in \Cref{fig:rank2}.
In type $A_2^{(1)}$, we have $J=\{1,2\}$ and one checks that $\Sigma=\{1,\sigma_1,\sigma_2\}$.
In the following three cases, we have $J=\{2\}$, and we see that 
we indeed have $\Sigma=\{1,\sigma_2\}$.
In the two remaining types, $G_2^{(1)}$ and $D_4^{(3)}$
have trivial fundamental group.
\begin{figure}
\centering
\begin{center}
\begin{tabular}{cc} 
Type $A_2^{(1)}$
& 
Type $C_2^{(1)}$
\\
\begin{tikzpicture}[scale=1]
\tiny
\clip (-3,-3) rectangle (3,3);
\path[fill=gray!30] (0,0) -- +(0:1) -- ++(60:1);
\draw[line width = 0.3mm, ->] (0,0) -- + (-30:{sqrt(3)});
\draw[line width = 0.3mm, ->] (0,0) -- + (90:{sqrt(3)});
\draw[line width = 0.3mm, ->] (0,0) -- + (30:{sqrt(3)});
\node[anchor = north west, scale=1.5] at ( 1.5, {-sqrt(3)/2}  ) {$\alpha_1$};
\node[anchor = south, scale=1.5] at ( 0, {sqrt(3)} ) {$\alpha_2$};
\node[anchor = south west, scale=1.5] at ( 1.5, {sqrt(3)/2} ) {$\theta$};
\draw[line width = 0.3mm, ->, color=purple!80] (0,0) -- + (0:1);
\node[anchor = south west, scale=1.5, color=purple!80] at
(1,0) {$\varpi_1$};
\draw[line width = 0.3mm, ->, color=teal] (0,0) -- + (60:1);
\node[anchor = south, scale=1.5, color=teal] at
(0.5,{sqrt(3)/2}) {$\varpi_2$}; 
\foreach \i in {-6,...,6}
{
\draw[color=gray, line width = 0.01mm] (\i,0) -- +(60:6) -- ++(60:-6);
\draw[color=gray, line width = 0.01mm] (-6,{\i*sqrt(3)/2}) -- (6,{\i*sqrt(3)/2});
\draw[color=gray, line width = 0.01mm] ({\i*3/4},{\i*sqrt(3)/4}) -- +(120:6) -- ++(120:-6);
}
\end{tikzpicture}
&  
\begin{tikzpicture}[scale=1]
\tiny
\clip (-3,-3) rectangle (3,3);
\path[fill=gray!30] (0,0) -- +(0:1) -- ++(45:{1*sqrt(2)});
\foreach \i in {-6,...,6}
{
\draw[color=gray, line width = 0.01mm] (\i,-\i) -- + (45:10) -- ++ (45:-10);
\draw[color=gray, line width = 0.01mm] (\i,\i) -- + (-45:10) -- ++ (-45:-10);
\draw[color=gray, line width = 0.01mm] (0,\i) -- + (0:10) -- ++ (0:-10);
\draw[color=gray, line width = 0.01mm] (\i,0) -- + (90:10) -- ++ (90:-10);
}
\draw[line width = 0.3mm, ->] (0,0) -- + (-45:{sqrt(2)});
\node[anchor = north west, scale=1.5] at
(1,-1) {$\alpha_1$};
\draw[line width = 0.3mm, ->] (0,0) -- + (90:2);
\node[anchor = south, scale=1.5] at
(0,2) {$\alpha_2$};
\draw[line width = 0.3mm, ->] (0,0) -- + (0:2);
\node[anchor = north west, scale=1.5] at
(2,0) {$\theta$};
\draw[line width = 0.3mm, ->, color=purple!80] (0,0) -- + (0:1.97);
\node[anchor = south west, scale=1.5, color=purple!80] at
(2,0) {$\varpi_1$};
\draw[line width = 0.3mm, ->, color=teal] (0,0) -- + (45:{sqrt(2)});
\node[anchor = south, scale=1.5, color=teal] at
(1,1) {$\varpi_2$}; 
\end{tikzpicture} 
\\
Type $D_3^{(2)}$
& 
Type $A_4^{(2)}$
\\
\begin{tikzpicture}[scale=1]
\tiny
\clip (-3,-3) rectangle (3,3);
\path[fill=gray!30] (0,0) -- +(0:1) -- ++(45:{sqrt(2)});
\foreach \i in {-4,...,4}
{
\draw[color=gray, line width = 0.01mm] (\i,-\i) -- + (45:10) -- ++ (45:-10);
\draw[color=gray, line width = 0.01mm] (\i,\i) -- + (-45:10) -- ++ (-45:-10);
\draw[color=gray, line width = 0.01mm] (0,\i) -- + (0:10) -- ++ (0:-10);
\draw[color=gray, line width = 0.01mm] (\i,0) -- + (90:10) -- ++ (90:-10);
}
\draw[line width = 0.3mm, ->] (0,0) -- + (-45:{2*sqrt(2)});
\node[anchor = north west, scale=1.5] at
(2,-2) {$\alpha_1$};
\draw[line width = 0.3mm, ->] (0,0) -- + (90:2);
\node[anchor = south, scale=1.5] at
(0,2) {$\alpha_2$};
\draw[line width = 0.3mm, ->] (0,0) -- + (0:2);
\node[anchor = south west, scale=1.5] at
(2,0) {$\theta$};
\draw[line width = 0.3mm, ->, color=purple!80] (0,0) -- + (0:1);
\node[anchor = south west, scale=1.5, color=purple!80] at
(1,0) {$\varpi_1$};
\draw[line width = 0.3mm, ->, color=teal] (0,0) -- + (45:{sqrt(2)});
\node[anchor = south, scale=1.5, color=teal] at
(1,1) {$\varpi_2$}; 
\end{tikzpicture}
&  
\begin{tikzpicture}[scale=1]
\tiny
\clip (-3,-3) rectangle (3,3);
\path[fill=gray!30] (0,0) -- +(0:0.5) -- ++(45:{0.5*sqrt(2)});
\foreach \i in {-6,...,6}
{
\draw[color=gray, line width = 0.01mm] (0.5*\i,-0.5*\i) -- + (45:10) -- ++ (45:-10);
\draw[color=gray, line width = 0.01mm] (0.5*\i,0.5*\i) -- + (-45:10) -- ++ (-45:-10);
\draw[color=gray, line width = 0.01mm] (0,0.5*\i) -- + (0:10) -- ++ (0:-10);
\draw[color=gray, line width = 0.01mm] (0.5*\i,0) -- + (90:10) -- ++ (90:-10);
}
\draw[line width = 0.3mm, ->] (0,0) -- + (-45:{sqrt(2)});
\node[anchor = north west, scale=1.5] at
(1,-1) {$\alpha_1$};
\draw[line width = 0.3mm, ->] (0,0) -- + (90:2);
\node[anchor = south, scale=1.5] at
(0,2) {$\alpha_2$};
\draw[line width = 0.3mm, ->] (0,0) -- + (45:{sqrt(2)});
\node[anchor = south west, scale=1.5] at
(1,1) {$\theta$};
\draw[line width = 0.3mm, ->, color=purple!80] (0,0) -- + (0:1);
\node[anchor = west, scale=1.5, color=purple!80] at
(1,0) {$\varpi_1$};
\draw[line width = 0.3mm, ->, color=teal] (0,0) -- + (45:{0.5*sqrt(2)});
\node[anchor = south, scale=1.5, color=teal] at
(0.5,0.5) {$\varpi_2$}; 
\end{tikzpicture}
\\
Type $G_2^{(1)}$
& 
Type $D_4^{(3)}$
\\
\begin{tikzpicture}[scale=0.6666]
\tiny
\clip (-4.5,-4.5) rectangle (4.5,4.5);
\path[fill=gray!30] (0,0) -- +(0:1.5) -- ++(30:{sqrt(3)});
\foreach \i in {-10,...,10}
{
\draw[color=gray, line width = 0.01mm] (3*\i,0) -- +(-30:10) -- ++(-30:-10);
\draw[color=gray, line width = 0.01mm] (0,{3*\i*sqrt(3)/2}) -- +(0:6) -- ++(0:-6);
\draw[color=gray, line width = 0.01mm] (3*\i,0) -- +(30:10) -- ++(30:-10);
\draw[color=gray, line width = 0.01mm] (3*\i,0) -- +(60:6) -- ++(60:-6);
\draw[color=gray, line width = 0.01mm] (1.5*\i,0) -- +(90:6) -- ++(90:-6);
\draw[color=gray, line width = 0.01mm] (3*\i,0) -- +(120:6) -- ++(120:-6);
}
\draw[line width = 0.3mm, ->] (0,0) -- + (-60:3);
\draw[line width = 0.3mm, ->] (0,0) -- + (90:{sqrt(3)});
\draw[line width = 0.3mm, ->] (0,0) -- + (0:3);
\node[anchor = south west, scale=1.5] at ( 1.5, {-3*sqrt(3)/2}  ) {$\alpha_1$};
\node[anchor = south, scale=1.5] at ( 0, {sqrt(3)} ) {$\alpha_2$};
\node[anchor = south east, scale=1.5] at ( 3,0 ) {$\theta$};
\draw[line width = 0.3mm, ->, color=purple!80] (0,0) -- + (0:2.97);
\node[anchor = north east, scale=1.5, color=purple!80] at (3,0) {$\varpi_1$};
\draw[line width = 0.3mm, ->, color=teal] (0,0) -- (4.5,{3*sqrt(3)/2});
\node[anchor = south east, scale=1.5, color=teal] at( 4.5, {3*sqrt(3)/2} ) {$\varpi_2$}; 
\end{tikzpicture}
&
\begin{tikzpicture}[scale=1]
\tiny
\clip (-3,-3) rectangle (3,3);
\path[fill=gray!30] (0,0) -- +(30:{sqrt(3)/2}) -- ++(60:1);
\foreach \i in {-10,...,10}
{
\draw[color=gray, line width = 0.01mm] (3*\i,0) -- +(-30:6) -- ++(-30:-6);
\draw[color=gray, line width = 0.01mm] (0,{\i*sqrt(3)/2}) -- +(0:6) -- ++(0:-6);
\draw[color=gray, line width = 0.01mm] (3*\i,0) -- +(30:6) -- ++(30:-6);
\draw[color=gray, line width = 0.01mm] (\i,0) -- +(60:6) -- ++(60:-6);
\draw[color=gray, line width = 0.01mm] (1.5*\i,0) -- +(90:6) -- ++(90:-6);
\draw[color=gray, line width = 0.01mm] (\i,0) -- +(120:6) -- ++(120:-6);
}
\draw[line width = 0.3mm, ->] (0,0) -- + (-30:{sqrt(3)});
\draw[line width = 0.3mm, ->] (0,0) -- + (120:3);
\draw[line width = 0.3mm, ->] (0,0) -- + (30:{sqrt(3)});
\node[anchor = south west, scale=1.5] at ( 1.5, {-sqrt(3)/2}  ) {$\alpha_1$};
\node[anchor = south, scale=1.5] at ( -1.5, {3*sqrt(3)/2} ) {$\alpha_2$};
\node[anchor = south west, scale=1.5] at (1.5,{sqrt(3)/2} ) {$\theta$};
\draw[line width = 0.3mm, ->, color=purple!80] (0,0) -- + (30:{sqrt(3)-0.03}); 
\node[anchor = north west, scale=1.5, color=purple!80] at (1.5,{sqrt(3)/2}) {$\varpi_1$};
\draw[line width = 0.3mm, ->, color=teal] (0,0) -- + (60:1); 
\node[anchor = south, scale=1.5, color=teal] at (0.5,{sqrt(3)/2}) {$\varpi_2$};
\end{tikzpicture}
\end{tabular}
\end{center}
\caption{All rank $2$ root systems.
We have shaded in grey the fundamental alcove.}
\label{fig:rank2}
\end{figure}
\end{Exa}

The following result is essential.
\bigskip
\begin{Prop}\label{L in terms of M}
\
\begin{enumerate}
\item The set $\{\varpi_j \,;\, j\in J\}$ is a complete set of coset representatives of $L/{{M}}$.
\item We have $L = \bigsqcup\limits_{j \in J}(\varpi_j+w_{0,j}w_0({{M}}))$.
\item We have $\widehat{W} = \Sigma \ltimes W$. In particular $\Sigma \simeq F$.
\end{enumerate}
\end{Prop}

\begin{proof}

To prove (1), we apply \cite[Corollary VI.2.3]{BOURB} to the root systems
appearing in \Cref{rem_duality_trick}. More precisely, in cases (i)
and (ii), we directly obtain that the set of elements $\varpi_i$, for $i
\in J$, forms a system of representatives of $L/M$. In case (iii),
\cite[Corollary VI.2.3]{BOURB} gives that the set of elements $\omega_i$
for $i\in J$ is a set of coset representatives of $L/M$, and the result
then follows from \Cref{lemma Garnier}(1).
Point (2) is an immediate consequence of Point (1) and the fact that
$w_{0,j}w_0$ is a bijection of $M$.
Finally, in order to prove (3), we 
show that  $W$ is a normal subgroup of $\widehat{W}$, that $\Sigma \cap W = \{1\}$ and that  $\widehat{W}=\Sigma W$. 
\begin{itemize}
    \item Since $\widehat{W}=T(L)\rtimes W_0$, to prove normality, it suffices to show that $t_r w t_r^{-1} \in W$ for all $w\in W$ and for all $r\in L$.
Write $w = t_q\overline{w}$ with $q \in {{M}}$ and $\overline{w} \in W_0$. Then
        $$
        t_r w t_r^{-1} = t_r t_q\overline{w}t_r^{-1} = t_{r+q-\overline{w}(r)}\overline{w},
        $$
and this belongs to $W$ 
if and only if $r-\overline{w}r\in {{M}}$.
It suffices to show this for $\overline{w}=s_i$ for $i = 1,2,\dots, n$. 
We have by \Cref{reflexion s_alpha}
\begin{equation}\label{r-s_i(r)}
            r-s_i(r) = \frac{2 }{(\al_i~|~\al_i)} (r~|~\al_i) \al_i.
\end{equation}
We have two cases to consider:
        \begin{itemize}
            \item Assume that we are in Type I.  First, rewrite \Cref{r-s_i(r)} as
            \begin{equation*}
                r-s_i(r) = (r~|~\al_i)\al_i^{\dagger}.
            \end{equation*}
            Since $r\in L = \bigoplus_{i=1}^n\mathbb{Z}\varpi_i$, we have $(r~|~\al_i) \in \mathbb{Z}$,
            and thus $r-s_i(r)  \in \bigoplus_{i=1}^n\mathbb{Z}\alpha_i^{\dagger} =  {{M}}$.
            
            \item Assume that we are in Type II.    First, rewrite \Cref{r-s_i(r)} as
            \begin{equation*}
                r-s_i(r) = (r~|~\al_i^\dagger)\al_i.
            \end{equation*}
            Since $r\in L = \bigoplus_{i=1}^n\mathbb{Z}\om_i$, we have $(r~|~\al_i^\dagger) \in \mathbb{Z}$,
            and thus $r-s_i(r)  \in \bigoplus_{i=1}^n\mathbb{Z}\alpha_i =  {{M}}$.
        \end{itemize}

    \item Let us show that $\Sigma \cap W = \{1\}$. Let $\sigma_j \in \Sigma$ and $w = t_q\overline{w} \in W$ with $q \in {{M}}$ and $\overline{w} \in W_0$. We have
    $$
    \sigma_j = w \Longleftrightarrow t_{\varpi_j}w_{0,j}w_0 = t_q\overline{w} \Longleftrightarrow t_{\varpi_j - q} = \overline{w}w_0w_{0,j}.
    $$
    However, $\overline{w}w_0w_{0,j} \in W_0$ and $\varpi_j - q \in L$. Since $T(L) \cap W_0 = \{1\}$, it follows that 
    $$
    \sigma_j = w \Longleftrightarrow \varpi_j = q \quad \text{and} \quad w_{0,j}w_0 = \overline{w}. 
    $$
    But now, by \Cref{lemma Garnier} (5), the equality $\varpi_j = q$ is only possible if $j=0$. Hence the result. 

    \item Finally, we show that $\Sigma W=\widehat W$. 
    Let $\sigma_jt_q\overline{w} \in \Sigma W$. This element is equal to $t_{\varpi_j + w_{0,j}w_0\overline{w}(q)}w_{0,j}w_0\overline{w}$, 
    and in order to show the inclusion $\Sigma W \subset \widehat W$ we only need to show that $\varpi_j + w_{0,j}w_0\overline{w}(q)$ belongs to $L$, which is clear by (2) and \Cref{MvsL}.
    For the reverse inclusion, (2) and \Cref{MvsL} ensure that any $r\in L$ writes $r=\varpi_j + q$ for some $j\in J$ and $q\in M$.
    Therefore $t_{r}=t_{\varpi_j}t_q$, which
        belongs to $\Sigma W$ since
        $t_{\varpi_j}\in \Sigma W$ and
        $t_q\in W$.
    
    \end{itemize}
\end{proof}

 \begin{Exa}[$\Sigma$ in type $A_n^{(1)}$]\label{section Sigma type A}
Using for instance \cite{BOURB}, we find the decompositions of the fundamental weights
\begin{equation}\label{fund_weights_A}
\begin{aligned}
\varpi_i & = \frac{n+1-i}{n+1} \left(\alpha_1 + 2\alpha_2 + \dots + i\alpha_i \right) + \frac{i}{n+1}\left(\alpha_n + 2\alpha_{n-1} + \dots + (n-i)\alpha_{i+1} \right)
\\
& = \left(1-\frac{i}{n+1}\right)e_1 + \dots + \left(1-\frac{i}{n+1}\right)e_i - \frac{i}{n+1}e_{i+1} - \dots - \frac{i}{n+1}e_{n+1}.
\end{aligned}
\end{equation}

An easy computation, for $j =1,\dots, n$, shows that
\begin{equation}
    \varpi_j - j\varpi_1 = \sum\limits_{k=1}^{j-1}\bigl(k-j\bigr)\alpha_k.
\end{equation}
Therefore, $[{\varpi_j}] = j[{\varpi_1}]$ in $L/{{M}}$, or equivalently $\sigma_j = \sigma_1^j$, showing with \Cref{lemma Garnier} (2) that $\Sigma = \langle \sigma_1 \rangle$.

\medskip

Moreover, for all $i=0,1,\ldots, n$, the matrix of $w_{0,i}w_0$ in the canonical basis $(e_1,\ldots, e_{n+1})$ of $\R^{n+1}$ is $A^i$ where
\begin{equation}\label{w01w0_can_A}
A = \mathrm{Mat}_{(e_1,\ldots, e_{n+1})} (w_{0,1}w_0) =
\begin{pmatrix}
0 &0 & \cdots & 0 & 1
\\
1&0&\cdots & 0&0
\\
0&1&\cdots & 0& 0
\\
\vdots & & \ddots & &\vdots
\\
0&0&\cdots &1&0
\end{pmatrix}.
\end{equation}

We find $J=\{0,1,\ldots,n\}$ and therefore
\begin{equation}\label{Sigma type A}
\Sigma = \{\sigma_j \, ; \, j=0,1,\dots, n\} = \{1,~t_{\varpi_j} w_{0,j}w_0\, ; \, j=1,\dots, n\}.
\end{equation}

Let us detail the case $n=2$.
We have $I = \{1,2\}$ and $J = \{0,1,2\}$.   
Then $\sigma_1 = t_{\varpi_1}w_{0,1}w_0 = t_{\varpi_1} \cdot s_2\cdot s_2s_1s_2  = t_{\varpi_1}s_1s_2$ 
and $\sigma_2 = t_{\varpi_2}w_{0,2}w_0 = t_{\varpi_2} \cdot s_1\cdot s_1s_2s_1  = t_{\varpi_2}s_2s_1$.  Hence
\begin{equation}\label{Sigma in A2}
\Sigma = \{1, ~ t_{\varpi_1}s_1s_2, ~t_{\varpi_2}s_2s_1\}.
\end{equation}

An easy computation shows that $\sigma_1^2 = \sigma_2$ and $\sigma_2^2 = \sigma_1$ and $\sigma_1^3 = \sigma_2^3 = 1$. Hence $\Sigma = \langle \sigma_1 \rangle \simeq \mathbb{Z}/3\mathbb{Z}$. 
\end{Exa}

 \begin{Exa}[$\Sigma$ in type $C_n^{(1)}$]\label{section Sigma type C}
We find $J = \{0,n\}$. It follows that
\begin{equation}\label{Sigma type C}
\Sigma = \{1,  \sigma_n\} = \{1, t_{\varpi_n} w_{0,n}w_0\}
\end{equation}

where $w_{0,n}$ is the longest element of $W_n = \langle s_1,s_2,\dots,s_{n-1}\rangle \simeq A_{n-1}$ and $w_0$ the longest element of $W(C_n)$.  

The matrix of $w_{0,n}w_0$ in the canonical basis $(e_1,\dots, e_n)$ of $\mathbb{R}^n $ is given by
\begin{equation}
    \mathrm{Mat}_{(e_1,\ldots, e_{n})} (w_{0,n}w_0) = 
  \begin{pmatrix}[cccc]
 0           & \cdots        & 0            &  -1  \\
 0           & \cdots        & -1         &   0 \\
 \vdots   & \iddots       &  \vdots    & \vdots  \\
 -1          & \cdots        &    0          & 0
\end{pmatrix}.
\end{equation}

  \end{Exa}

\subsection{Extended affine Grassmannian}
\label{sec_ext_AG}

We define here a natural notion of extended affine Grassmannian.  
Since we cannot build the semidirect product of $\Sigma$ with $W^0$ (because $W^0$ does not enjoy a group structure),  
we consider instead the quotient $\widehat{W}/W_0$, which naturally extends the original quotient $W/W_0$,
and call it the  \textit{(right) extended affine Grassmannian}.

\medskip

Recall that any element in $\Sigma \ltimes W$ decomposes uniquely as $\sigma_j w$ for $w = t_qw_0 \in W$, where $q \in {{M}}$, $w_0 \in W_0$, 
and $j\in J$. 
As $\Sigma$ is the set of elements $\sigma$ of length $\ell(\sigma) = 0$, 
the following maps are natural bijections

\begin{equation}\label{bij can extended}
\begin{array}{ccccccc}
\faktor{(\Sigma \ltimes W)}{W_0}   & \longrightarrow  & \Sigma \times  \faktor{W}{W_0} &    & \longrightarrow       & \Sigma W^0 \\
 \left[\sigma_j w \right]            & \longmapsto      & (\sigma_j,[w])                   &    & \longmapsto           &\sigma_j w^0
\end{array}
\end{equation}

and we can consider the composition map $\widehat{\chi}$, which is the extended analogue of \Cref{def chi}

\begin{equation}\label{bij 2 can extended}
\begin{array}{ccccc}
\widehat{\chi} & : & \faktor{\widehat{W}}{W_0} & \longrightarrow    & \Sigma W^0 \\
     &   &  \left[\sigma w \right]   & \longmapsto        &\sigma w^0.
\end{array}
\end{equation}

The elements of $\Sigma W^0$ are called \textit{extended affine Grassmannian elements} and will be of central interest in the rest of the paper. 

\medskip

Moreover, by \Cref{L in terms of M} (2), we have the following bijection (similar to \Cref{map_f})
$$
\begin{array}{cccc}
\widehat f : & \faktor{\widehat W}{W_0}  & \longrightarrow & L    \\
  &    \left[\widehat w\right]       & \longmapsto       &   \varpi_j+w_{0,j}w_0(q)
  \end{array}
$$
where $\widehat w=\si_j t_q\overline{w} = t_{\varpi_j + w_{0,j}w_0(q)} w_{0,j}w_0\overline{w}$ with $j\in J, q\in {{M}}$ and $\overline{w}\in W_0$.


\section{Atomic length and extended atomic length}
\label{sec_AL}

In \cite{CG2022}, a new statistic on finite and affine Weyl groups, the \textit{atomic length}, 
was introduced and studied. In this section we recall the general definition of this statistic and we extend it to the extended affine Weyl group. We then generalise some  results obtained in \cite{CG2022} (\Cref{lem_aff_AL_ext} and \Cref{analogue cor 8.2}) and we give a new theorem that will be of importance in the following sections, namely \Cref{decompo atomic D}. 
Finally, the quadratic formula established in \Cref{section gaussian equation gen} is one of the main characters of this paper,
and will be thoroughly reinvested in \Cref{sec_pell} and \Cref{sec_pell_A3}.

\subsection{Atomic length}

 Let $W$ be an affine Weyl group of rank $n$ and let $\La$
be a dominant weight.
We set, for all $w\in W$,
\begin{equation}\label{def atomic length}
\sL_\La (w) = \langle \La- w\La , \rho^\vee \rangle.
\end{equation}

Recall that any $w \in W$ decomposes uniquely as $w=t_q \overline{w}$, where $\be\in {{M}}, \overline{w}\in W_0$. Recall also that $\h(q)$ is the \textit{height} of $\be$.
The following result is \cite[Cor. 8.2]{CG2022}

\begin{Th}\label{thm_LA}
For all $w=t_q \overline{w}\in W$, we have
\begin{equation}\label{formula thm_LA}
    \sL_{\Lambda_0}(w) = \frac{h}{2}|\be|^2 - \h(q).
\end{equation}

In particular $\sL_{\Lambda_0}(w) =  \sL_{\Lambda_0}(t_q)$.
\end{Th}

Note that $\sL_{\Lambda_0}(w)$ depends only on $\be$.
It turns out that this coincides with the formula established independently by Stucky, Thiel and Williams 
in their study of generalised cores \cite{STW2023}.

\medskip 

We will also need a similar formula for $\sL_{\Lambda_i}$.
Recall that $$\La_i=\frac{a_i^\vee}{a_0^\vee} \La_0 + \om_i$$
where $\om_i$ is the $i$-th fundamental weight of the corresponding finite root system.

\begin{Prop}\label{thm_LAi}
For all $w=t_q \overline{w}\in W$, we have
$$\sL_{\Lambda_i}(w) = \frac{a_i^\vee}{a_0^\vee}\left(\frac{h}{2}|\be|^2 - \h(q) \right)+h(\overline{w}(\om_i)\,|\,q) + 
\h(\om_i-\overline{w}(\om_i)).$$
In particular, we have
$$\sL_{\Lambda_i}(t_q) = \frac{a_i^\vee}{a_0^\vee}\left(\frac{h}{2}|\be|^2 - \h(q) \right) +h(\om_i\,|\,q).$$
\end{Prop}

\begin{proof}
We have
\begin{align*}
\sL_{\La_i}(w) 
& = \left\langle \La_i-w(\La_i) , \rho^\vee\right\rangle
\\
& = \left\langle \frac{a_i^\vee}{a_0^\vee}\La_0+\om_i-w\left(\frac{a_i^\vee}{a_0^\vee}\La_0+\om_i\right) , \rho^\vee\right\rangle
\\
& = \frac{a_i^\vee}{a_0^\vee}\left(
\left\langle \La_0-w(\La_0) , \rho^\vee\right\rangle
\right)
+ 
\left\langle \om_i-w(\om_i) , \rho^\vee\right\rangle
\\
& = \frac{a_i^\vee}{a_0^\vee}\sL_{\La_0}(w)
+ 
\left\langle 
\om_i-\left[ \overline{w}(\om_i) - (\overline{w}(\om_i) \,|\, q)\delta \right] , \rho^\vee 
\right\rangle
\\
& = \frac{a_i^\vee}{a_0^\vee}\sL_{\La_0}(w)
+ h(\overline{w}(\om_i) \,|\, q) + \h(\om_i-\overline{w}(\om_i))
\end{align*}
where the penultimate identity is obtained by applying Formula \Cref{translation}
to expand $w(\om_i)=t_q(\overline{w}(\om_i))$.
We conclude using \Cref{thm_LA}.
\end{proof}

\begin{Rem}\label{rem_AL_AG}
The $\La_0$-atomic length has a particularly easy description on the affine Grassmannian. 
Recall that $W/W_0$ is in bijection with $W^0$ via the map $\chi: [w]\mapsto w^0$ (see \Cref{def chi}).
Now, let $w = t_q\overline{w} = w^0w_0\in W $ where $q\in {{M}}$ and $\overline{w} \in W_0$. Therefore, $w^0 = (t_q)^0$ with $w_0 = (t_q)_0\overline{w}$, and it follows that
\begin{equation}\label{LA_GA}
\sL_{\La_0}( w^0) = \sL_{\La_0}( (t_q)^0)) = \sL_{\La_0} (t_q(t_q)_0^{-1}) = \sL_{\La_0}(t_q)
\end{equation}
where the last equality follows from \Cref{thm_LA}.
In fact, \Cref{thm_LA} provides a polynomial formula in $q$, so it is natural to consider
the map $\mu:{{M}}\to \N, q\mapsto \sL_{\La_0}(t_q)$.
In other terms, combining with \Cref{LA_GA}, the following diagram commutes.
\begin{equation}
\label{link Q et L}
\begin{array}{c}
\begin{tikzpicture}
\node at (0,0) {$ \faktor{W}{W_0}$} ;
\node at (3,0) {$W^0$} ;
\node at (3,-2) {$\mathbb{N}$.} ;
\node at (0,-2) {${{M}}$} ;

\node at (0.5,0) (1)  {} ;
\node at (2.5, 0) (2) {} ; 
\draw [>=stealth,->] (1) to (2);
\node at (1.5,0.3) {$\chi$};

\node at (0,-0.3) (3)  {} ;
\node at (0, -1.8) (4) {} ; 
\draw [>=stealth,->] (3) to (4);
\node at (-0.4,-1) {$f$};

\node at (3,-0.3) (5)  {} ;
\node at (3, -1.8) (6) {} ; 
\draw [>=stealth,->] (5) to (6);
\node at (3.5,-1) {$\sL_{\Lambda_0}$};

\node at (0.5,-2) (7)  {} ;
\node at (2.5, -2) (8) {} ; 
\draw [>=stealth,->] (7) to (8);
\node at (1.5,-2.3) {$\mu$};
\end{tikzpicture}
\end{array}
\end{equation}

\end{Rem}

Finally, the following set will be of importance in \Cref{sec_param}.

\begin{Def}
\label{definition B sans hat}
    Let $N \in \mathbb{N}$. We define the set $\mathcal{B}(N)$ as follows 
    \begin{equation}\label{extended B(N) sans hat}
\mathcal{B}(N) = \{q \in {{M}} \mid \sL_{\Lambda_0}(t_q) = N \}.
\end{equation}
\end{Def}

\subsection{Extended atomic length}
\label{sec_ext_AL}

Let $\Lambda$ be a dominant weight. Formula \Cref{def atomic length} also makes sense for any element $g \in \text{GL}(\mathfrak{h}^*)$:
$$
\sL_{\Lambda}(g) = \langle \Lambda - g\Lambda,\rho^{\vee} \rangle.
$$
We will be particularly interested in studying the atomic length on the extended Weyl group $\widehat{W} \subset \text{GL}(\mathfrak{h}^*)$. 
Similarly, we can define the height of any element of $\mathfrak{h}^*$,
which will be a real number but not necessarily an integer.

\medskip

Let $\Lambda \in P^+$, so that
\begin{equation}\label{def Lambda}
\Lambda = \lambda + \ell\Lambda_0 + z\delta,
\end{equation}
for some $\lambda \in P_0^+$, $\ell\in \N$ and $z \in a_0^{-1}\mathbb{Z}$.
Note that $\ell$, the level of $\Lambda$, equals $\langle \Lambda, c\rangle$. 

\medskip

The following Lemma is the extended version of \cite[Lemma 8.1]{CG2022}.

\begin{Lem}\label{lem_aff_AL_ext}
Let $x\in \mathfrak{h}^*$, $\overline{w} \in W_0 $ and set $w=t_x \overline{w} \in \mathrm{GL}(\mathfrak{h}^*)$.
We have
$$
\sL_\La( w) 
=
\sL_{\la}(\overline{w})
- \ell \ \h (x)
+ h \left( ( \la \mid \overline{w}^{-1}(x)) + \frac{1}{2}|x|^2\ell\right).
$$
\end{Lem}

\begin{proof}
    The proof is exactly the same as that of \cite[Lemma 8.1]{CG2022}. The only difference is that $x \in \mathfrak{h}^*$ instead of ${{M}}$ but this does not affect the properties used.
\end{proof}

We define the following map
$$
\begin{array}{ccccc}
D_{\Lambda} & : & W_0 \times \mathfrak{h}^* \times \mathfrak{h}^* & \longrightarrow & \mathbb{R}\\
                      &   & (u,x,y) & \longmapsto &  h\cdot\ell~(x~|~y) - \sL_{\lambda}(u).
\end{array}
$$

Write $ b_i = a_i^{\vee}/a_0^{\vee}$.  We see directly, with the convention $\omega_0 = 0$,  that for any $0 \leq i \leq n$ one has
\begin{align*}
\mathcal{D}_{\Lambda_i}(u,x,y) & =   h\cdot b_i(x~|~y) - \sL_{\omega_i}(u) \\
& =  h\cdot b_i(x~|~y)  - \langle \omega_i - u(\omega_i),\rho^{\vee}\rangle \\
& =  h\cdot b_i(x~|~y)  - (\langle \omega_i, \rho^{\vee}\rangle - \langle u(\omega_i), \rho^{\vee}\rangle) \\
& =  h\cdot b_i(x~|~y)   + \text{ht}(u(\omega_i)) - \text{ht}(\omega_i). 
\end{align*}

In particular for $\Lambda_0$ we get

\begin{equation}\label{defect term for Lambda0}
\mathcal{D}_{\Lambda_0}(u,x,y)  =  h\cdot \frac{a_0^{\vee}}{a_0^{\vee}}(x~|~y)   + \text{ht}(u(\omega_0)) - \text{ht}(\omega_0)  = h\cdot (x~|~y).
\end{equation}

\bigskip

\begin{Th}\label{decompo atomic D}
Let $x,y \in \mathfrak{h}^*$ and $u \in W_0$.  We have
$$
\sL_{\Lambda}(t_{x+y}u) = \sL_{\Lambda}(t_{x}u) + \sL_{\Lambda}(t_{y}u)  + \mathcal{D}_{\Lambda}(u,x,y).
$$
\end{Th}

\begin{proof}
This follows from the following computation, starting with \Cref{lem_aff_AL_ext}.
\begin{align*}
\sL_{\Lambda}(t_{x+y}u)   & =  \sL_{\lambda}(u) - \ell~\text{ht}(x+y)  + h\Big((\lambda~|~\overline{w}^{-1}(x+y)) + \frac{1}{2}~\ell~|x+y|^2\Big) \\
& =  \sL_{\lambda}(u) -\ell ~\text{ht}(x)- \ell ~\text{ht}(y) + h\Big((\lambda~|~\overline{w}^{-1}(x)) + (\lambda~|~\overline{w}^{-1}(y)) + \frac{1}{2}~\ell~|x|^2 +\ell~(x~|~y) + \frac{1}{2}~\ell~|y|^2\Big)  \\
& = \sL_{\Lambda}(t_{x}u) -\ell~\text{ht}(y) + h\Big((\lambda~|~\overline{w}^{-1}(y))  + \frac{1}{2}~\ell~|y|^2\Big) + h~\ell~(x~|~y) \\
& = \sL_{\Lambda}(t_{x}u) +\Big(\sL_{\lambda}(u) -\ell~\text{ht}(y) + h\Big[(\lambda~|~\overline{w}^{-1}(y))  + \frac{1}{2}~\ell~|y|^2\Big]\Big) - \sL_{\lambda}(u) + h~\ell~(x~|~y) \\ 
& =  \sL_{\Lambda}(t_{x}u) +  \sL_{\Lambda}(t_{y}u) +  h~\ell~(x~|~y) - \sL_{\lambda}(u) \\
& = \sL_{\Lambda}(t_{x}u) + \sL_{\Lambda}(t_{y}u)  + \mathcal{D}_{\Lambda}(u,x,y).
\end{align*}
\end{proof}

In particular,  \Cref{decompo atomic D} shows that for any $0 \leq i \leq n$ one has
\begin{equation}
\sL_{\Lambda_i}(t_{x+y}) = \sL_{\Lambda_i}(t_x) + \sL_{\Lambda_i}(t_y)  + h\cdot \frac{a_i^{\vee}}{a_0^{\vee}} (x~|~y).
\end{equation}
In other terms, up to a defect term $D_{\Lambda}(1,x,y)$, the $\Lambda$-atomic length $\sL_{\Lambda}$ restricted to $T(\mathfrak{h}^*)$ is a group morphism from $T(\mathfrak{h}^*)$ to $\mathbb{Z}$.  

\bigskip

\begin{Prop} \label{atomic extended}
Let $\Lambda\in P^+$ be a dominant weight of level $\ell$, and $x\in \mathfrak{h}^\ast$.  
We have
$$
\sL_{\Lambda}(t_x) =  h(\Lambda~|~x) + \ell \left(\frac{h}{2}|x|^2 - \mathrm{ht}(x)\right).
$$
\end{Prop}

\begin{proof}
We have
\begin{align*}
\sL_{\Lambda}(t_{ x}) & =   \langle \Lambda - t_{ x}(\Lambda),\rho^{\vee} \rangle \\
	  & =   \langle \Lambda -\left( \Lambda + \langle \Lambda,c \rangle x - \left( (\Lambda~|~ x)+\frac{1}{2}| x|^2\langle \Lambda,c\rangle \right)\delta  \right),\rho^{\vee} \rangle \\
  & =   \langle - \langle \Lambda,c \rangle x + \left( (\Lambda~|~ x)+\frac{1}{2}| x|^2\langle \Lambda,c\rangle \right)\delta,\rho^{\vee} \rangle \\
  & =  \langle - \ell x + \left( (\lambda~|~ x)+\frac{1}{2}| x|^2\ell \right)\delta,\rho^{\vee} \rangle \\
  & = -\ell\langle x,\rho^{\vee} \rangle + \left( (\Lambda~|~ x)+\frac{1}{2}| x|^2\ell \right) \langle \delta, \rho^{\vee} \rangle \\ 
  & = -\ell~\langle x,\rho^{\vee} \rangle + \left( (\Lambda~|~ x)+\frac{1}{2}| x|^2\ell \right)h \\ 
  & = -\ell~\text{ht}( x) + \left( (\Lambda~|~ x)+\frac{1}{2}| x|^2\ell \right)h. 
\end{align*}
\end{proof}

The following proposition is the extended version of \Cref{thm_LA}.

\begin{Prop}\label{analogue cor 8.2}
Let $\overline{w} \in W_0$.
\begin{enumerate}
    \item Let $x \in V_0$ and $w := t_x\overline{w} \in T(V_0)\rtimes W_0$. We have 
$$
\sL_{\Lambda_0}(w) = \frac{h}{2}|x|^2 - \mathrm{ht}(x).
$$
 In particular $\sL_{\Lambda_0}(w) = \sL_{\Lambda_0}(t_x)$.
\item
Let $q \in {{M}}$ and $w := t_q\overline{w} \in W$. For any $\sigma_j \in \Sigma$ we have
$$
\sL_{\Lambda_0}(\sigma_jw) = \sL_{\Lambda_0}(\sigma_jt_q) = \sL_{\Lambda_0}(t_{\varpi_j + w_{0,j}w_0(q)}).
$$
\end{enumerate}
\end{Prop}

\begin{proof}
The first part 
is a straightforward consequence of \Cref{atomic extended}. Indeed, the level $\ell$ of $\Lambda_0$ is $1$ and since $r \in L$ we have $(r~|~\Lambda_0) = 0$. Let us prove now that $\sL_{\Lambda_0}(\sigma_jw) = \sL_{\Lambda_0}(\sigma_jt_q)$. 
On the one hand we have $\sigma_jw = \sigma_jt_q\overline{w} = t_{\varpi_j}w_{0,j}w_0t_q\overline{w} = t_{\varpi_j + w_{0,j}w_0(q)}w_{0,j}w_0\overline{w}$. 
Since $\varpi_j + w_{0,j}w_0(q) \in L$ by \Cref{MvsL} and \Cref{lemma Garnier} (4), Point (1) ensures that
$\sL_{\Lambda_0}(\sigma_jw) = \sL_{\Lambda_0}(t_{\varpi_j + w_{0,j}w_0(q)})$. 
On the other hand we have $\sigma_jt_q =  t_{\varpi_j + w_{0,j}w_0(q)}w_{0,j}w_0$, and by the same argument as before we get $\sL_{\Lambda_0}(\sigma_jt_q) = \sL_{\Lambda_0}(t_{\varpi_j + w_{0,j}w_0(q)})$. This ends the proof.
\end{proof}

\subsection{Atomic length vs fundamental group}
\label{sec_ext_AL_A_C}

In this section, we establish an important property about the $\Lambda_0$-atomic length in its extended version. 
Let $W$ be the affine Weyl group of an affine Kac-Moody algebra, ${{M}}$ be the lattice introduced in \Cref{orbit theta} so that $W = T({{M}})\rtimes W_0$ with $W_0$ the finite part of $W$, and let $\Sigma$ be the associated fundamental group \Cref{sec_ext_AWG}. Let $\widehat{W} = \Sigma \ltimes W$ be the corresponding extended affine Weyl group. Recall that any element $\widehat{w} \in \widehat{W}$ decomposes uniquely as $\widehat{w} = \sigma w$ for some $\sigma \in \Sigma $ and some $w \in W$.

\medskip

It has been observed that the atomic length, when considered on the Weyl group $W$, enjoys some striking similarities with the usual Coxeter length $\ell$ \cite{CG2022}. 
Now, the Coxeter length $\ell$ is not affected while going to the extended affine Weyl group $\widehat{W}$, 
that is for any $w \in W$ and any $\sigma \in \Sigma$ one has 
\begin{equation}\label{ell not affected extended}
    \ell(w) = \ell(\sigma w) = \ell(w \sigma).
\end{equation}

This property is straightforward by definition. On the contrary, it is a priori not clear at all that a similar behaviour should hold for the $\Lambda$-atomic length. 
The main goal of this section is to prove that this is the case for $\La=\Lambda_0$,
as stated in the following theorem.

\bigskip

\begin{Th}\label{Theorem uniform extended atomic length}
Assume that we are not in type $A_{2n}^{(2)}$.
Let $\sigma \in \Sigma$ and $w \in W$. We have \begin{equation}
\sL_{\Lambda_0}(\sigma w) = \sL_{\Lambda_0}(w).
\end{equation}
\end{Th}

We delay the proof of this theorem until the end of this section, as it requires some preliminary results.

\begin{Lem}\label{most important lemma}
Assume that $a_0 = 1$. Let $r=\varpi_j+w_{0,j}w_0(q) \in L$ where $j\in J$ and $q\in {{M}}$. 
We have $\mathrm{ht}(w_{0,j}w_0(q)) = \mathrm{ht}(q) + h\Big(\varpi_j~|~w_{0,j}w_0(q)\Big)$.
\end{Lem}

\begin{proof}
     Since $w_0(\Delta_0) = -\Delta_0$,  there is only one element in $\Delta_0$, say $\alpha_{j^*}$, such that $w_0(\alpha_{j^*}) = -\alpha_j$.  Let $q =\sum\limits_{\alpha_k \in \Delta_0}a_k\alpha_k$ be the decomposition of $q$ in the simple basis. It follows that
\begin{align*}
    -a_{j^*} &  = \big( \varpi_j~|~w_0(q) \big) \\
            & = \big(w_{0,j}(\varpi_j)~|~w_{0,j}w_0(q) \big) \quad \text{by}\quad \Cref{invariance form ()} \\
            & = \big( \varpi_j~|~w_{0,j}w_0(q) \big) \quad \text{by} \quad \text{\Cref{lemma Garnier} (5)}
\end{align*}

Write now $X_j = w_{0,j}w_0(\Delta_0 \setminus \{\alpha_{j^*}\})$. By  \Cref{lemma Garnier} we easily have that $X_j \subset \Delta_0$ and we also know that $w_{0,j}w_0(\al_{j^*}) = - \theta$. Therefore
\begin{align}\label{decomp element important}
w_{0,j}w_0(q)
 = w_{0,j}w_0\Big(\sum\limits_{\alpha_k \in \Delta_0 \setminus \{\alpha_{j^*}\}}a_k\alpha_k\Big) -a_{j^*}\theta  
 = \sum\limits_{\beta_k \in X_j}a_k\beta_k -a_{j^*}\theta. 
\end{align}

Since $a_0 = 1$, by \Cref{root theta} we have that $h = 1+\text{ht}(\theta)$. It follows then that 
\begin{align*}
\text{ht}(w_{0,j}w_0(q)) &= \sum_{a_k \neq a_{j^*}} a_k - a_{j^*}\text{ht}(\theta) \\
&=  \sum_{k \in I} a_k  - a_{j^*} - a_{j^*}\text{ht}(\theta)\\
&= \text{ht}(q) - (1 + \text{ht}(\theta))a_{j^*} \\
& = \text{ht}(q) - ha_{j^*}.
\end{align*}

Since $\Big(\varpi_j~|~w_{0,j}w_0(q)\Big) = -a_{j^*}$, it follows that $\text{ht}(w_{0,j}w_0(q)) = \text{ht}(q) + h\Big(\varpi_j~|~w_{0,j}w_0(q)\Big)$.
\end{proof}

The following lemma is crucial. 
Note that we provide a case-by-case proof, and that it would be interesting to have a more uniform proof.

\begin{Lem}\label{Lem_ext} 
Assume that we are not in type $A_{2n}^{(2)}$.
Let $i \in J$. We have $\sL_{\Lambda_0}(t_{\varpi_i}) = 0$.
\end{Lem}

\begin{proof}
We only prove it in type $A_n^{(1)}$, the other cases being completely similar.
In this case, recall from \Cref{table} that the Coxeter number is $ h = n+1$ and  from \Cref{Sigma type A} that $J = \{1,\dots, n\}$.
Moreover, recall the formulas for the fundamental weights given in \Cref{fund_weights_A}.
Therefore, one has
\begin{align*}
\text{ht}(\varpi_i) = \frac{n+1-i}{n+1} \cdot \frac{i(i+1)}{2} +  \frac{i}{n+1}\cdot  \frac{(n-i)(n-i+1)}{2} 
					= \frac{i(n+1-i)}{2(n+1)}(i+1+n-i)
                     =  \frac{i(h-i)}{2},
\end{align*}
and
\begin{align*}
|\varpi_i|^2 & = i  \left(1-\frac{i}{n+1}\right)^2 + (n+1-i)\left(\frac{i}{n+1}\right)^2 
 = i\left(\frac{n+1-i}{n+1}\right)^2 + (n+1-i)\left(\frac{i}{n+1}\right)^2 \\
& = \frac{i(n+1-i)}{(n+1)^2}(n+1-i + i) 
 = \frac{i(n+1-i)}{n+1} \\
& =  \frac{i(h-i)}{h}.
\end{align*}

Now, by \Cref{analogue cor 8.2}, we have 
\begin{align*}
 \sL_{\Lambda_0}(t_{\varpi_i})  = \frac{h}{2}|\varpi_i|^2 - \text{ht}(\varpi_i)  
 = \frac{h}{2} \frac{i(h-i)}{h} -\frac{i(h-i)}{2} 
= 0.
\end{align*}
\end{proof}

\begin{Rem}
Excluding type $A_{2n}^{(2)}$ is necessary.
Indeed, for instance in type $A_4^{(2)}$, we have $J=\{2\}$, and 
$\varpi_2=\frac{1}{2}(e_1+e_2)=\frac{1}{2}(\al_1+\al_2)$.
Thus $\sL_{\La_0}(t_{\varpi_2}) = \frac{5}{2}|\varpi_2|^2-\h(\varpi_2)= -\frac{3}{8}$.
\end{Rem}

We are now ready to prove the  main result of this section.

\begin{proof}[Proof of \Cref{Theorem uniform extended atomic length}]
First of all, the assumption that we are not in type $A_{2n}^{(2)}$ ensures that we can apply \Cref{most important lemma}
and \Cref{Lem_ext}.
Let $\sigma=\sigma_j$ for some $j\in J$ and 
$w= t_q \overline{w}$ for some $q\in {{M}}$ and $\overline{w}\in W_0$.
Then we have
$$
\begin{array}{rcll}
\sL_{\La_0}(\sigma_j t_q \overline{w})
& = & \sL_{\La_0}(t_{\varpi_j+w_{0,j}w_0(q)}) & \text{by \Cref{analogue cor 8.2} (2)}
\\
& = & \sL_{\La_0}(t_{\varpi_j})+\sL_{\La_0}(t_{w_{0,j}w_0(q)})+ h\left( \varpi_j\mid w_{0,j}w_0(q)\right) & \text{by \Cref{decompo atomic D} and \Cref{defect term for Lambda0}}
\\
& = & \sL_{\La_0}(t_{w_{0,j}w_0(q)})+ h\left( \varpi_j\mid w_{0,j}w_0(q)\right) & \text{by  \Cref{Lem_ext}}
\\
& = & \frac{h}{2}|w_{0,j}w_0(q)|^2-\h(w_{0,j}w_0(q)) + h\left( \varpi_j\mid w_{0,j}w_0(q)\right) & \text{by \Cref{analogue cor 8.2} (1)}
\\
& = & \frac{h}{2}|q|^2-\h(w_{0,j}w_0(q)) + h\left( \varpi_j\mid w_{0,j}w_0(q)\right) & \text{since $w_{0,j}w_0$ is an isometry}
\\
& = & \frac{h}{2}|q|^2- (\h(q) + \h\left( \varpi_j\mid w_{0,j}w_0(q)\right)  ) + h\left( \varpi_j\mid w_{0,j}w_0(q)\right) & \text{by \Cref{most important lemma}}
\\
& = & \frac{h}{2}|q|^2- \h(q) &
\\
& = & \sL_{\La_0}(t_q)& \text{ by \Cref{thm_LA}}.
\end{array}
$$
\end{proof}

Note that we have proved along the way that the following diagram commutes 
(assuming we are not in type $A_{2n}^{(2)}$).
\begin{center}
\begin{tikzpicture} 
\node at (0,0) {$ \faktor{\widehat{W}}{W_0}$} ;
\node at (3,0) {$\Sigma W^0$} ;
\node at (3,-2) {$\mathbb{N}.$} ;
\node at (0,-2) {$L$} ;

\node at (0.5,0) (1)  {} ;
\node at (2.5, 0) (2) {} ; 
\draw [>=stealth,->] (1) to (2);
\node at (1.5,0.3) {$\widehat{\chi}$};

\node at (0,-0.3) (3)  {} ;
\node at (0, -1.8) (4) {} ; 
\draw [>=stealth,->] (3) to (4);
\node at (-0.7,-1) {$\widehat{f}$};

\node at (3,-0.3) (5)  {} ;
\node at (3, -1.8) (6) {} ; 
\draw [>=stealth,->] (5) to (6);
\node at (3.5,-1) {$\sL_{\Lambda_0}$};

\node at (0.5,-2) (7)  {} ;
\node at (2.5, -2) (8) {} ; 
\draw [>=stealth,->] (7) to (8);
\node at (1.5,-2.4) {$\widehat{\mu}$};
\end{tikzpicture}
\end{center}
where $\widehat{\mu}(\varpi_j+w_{0,j}w_0(q)) = \mu(w_{0,j}w_0(q)) + h\Big( \varpi_j~|~w_{0,j}w_0(q) \Big)$.
This shall be compared with \Cref{link Q et L}, and  shows that $\widehat{\mu}$ is the extended analogue of $\mu$.
Accordingly, we give the analogue of \Cref{definition B sans hat} for extended affine Grassmannian elements.

\begin{Def}\label{definition B hat}
Let $N \in \mathbb{N}$. We define the extended version of $\mathcal{B}(N)$ by
\begin{equation}\label{extended B(N)}
\widehat{\mathcal{B}}(N) := \{ q\in L \mid \sL_{\Lambda_0}(t_q) = N\}.
\end{equation}
\end{Def}

\medskip 

A special case of \Cref{Theorem uniform extended atomic length} (taking $w=t_q$)
is the following corollary, which will be useful later.

\begin{Cor}\label{AL_r_vs_q}
Let $j \in J$ and $r=\varpi_j+w_{0,j}w_0(q) \in L$. 
Then $\sL_{\Lambda_0}(t_r) = \sL_{\Lambda_0}(t_q)$.
\end{Cor}

\begin{Exa}[\Cref{Theorem uniform extended atomic length} in action] 
Let $W$ be the Weyl group of type $A_2^{(1)}$. Let $w=t_q\overline{w} \in W$. We have $\Sigma = \{1, \sigma_1, \sigma_2\}$. In this example we compute explicitly the $\Lambda_0$-atomic length of $\sigma_1 w$ and $\sigma_2 w$. Write $q = (q_1,q_2,q_3)$ with $q_1 + q_2 + q_3 =0$. 

\begin{itemize}
    \item A direct computation using \Cref{analogue cor 8.2} (1) gives that
$
\sL_{\Lambda_0}(w)=\sL_{\Lambda_0}(t_q) = 3(q_1^2 + q_2^2 + q_1q_2) - (2q_1 + q_2).
$ A more detailed computation is provided in \Cref{exa_Pell_A2}.
    \item One computes $\sL_{\Lambda_0}(\sigma_1 w) = \frac{3}{2}|w_{0,1}w_0(q)|^2 - \text{ht}(w_{0,1}w_0(q)) + 3\Big(\omega_1~|~w_{0,1}w_0(q)\Big) $. Moreover, a short computation gives $w_{0,1}w_0(q) = (q_3,q_1,q_2)$ and it follows, first that $|w_{0,1}w_0(q)|^2 = |q|^2$,  and second that $\text{ht}(w_{0,1}w_0(q)) = 2q_3 + q_1 = 2(-q_1 - q_2) + q_1 = -q_1 - 2q_2$. Finally, an easy computations gives that $\Big(\varpi_1~|~w_{0,1}w_0(q)\Big) = q_3 = -(q_1 + q_2)$.
    Therefore $\sL_{\Lambda_0}(\sigma_1 w) = 3(q_1^2 + q_2^2 + q_1q_2) + (q_1 + 2q_2) -3(q_1 + q_3) = 3(q_1^2 + q_2^2 + q_1q_2) - (2q_1 + q_2)$.

    \item Following the exact same process as in the previous point, we obtain that $\sL_{\Lambda_0}(\sigma_2 w) = 3(q_1^2 + q_2^2 + q_1q_2) - (q_2 - q_1)  -3q_1 = 3(q_1^2 + q_2^2 + q_1q_2) - (2q_1 + q_2)$. 
\end{itemize}
\end{Exa}

\subsection{Gaussian reduction and a quadratic Diophantine equation}\label{section gaussian equation gen}

Denote $q=q_1e_1+\ldots +q_{d}e_{d}\in {{M}}$,
where $n\leq d\leq n+2$ is the dimension of the ambient space.
Recall that since the rank of $A$ is $n$, $q$ is determined by $q_i$ for $1\leq i\leq n$.
Now, \Cref{thm_LAi} gives a polynomial formula of degree 2 in the variables
$q_1,\ldots, q_{n}$ for $\sL_{\Lambda_i}(t_q)$.
Indeed, $\frac{h}{2}|q|^2$ and $-\h(q)=-\langle q,\rho^\vee\rangle$ 
give a homogeneous part of degree $2$ and $1$ respectively.

\begin{Not}\label{def_pol_Q}
We denote by $Q_{\Lambda_i}\in\Q[X_1,\ldots, X_{n}]$ 
the corresponding polynomial. 
We write $Q_{\Lambda_i}(q)$ for the evaluation of $Q_{\La_i}$ at $(q_1,\ldots, q_n)$.
In particular, the case $i=0$ gives (by \Cref{thm_LA} and  \Cref{link Q et L})
\begin{equation}\label{formule Q et theo La}
    Q_{\Lambda_0}(q) = \mu(q) = \frac{h}{2}|q|^2-\text{ht}(q)= \sL_{\La_0}(t_q).
\end{equation}
\end{Not}

We now observe the following key phenomenon.
The right-hand side of \Cref{formule Q et theo La} 
is the sum of a quadratic form and a linear form in the variables $q_1,\ldots, q_n$ with integer coefficients.
Therefore, we can perform a Gaussian
 reduction (for example with respect to $q_1$, $q_2$ and so on) and we obtain an equivalent identity of the form
\begin{equation}
\label{eq_gen_0}
    a\sL_{\Lambda_0}(t_q) + b = \sum_{i=1}^n d_i\varphi_i(q_1,\ldots, q_{n})^2,
\end{equation}
where $a,b,d_i\in\Z$ and
where $\varphi_i :  \R^{n}\to\R$, $i=1,\ldots, n$ are linearly independent linear forms with integer coefficients.
We set accordingly
\begin{equation}\label{varphi general}
    \begin{array}{ccccc}\
\varphi & : & V_0 & \longrightarrow & \mathbb{R}^n\\
&   & v  & \longmapsto &  (\varphi_1(v),\ldots, \varphi_n(v)).
\end{array}
\end{equation}
\medskip

\newcommand{\bd}{\mathbf{d}}

Consider for $k\in\N$ the set
\begin{equation}\label{eq atomique}
    \cU^\Q_\bd (k) = \left\{ (x_1,\ldots, x_n)\in \Q^n \mid \sum_{i=1}^n d_i x_i^2 = k \right\}
\end{equation}
where $\bd=(d_1,\ldots, d_n)$,
and let $\cU_\bd(k)$ be the integral points of $ \cU^\Q_\bd(k)$.
In most of the rest of the paper, $\bd$ will be fixed 
and we will drop the subscript (but it will be useful in \Cref{exa_Pell_A2}). 
With this notation, we can reformulate \Cref{eq_gen_0} as follows.
\begin{Lem}
Let $q\in {{M}}$ and set $N=\sL_{\La_0}(t_q)$.
    Then $\varphi(q)\in\cU_{\bd}(aN+b)$.
\end{Lem}

In other terms, each lattice point $q$ (or equivalently, each affine Grassmannian element)
of $\La_0$-atomic length $N$ yields an integral solution $(\varphi_1(q),\ldots, \varphi_n(q))$ 
of the Diophantine equation
\begin{equation}
\label{eq_gen}
    d_1 x_1^2 + \cdots + d_n x_n^2 = aN + b.
\end{equation}
This phenomenon will be investigated in detail in \Cref{sec_sols_rank2}, \Cref{sec_pell_A3} and \Cref{sec_Pell_hyperoct}.

\medskip

In addition, the following proposition gives a way to construct new solutions of \Cref{eq_gen}
by considering extended affine Grassmannian elements.

\begin{Prop}\label{new_sols_ext}
Assume that we are not in type $A_{2n}^{(2)}$.
Let $i=1,\dots, n$, $q \in {{M}}$ and set $N=\sL_{\Lambda_0}(t_q)$. 
\begin{equation}
    \varphi(q)\in\mathcal{U}^{\Q}_\bd(aN+b) \quad \Longleftrightarrow \quad\varphi(\varpi_i+w_{0,i}w_0(q)) 
\in \mathcal{U}^{\Q}_\bd(aN+b).
\end{equation}
\end{Prop}

\begin{proof}
Denote $r=\varpi_j+w_{0,j}w_0(q)$.
We have 
$$\begin{array}{rcll}
Q_{\La_0}(r)
& = & \sL_{\La_0}(t_r) & \text{by \Cref{analogue cor 8.2} (1) and \Cref{def_pol_Q}}
\\
& = & \sL_{\La_0}(t_q) & \text{by \Cref{AL_r_vs_q}}
\\
& = & N. & 
\end{array}
$$
This yields
$d_1\varphi_1(r)^2+d_2\varphi_2(r)^2+\cdots + d_n\varphi_n(r)^2 =aN+b.$
Since the elements $\varphi_i(r)$ obviously belong to $\Q$, we conclude that $\varphi(r)\in\cU^{\Q}_{\bd}(aN+b)$.
\end{proof}

\medskip

\begin{Rem}\label{remark point Q to Z}
\Cref{new_sols_ext} enables us to construct new \textit{rational} solutions of \Cref{eq_gen},
but it is not clear whether these are actually \textit{integral}.
In \Cref{new_sols_A}, we will prove that these indeed lie in $\cU_{\bd}(aN+b)$ in type $A_n^{(1)}$. 
\end{Rem}


\section{Type $A_n^{(1)}$: core size and atomic length}

\label{sec_cores}

In this section, we focus on type $A_n^{(1)}$.
More precisely, let $W_0=\langle s_1,\cdots, s_n\rangle$ be the Weyl group of type $A_n$ 
and $W = \langle s_0, s_1, \cdots, s_n\rangle = T({{M}}) \rtimes W_0$ (see \Cref{table}).
In this setting, it will be convenient to use the $\Z$-basis $\cC=(\eps_1,\ldots, \eps_n)$ of ${{M}}$ (see \Cref{table}) defined by
\begin{equation}\label{basis_eps}
\eps_i = e_i-e_{n+1}.
\end{equation}
Note that if $q=q_1\eps_1 +\cdots + q_n\eps_n\in {{M}}$, then the $q_i$'s are also the coefficients of $q$
written in the basis $(e_1,\ldots, e_{n+1})$ of $\R^{n+1}$.

\medskip

We will now recall Lascoux's recursive construction of the $(n+1)$-cores via the action of
$W$ \cite{Lascoux2001}.
Then, we will use the quadratic formula for the atomic length of \Cref{thm_LA}, similar to that of \cite{GKS1990},
to establish the special Diophantine equation evoked in \Cref{section gaussian equation gen}.

\subsection{Core partitions via alcove geometry}
\label{sec_cores_alcoves}

Recall that we have bijections between the affine Grassmannian $W/W_0$, the set $W^0$, the lattice $M$, and alcoves of the fundamental chamber, see \Cref{sec_AG}.

\begin{Exa}
We illustrate these bijections in the case $n=2$ in \Cref{affine grass elements A2}.
There, $0$ corresponds to the color black, $1$ to blue and $2$ to red.
In the left figure, to shorten the notation, each alcove $\mathcal{A}_w$ is filled with the word 
$i_1i_2\dots i_n$, where $w=s_{i_1}s_{i_2}\dots s_{i_n}$ is a chosen reduced expression (recall that $i_1=0$ by \Cref{sec_AG}).
Recall from \Cref{rem_wall_crossing} that crossing an $i$-wall corresponds to adding the letter $i$ (on the right).
In the right figure, each alcove $\cA_w$ is filled with $q\in {{M}}$ where $w^{-1}=t_q\overline{w}$  with $\overline{w}\in W_0$

\begin{figure}[h!]
\centering
\begin{tikzpicture}[scale=1.7]
\footnotesize
\clip (0,0) to + (4.5,0) -- + (60:4.5) to (60:4) -- cycle;

\foreach \i in {0,...,4}
{
\foreach \j in {0,...,4}
{
\draw[purple!80] (3*\i,{\j*2*sqrt(3)/2}) -- + (0:1);
\draw[purple!80] (3/2+3*\i,{sqrt(3)/2+\j*2*sqrt(3)/2}) -- + (0:1);
\draw[purple!80] (1+3*\i,{\j*2*sqrt(3)/2}) -- +(60:1);
\draw[purple!80] (5/2+3*\i,{sqrt(3)/2+\j*2*sqrt(3)/2}) -- +(60:1);
\draw[purple!80] (1+3*\i,{\j*2*sqrt(3)/2}) -- +(-60:1);
\draw[purple!80] (5/2+3*\i,{sqrt(3)/2+\j*2*sqrt(3)/2}) -- +(-60:1);

\draw[teal] (-1+3*\i,{\j*2*sqrt(3)/2}) -- + (0:1);
\draw[teal] (1/2+3*\i,{sqrt(3)/2+\j*2*sqrt(3)/2}) -- + (0:1);
\draw[teal] (3*\i,{\j*2*sqrt(3)/2}) -- +(60:1);
\draw[teal] (3/2+3*\i,{sqrt(3)/2+\j*2*sqrt(3)/2}) -- +(60:1);
\draw[teal] (3/2+3*\i,{-sqrt(3)/2+\j*2*sqrt(3)/2}) -- +(-60:1);
\draw[teal] (3+3*\i,{\j*2*sqrt(3)/2}) -- +(-60:1);

\draw[black] (1+3*\i,{\j*2*sqrt(3)/2}) -- + (0:1);
\draw[black] (5/2+3*\i,{sqrt(3)/2+\j*2*sqrt(3)/2}) -- + (0:1);
\draw[black] (2+3*\i,{\j*2*sqrt(3)/2}) -- +(60:1);
\draw[black] (1/2+3*\i,{-sqrt(3)/2+\j*2*sqrt(3)/2}) -- +(60:1);
\draw[black] (1/2+3*\i,{sqrt(3)/2+\j*2*sqrt(3)/2}) -- +(-60:1);
\draw[black] (2+3*\i,{-sqrt(3)/2+sqrt(3)/2+\j*2*sqrt(3)/2}) -- +(-60:1);
}
}
\node[anchor = mid, scale=1.2] at ( 0.5, 0.3 ) {$e$};
\node[anchor = mid, scale=1] at ( 1, 0.55 ) {$0$};
\node[anchor = mid, scale=1] at ( 1, 1.1 ) {$01$};
\node[anchor = mid, scale=1] at ( 1.5, 0.3 ) {$02$};
\node[anchor = mid, scale=1] at ( 1.5, 1.5 ) {$012$};
\node[anchor = mid, scale=1] at ( 2, 0.6 ) {$021$};
\node[anchor = mid, scale=1] at ( 2, 1.1 ) {$0121$};
\node[anchor = mid, scale=1] at ( 1.5, 1.9 ) {$0120$};
\node[anchor = mid, scale=1] at ( 2, 2.4 ) {$01201$};
\node[anchor = mid, scale=1] at ( 2, 2.75 ) {$012012$};
\node[anchor = mid, scale=1] at ( 2.5, 0.3 ) {$0210$};
\node[anchor = mid, scale=1] at ( 2.5, 1.5 ) {$01210$};
\node[anchor = mid, scale=1] at ( 2.5, 1.9 ) {$012101$};
\node[anchor = mid, scale=1] at ( 3, 0.6 ) {$02102$};
\node[anchor = mid, scale=1] at ( 3, 1.1 ) {$021020$};
\node[anchor = mid, scale=1] at ( 3.5, 0.3 ) {$021021$};
\end{tikzpicture}
\qquad
\begin{tikzpicture}[scale=1.7]
\footnotesize
\clip (0,0) to + (4.5,0) -- + (60:4.5) to (60:4) -- cycle;

\foreach \i in {0,...,4}
{
\foreach \j in {0,...,4}
{
\draw[purple!80] (3*\i,{\j*2*sqrt(3)/2}) -- + (0:1);
\draw[purple!80] (3/2+3*\i,{sqrt(3)/2+\j*2*sqrt(3)/2}) -- + (0:1);
\draw[purple!80] (1+3*\i,{\j*2*sqrt(3)/2}) -- +(60:1);
\draw[purple!80] (5/2+3*\i,{sqrt(3)/2+\j*2*sqrt(3)/2}) -- +(60:1);
\draw[purple!80] (1+3*\i,{\j*2*sqrt(3)/2}) -- +(-60:1);
\draw[purple!80] (5/2+3*\i,{sqrt(3)/2+\j*2*sqrt(3)/2}) -- +(-60:1);

\draw[teal] (-1+3*\i,{\j*2*sqrt(3)/2}) -- + (0:1);
\draw[teal] (1/2+3*\i,{sqrt(3)/2+\j*2*sqrt(3)/2}) -- + (0:1);
\draw[teal] (3*\i,{\j*2*sqrt(3)/2}) -- +(60:1);
\draw[teal] (3/2+3*\i,{sqrt(3)/2+\j*2*sqrt(3)/2}) -- +(60:1);
\draw[teal] (3/2+3*\i,{-sqrt(3)/2+\j*2*sqrt(3)/2}) -- +(-60:1);
\draw[teal] (3+3*\i,{\j*2*sqrt(3)/2}) -- +(-60:1);

\draw[black] (1+3*\i,{\j*2*sqrt(3)/2}) -- + (0:1);
\draw[black] (5/2+3*\i,{sqrt(3)/2+\j*2*sqrt(3)/2}) -- + (0:1);
\draw[black] (2+3*\i,{\j*2*sqrt(3)/2}) -- +(60:1);
\draw[black] (1/2+3*\i,{-sqrt(3)/2+\j*2*sqrt(3)/2}) -- +(60:1);
\draw[black] (1/2+3*\i,{sqrt(3)/2+\j*2*sqrt(3)/2}) -- +(-60:1);
\draw[black] (2+3*\i,{-sqrt(3)/2+sqrt(3)/2+\j*2*sqrt(3)/2}) -- +(-60:1);
}
}
\node[anchor = mid, scale=1] at ( 0.5, 0.3 ) {$e$};
\node[anchor = mid, scale=0.8] at ( 1, 0.7 ) {${\alpha_1 + \alpha_2}$};
\node[anchor = mid, scale=0.9] at ( 1, 1.1 ) {${\alpha_2}$};
\node[anchor = mid, scale=0.9] at ( 1.5, 0.3 ) {${\alpha_1}$};
\node[anchor = mid, scale=0.9] at ( 1.5, 1.5 ) {${-\alpha_2}$};
\node[anchor = mid, scale=0.9] at ( 2, 0.6 ) {${-\alpha_1}$};
\node[anchor = mid, scale=0.9] at ( 2, 1.1 ) {${-\alpha_1-\alpha_2}$};
\node[anchor = mid, scale=0.8] at ( 1.5, 1.9 ) {${2\alpha_1+\alpha_2}$};
\node[anchor = mid, scale=0.8] at ( 2, 2.4 ) {${-\alpha_1+\alpha_2}$};
\node[anchor = mid, scale=0.9] at ( 2, 2.75 ) {${-\alpha_1-2\alpha_2}$};
\node[anchor = mid, scale=0.9] at ( 2.5, 0.12 ) {${\alpha_1+2\alpha_2}$};
\node[anchor = mid, scale=0.8] at ( 2.5, 1.6 ) {${2\alpha_1 + 2\alpha_2}$};
\node[anchor = mid, scale=0.9] at ( 2.5, 2 ) {${2\alpha_2}$};
\node[anchor = mid, scale=0.9] at ( 3, 0.7 ) {${\alpha_1-\alpha_2}$};
\node[anchor = mid, scale=0.9] at ( 3, 1.1 ) {${2\alpha_1}$};
\node[anchor = mid, scale=0.92] at ( 3.5, 0.12 ) {${-2\alpha_1-\alpha_2}$};
\end{tikzpicture}
\caption{Bijection between the alcoves of the fundamental chamber in type $A_2^{(1)}$
and affine Grassmannian elements (on the left), and the lattice $M$ (on the right).}
\label{affine grass elements A2}
\end{figure}

\end{Exa}

Recall that a partition is an $(n+1)$-core if it has no hook (or equivalently, no rim-hook) of size $(n+1)$ 
\cite[Section 2.7]{JamesKerber1984}.
We write
$\cC_{n+1}(N)$ for the set of $(n+1$)-cores of size $N$, and $\cC_{n+1}=\bigsqcup_{N\in\N} \cC_{n+1}(N).$
We will identify integer partitions with their Young diagram.
For $i\in\{0,\ldots, n\}$, an \textit{$i$-box} is a box $b$ of a Young diagram
such that the difference between the column and row index of $b$ equals $i$ modulo $n+1$.
One can show that $W$ acts transitively on $\cC_{n+1}$,
where the generator $s_i$
acts by adding or removing all $i$-boxes
(it is easy to see that these are either all addable or all removable) \cite{Lascoux2001}.
This yields a bijection
\begin{equation}
\label{bij_lascoux}
\begin{array}{ccc}
W^0 &  \overset{\sim}{\longrightarrow} & \cC_{n+1}
\\
w & \longmapsto & w\emptyset,
\end{array}
\end{equation}
which intertwines the Bruhat order and the Young order (inclusion of partitions),
see also \cite[Proposition 40]{LapointeMorse2005}.
This bijection is usually realised using the alcove geometry of $W$, by
filling the alcoves of the fundamental chamber by cores recursively as follows:
\begin{enumerate}
\item Put the empty partition $\emptyset$ inside the fundamental alcove.
\item Fill in all other alcoves of the fundamental chamber by $(n+1)$-cores recursively  by adding all possible $i$-boxes each time a new $i$-wall is crossed. 
\end{enumerate}

\begin{Exa}
The case $n=2$ is illustrated in \Cref{alcove_cores_A2}.
\begin{figure}[h!]
\centering
\begin{tikzpicture}[scale=2]
\clip (0,0) to + (4.5,0) -- + (60:4.5) to (60:4) -- cycle;

\foreach \i in {0,...,4}
{
\foreach \j in {0,...,4}
{
\draw[purple!80] (3*\i,{\j*2*sqrt(3)/2}) -- + (0:1);
\draw[purple!80] (3/2+3*\i,{sqrt(3)/2+\j*2*sqrt(3)/2}) -- + (0:1);
\draw[purple!80] (1+3*\i,{\j*2*sqrt(3)/2}) -- +(60:1);
\draw[purple!80] (5/2+3*\i,{sqrt(3)/2+\j*2*sqrt(3)/2}) -- +(60:1);
\draw[purple!80] (1+3*\i,{\j*2*sqrt(3)/2}) -- +(-60:1);
\draw[purple!80] (5/2+3*\i,{sqrt(3)/2+\j*2*sqrt(3)/2}) -- +(-60:1);

\draw[teal] (-1+3*\i,{\j*2*sqrt(3)/2}) -- + (0:1);
\draw[teal] (1/2+3*\i,{sqrt(3)/2+\j*2*sqrt(3)/2}) -- + (0:1);
\draw[teal] (3*\i,{\j*2*sqrt(3)/2}) -- +(60:1);
\draw[teal] (3/2+3*\i,{sqrt(3)/2+\j*2*sqrt(3)/2}) -- +(60:1);
\draw[teal] (3/2+3*\i,{-sqrt(3)/2+\j*2*sqrt(3)/2}) -- +(-60:1);
\draw[teal] (3+3*\i,{\j*2*sqrt(3)/2}) -- +(-60:1);

\draw[black] (1+3*\i,{\j*2*sqrt(3)/2}) -- + (0:1);
\draw[black] (5/2+3*\i,{sqrt(3)/2+\j*2*sqrt(3)/2}) -- + (0:1);
\draw[black] (2+3*\i,{\j*2*sqrt(3)/2}) -- +(60:1);
\draw[black] (1/2+3*\i,{-sqrt(3)/2+\j*2*sqrt(3)/2}) -- +(60:1);
\draw[black] (1/2+3*\i,{sqrt(3)/2+\j*2*sqrt(3)/2}) -- +(-60:1);
\draw[black] (2+3*\i,{-sqrt(3)/2+sqrt(3)/2+\j*2*sqrt(3)/2}) -- +(-60:1);
}
}
\node[anchor = mid, scale=1] at ( 0.5, 0.3 ) {$\emptyset$};
\node[anchor = mid, scale=0.3] at ( 1, 0.5 ) {$\yng(1)$};
\node[anchor = mid, scale=0.3] at ( 1, 1.2 ) {$\yng(2)$};
\node[anchor = mid, scale=0.3] at ( 1.5, 0.3 ) {$\yng(1,1)$};
\node[anchor = mid, scale=0.3] at ( 1.5, 1.5 ) {$\yng(3,1)$};
\node[anchor = mid, scale=0.3] at ( 2, 0.6 ) {$\yng(2,1,1)$};
\node[anchor = mid, scale=0.3] at ( 2, 1.2 ) {$\yng(3,1,1)$};
\node[anchor = mid, scale=0.3] at ( 1.5, 2 ) {$\yng(4,2)$};
\node[anchor = mid, scale=0.3] at ( 2, 2.4 ) {$\yng(5,3,1)$};
\node[anchor = mid, scale=0.3] at ( 2, 2.9 ) {$\yng(6,4,2)$};
\node[anchor = mid, scale=0.3] at ( 2.5, 0.4 ) {$\yng(2,2,1,1)$};
\node[anchor = mid, scale=0.3] at ( 2.5, 1.6 ) {$\yng(4,2,1,1)$};
\node[anchor = mid, scale=0.3] at ( 2.5, 2.1 ) {$\yng(5,3,1,1)$};
\node[anchor = mid, scale=0.3] at ( 3, 0.7 ) {$\yng(3,2,2,1,1)$};
\node[anchor = mid, scale=0.3] at ( 3, 1.3 ) {$\yng(4,2,2,1,1)$};
\node[anchor = mid, scale=0.3] at ( 3.5, 0.47 ) {$\yng(3,3,2,2,1,1)$};

\end{tikzpicture}

\caption{The fundamental chamber in type $A_2^{(1)}$ and the $3$-cores.}
\label{alcove_cores_A2}
\end{figure}

\end{Exa}

In \cite[Proposition 1.9]{LLMSSZ2014}, an explicit formula for the inverse map (that is, to construct an affine Grassmannian element from a core) is given, and relies on the notion of \textit{bounded partition}.

\medskip

Let $w\in W$ and  write $w=t_q \overline{w}$ with $\be\in {{M}}$, $\overline{w}\in W_0$. 
To compute $\be$ starting from $w$, one can use a reduced expression of $w$ and
bring it to the desired form.
Alternatively, there is a more combinatorial way to figure out $\be$,
starting from the core partition $\la=w\emptyset$, let us recall this briefly.
We refer again to
\cite{JamesKerber1984} for the definition of the following notions.
Since $\la$ is an $(n+1)$-core,
the \textit{$(n+1)$-quotient} of $\la$ is empty,
and the position of the last bead of each row in the corresponding \textit{$(n+1)$-abacus}
determines a vector of size $n+1$, sometimes called the
\textit{$(n+1)$-charge}, which coincides with $(q_1,\ldots,q_n, -(q_1+\cdots+ q_n))$
if $q=q_1\eps_1+\cdots q_n\eps_n$ (see \Cref{basis_eps}).
The fact that these two procedures give the same vector $\be$
is not hard to prove, see for instance \cite{Lascoux2001, LapointeMorse2005, LLMSSZ2014}.
We decide not to recall in detail these combinatorial constructions, referring to the above literature instead, 
but we still find it useful to give below an example of computations with $(n+1)$-cores.
This will be reinvested in \Cref{exa_C2}, \Cref{exa_D32_2} and  \Cref{exa_D43}.

\begin{Exa}\
\label{exa_cores_A}
\begin{enumerate}
\item One checks that the partition $(8,5,5,2,2,1)$ is a $5$-core of size $23$.
It is represented by the following abacus.
\begin{center}
\begin{tikzpicture}[scale=0.5, bb/.style={draw,circle,fill,minimum size=2.5mm,inner sep=0pt,outer sep=0pt}, wb/.style={draw,circle,fill,minimum size=0.5mm,inner sep=0pt,outer sep=0pt}]
exa
\draw [line width = 0.1mm] (15.5,0) -- (-14.5,0);

\node [wb] at (15,0) {};
\node [wb] at (14,0) {};
\node [wb] at (13,0) {};
\node [wb] at (12,0) {};
\node [wb] at (11,0) {};
\node [wb] at (10,0) {};
\node [wb] at (9,0) {};
\node [bb] at (8,0) {};
\node [wb] at (7,0) {};
\node [wb] at (6,0) {};
\node [wb] at (5,0) {};
\node [bb] at (4,0) {};
\node [bb] at (3,0) {};
\node [wb] at (2,0) {};
\node [wb] at (1,0) {};
\node [wb] at (0,0) {};
\node [bb] at (-1,0) {};
\node [bb] at (-2,0) {};
\node [wb] at (-3,0) {};
\node [bb] at (-4,0) {};
\node [wb] at (-5,0) {};
\node [bb] at (-6,0) {};
\node [bb] at (-7,0) {};
\node [bb] at (-8,0) {};
\node [bb] at (-9,0) {};
\node [bb] at (-10,0) {};
\node [bb] at (-11,0) {};
\node [bb] at (-12,0) {};
\node [bb] at (-13,0) {};
\node [bb] at (-14,0) {};

\node [scale = 0.7] at (15,-1) {15};
\node [scale = 0.7] at (14,-1) {14};
\node [scale = 0.7] at (13,-1) {13};
\node [scale = 0.7] at (12,-1) {12};
\node [scale = 0.7] at (11,-1) {11};
\node [scale = 0.7] at (10,-1) {10};
\node [scale = 0.7] at (9,-1) {9};
\node [scale = 0.7] at (8,-1) {8};
\node [scale = 0.7] at (7,-1) {7};
\node [scale = 0.7] at (6,-1) {6};
\node [scale = 0.7] at (5,-1) {5};
\node [scale = 0.7] at (4,-1) {4};
\node [scale = 0.7] at (3,-1) {3};
\node [scale = 0.7] at (2,-1) {2};
\node [scale = 0.7] at (1,-1) {1};
\node [scale = 0.7] at (0,-1) {0};
\node [scale = 0.7] at (-1,-1) {-1};
\node [scale = 0.7] at (-2,-1) {-2};
\node [scale = 0.7] at (-3,-1) {-3};
\node [scale = 0.7] at (-4,-1) {-4};
\node [scale = 0.7] at (-5,-1) {-5};
\node [scale = 0.7] at (-6,-1) {-6};
\node [scale = 0.7] at (-7,-1) {-7};
\node [scale = 0.7] at (-8,-1) {-8};
\node [scale = 0.7] at (-9,-1) {-9};
\node [scale = 0.7] at (-10,-1) {-10};
\node [scale = 0.7] at (-11,-1) {-11};
\node [scale = 0.7] at (-12,-1) {-12};
\node [scale = 0.7] at (-13,-1) {-13};
\node [scale = 0.7] at (-14,-1) {-14};
\end{tikzpicture} 
\end{center}
Its $5$-quotient is represented by the abacus
\begin{center}
\begin{tikzpicture}[scale=0.5, bb/.style={draw,circle,fill,minimum size=2.5mm,inner sep=0pt,outer sep=0pt}, wb/.style={draw,circle,fill,minimum size=0.5mm,inner sep=0pt,outer sep=0pt}]
\foreach \i in {0,1,2,3,4}
{
\draw [line width = 0.1mm] (-2.5,\i) -- (3.5,\i);
}
\node [wb] at (3,0) {};
\node [wb] at (2,0) {};
\node [wb] at (1,0) {};
\node [wb] at (0,0) {};
\node [wb] at (-1,0) {};
\node [bb] at (-2,0) {};

\node [wb] at (3,1) {};
\node [wb] at (2,1) {};
\node [bb] at (1,1) {};
\node [bb] at (0,1) {};
\node [bb] at (-1,1) {};
\node [bb] at (-2,1) {};

\node [wb] at (3,2) {};
\node [bb] at (2,2) {};
\node [bb] at (1,2) {};
\node [bb] at (0,2) {};
\node [bb] at (-1,2) {};
\node [bb] at (-2,2) {};

\node [wb] at (3,3) {};
\node [wb] at (2,3) {};
\node [wb] at (1,3) {};
\node [wb] at (0,3) {};
\node [bb] at (-1,3) {};
\node [bb] at (-2,3) {};

\node [wb] at (3,4) {};
\node [wb] at (2,4) {};
\node [wb] at (1,4) {};
\node [bb] at (0,4) {};
\node [bb] at (-1,4) {};
\node [bb] at (-2,4) {};

\node [scale = 0.7] at (3,-1) {3};
\node [scale = 0.7] at (2,-1) {2};
\node [scale = 0.7] at (1,-1) {1};
\node [scale = 0.7] at (0,-1) {0};
\node [scale = 0.7] at (-1,-1) {-1};
\node [scale = 0.7] at (-2,-1) {-2};
\end{tikzpicture} 
\end{center}
which in particular shows that $\la$ is a $5$-core (since in each row, there are no gaps).
We can read off the $5$-charge by looking at the index of the last bead of each row. 
We obtain the vector $(0,-1,2,1,-2)$.
\item Consider the $3$-core $(3,1)$. On the one hand, one can check that
its $3$-quotient has associated charge $(0,-1,1)$.
On the other hand, from the alcove construction of cores, we see that $(3,1) = w\emptyset$ with $w=s_2s_1s_0$.
Indeed, in \Cref{alcove_cores_A2}, we have crossed walls of color $0$ (black), then $1$ (blue), 
then $2$ (red) to reach the alcove containing the partition $(3,1)$.
Now, one checks that $w=s_2s_1s_0 = \tau_{-\al_2} s_1$,
hence $q=-\al_2 = -e_2+e_3$, and we recover the charge computed above.
\end{enumerate}
\end{Exa}

\subsection{Size of core partitions and link with the Garvan--Kim--Stanton formula}
\label{sec_sizeAL}

Using either a representation-theoretic \cite{CG2022} or a combinatorial \cite{STW2023} argument, 
we are ensured that the formula for the atomic length given in
\Cref{sec_AL} gives the size of the corresponding $(n+1)$-core. 
More precisely, if $w=t_q \overline{w}$ for $\beta\in {{M}}$ and $\overline{w}\in W_0$ is any element of $W$,
the size of the core $\la= w \emptyset$ is given by the formula

$$|\la| = \sL_{\La_0}(t_q) = \frac{n+1}{2}|\be|^2 - \h(q).$$
In order to compute $\sL_{\La_0}(t_q)$,
we can use the decomposition 
$q = q_1\eps_1 + \cdots + q_{n}\eps_{n}$
of $\be$ in the basis $\cC$ of ${{M}}$.
By the discussion in \Cref{rk_norm_1}, we have
$$|\be|^2 = \sum_{i=1}^{n+1} q_i^2 =
q_1^2+\cdots+q_n^2+(q_1+\cdots + q_n)^2 = 2\sum_{i=1}^nq_i^2 + 2\sum_{1\leq i <j\leq n} q_iq_j$$

since $q_{n+1}= -(q_1+\cdots + q_n)$, and
\begin{equation}\label{formule hauteur type A}
\h(q) = \sum_{i=1}^{n} (n-i+1)q_i.
\end{equation}

This yields the formula
\begin{equation}
\label{formula_AL_An}
\sL_{\La_0}(t_q) = (n+1)\left(\sum_{i=1}^{n} q_i^2 + \sum_{1\leq i<j\leq n} q_iq_j \right) - \sum_{i=1}^{n}(n+1-i)q_i.
\end{equation}

The quadratic expression on the right-hand side is the polynomial $Q_{\Lambda_0}(q)$ introduced in \Cref{def_pol_Q}. 

\begin{Rem}\label{rem_GKS}
Note that setting $q_{n+1}=-(q_1+\cdots +q_n)$, the above identity writes
$$|\la| =
(n+1)\left(\sum_{i=1}^{n} q_i^2 + \sum_{1\leq i<j\leq n} q_iq_j \right) + \sum_{i=1}^{n+1}(i-1)q_i,$$
and we recover precisely the formula established by Garvan, Kim and Stanton \cite[Section 2]{GKS1990}.
\end{Rem}

Recall the procedure introduced in \Cref{section gaussian equation gen},
which allows us to rewrite
Identity \Cref{formula_AL_An} by perfoming a Gaussian reduction.
Let us detail the computations in the following essential example.

\subsection{Quadratic Gauss reduction of the $\Lambda_0$-atomic length}
\label{sec_gauss_A}

Let $q = (q_1,\cdots,q_n,q_{n+1}) \in {{M}}$ given in the standard basis. Recall that $h = n+1$ is the Coxeter number of $A_n^{(1)}$. Let us denote, for any $i =1,\dots, n+1$

\begin{equation}
\left\{
\begin{array}{l}
t_i = 2i(i+1)   \\
S_i = \sum\limits_{k=1}^i k = \frac{i(i+1)}{2} \\
P_i(q) = \sum\limits_{k=i}^n q_k,\quad\text{with the convention} \quad P_{n+1}(q) =0.
\end{array}
\right.
\end{equation}

\medskip

On can check from \Cref{formula_AL_An} the following Lemma, which gives a quadratic Gauss reduction of the atomic length  with respect to the order $q_1, q_2,\dots, q_{n}$. 

\begin{Lem} We have 
    \begin{equation}\label{reduc gauss general}
    \sL_{\Lambda_0}(t_q) = \sum\limits_{i=1}^n \frac{1}{t_ih} \Big( (i+1)hq_i + hP_{i+1}(q) +    S_i - h \Big)^2 - \frac{1}h\sum\limits_{i=1}^n\frac{(S_i - h)^2}{t_i}.
\end{equation}
\end{Lem}

But now, we have the following formula, which gives a better expression of the last term of \Cref{reduc gauss general}.

\begin{Lem}\label{lemma magic number (h-1)h(h+1)/24}
We have
    \begin{equation}
       \sum\limits_{i=1}^n\frac{(S_i - h)^2}{t_i} = \frac{(h-1)h(h+1)}{24}.
    \end{equation}
\end{Lem}

\begin{proof}
Note first that $t_i = 4S_i$. We have
    \begin{align*}
        \frac{1}{h}\sum\limits_{i=1}^n\frac{(S_i - h)^2}{t_i} & = \frac{1}{h}\sum\limits_{i=1}^n\frac{S_i^2}{t_i} - 2\sum\limits_{i=1}^n\frac{S_i}{t_i} + h\sum\limits_{i=1}^n\frac{1}{t_i}    = \frac{1}{h}\sum\limits_{i=1}^n\frac{S_i^2}{4S_i} - 2\sum\limits_{i=1}^n\frac{S_i}{4S_i} + \frac{h}{2}\sum\limits_{i=1}^n\frac{1}{i(i+1)} \\ 
                                                             & = \frac{1}{4h}\sum\limits_{i=1}^nS_i - \frac{1}{2}\sum\limits_{i=1}^n 1 + \frac{h}{2}\sum\limits_{i=1}^n\Big(\frac{1}{i} - \frac{1}{i+1} \Big) 
                                                             = \frac{1}{4h}\sum\limits_{i=1}^nS_i - \frac{1}{2}(h-1) + \frac{h}{2}\Big(1 - \frac{1}{n+1}\Big)\\ 
                                                             & = \frac{1}{4h}\sum\limits_{i=1}^nS_i - \frac{1}{2}(h-1) + \frac{h}{2}\Big(1 - \frac{1}{h}\Big) = \frac{1}{8h}\Big(\sum\limits_{i=1}^n i^2 + \sum\limits_{i=1}^n i) \\
                                                             & = \frac{1}{8h}\frac{n(n+1)(2n+1)}{6} + \frac{1}{8h}\frac{n(n+1)}{2}  
                                                              = \frac{1}{8h}\frac{n(n+1)}{2} \Big(\frac{2n+1}{3} + 1 \Big) \\ 
                                                             & = \frac{1}{8h}\frac{n(n+1)}{2}\frac{2n+4}{3}  
                                                              = \frac{1}{8h}\frac{(h-1)h}{2}\frac{2(h+1)}{3} \\ 
                                                             & = \frac{(h-1)(h+1)}{24}.
    \end{align*}
\end{proof}

Let us define
\begin{equation}\label{def partial}
    \partial = \mathrm{lcm}\big(t_i\, ; \, i=1,\dots,n\big) \quad \text{and} \quad \partial_i = \frac{\partial}{t_i}.
\end{equation}

Note that $\partial\frac{(h-1)h(h+1)}{24} \in \mathbb{Z}$ by \Cref{lemma magic number (h-1)h(h+1)/24}. Therefore, if $q$ is a $(n+1)$-core of size $N$, Equation \Cref{eq_gen_0} is exactly given by (setting $a = \partial h$, $b = \partial\frac{(h-1)h(h+1)}{24}$ and $\partial_i = d_i$) 

\begin{equation}\label{gaussA}
    \partial h N + \partial\frac{(h-1)h(h+1)}{24}= \sum\limits_{i=1}^n \partial_i \Big( (i+1)hq_i + hP_{i+1}(q) +    S_i - h \Big)^2.
\end{equation}

\medskip

The affine map $\varphi$ introduced in \Cref{varphi general} is given by

\begin{equation}\label{map varphi An}
\begin{array}{ccccc}
\varphi & : & V_0 & \longrightarrow & \mathbb{R}^n \\
   & & q &\longmapsto & \Big( (i+1)hq_i + hP_{i+1}(q) +    S_i - h \Big)_{1\leq i\leq n}.
   \end{array}
\end{equation}

Let $\vec{\varphi}$ be the linear part of $\varphi$ so that 
\begin{equation}
\label{eq:matrice_varphi_An}
\varphi(q)=\vec{\varphi}(q)+u_0,
\end{equation}
where $u_0 = \Big(S_i - h \Big)_{1\leq i\leq n} \in \mathbb{Z}^n$.
The matrix of $\vec{\varphi}$ with respect to
the basis $\cC$ of $V_0$ and the standard basis of $\mathbb R^n$, is given by
\begin{equation}\label{matrix Q_0 gen}
    Q_0 = h
    \begin{pmatrix}
2&1&\cdots & 1
\\
0&3& 1 \quad \cdots & 1
\\
\vdots & & \ddots & \vdots
\\
0&0&\cdots &h
\end{pmatrix}.
\end{equation}

\medskip

In the following proposition we establish the property mentioned in \Cref{remark point Q to Z}.

\begin{Prop}\label{new_sols_A}
   We keep the notation of \Cref{def partial} and \Cref{gaussA} and we let $\bd = (\partial_i)_{1\leq i \leq n}$. Let $q \in {{M}}$. For any $i = 1,\dots, n$ we have 
\begin{equation}
    \varphi(q)\in\mathcal{U}_\bd(aN+b) \quad \Longleftrightarrow \quad\varphi(\varpi_i+w_{0,i}w_0(q)) 
\in \mathcal{U}_\bd(aN+b).
\end{equation}
\end{Prop}

\begin{proof}
     We have already seen in \Cref{new_sols_ext} that 
\begin{equation}
\varphi(q)\in\mathcal{U}^{\mathbb{Q}}_\bd(aN+b) \quad \Longleftrightarrow \quad\varphi(\varpi_i+w_{0,i}w_0(q)) 
\in \mathcal{U}^{\mathbb{Q}}_\bd(aN+b).
\end{equation}
Moreover
    \begin{equation*}
    \varphi(\varpi_i+w_{0,i}w_0(q))  = \vec{\varphi}(\varpi_i) + \vec{\varphi}(w_{0,i}w_0(q)) + u_0.
    \end{equation*}
    
    Since $u_0 \in \mathbb{Z}^n$ and $q \in {{M}}$ (implying that $q_i \in \mathbb{Z}$ and then also $w_{0,i}w_0(q)$), it follows that 
    \begin{equation*}
         \varphi(\varpi_i+w_{0,i}w_0(q)) \in \mathbb{Z}^n \quad \Longleftrightarrow \quad \vec{\varphi}(\varpi_i) \in \mathbb{Z}^n \quad \Longleftrightarrow \quad Q_0\varpi_i \in \mathbb{Z}^n.
    \end{equation*}
    However, by \Cref{fund_weights_A} we see that $\varpi_i \in \frac{1}{h}\mathbb{Z}^n$ in the standard basis and then, also in the basis $\mathcal{C}$, which by \Cref{matrix Q_0 gen} shows that $Q_0\varpi_i \in \mathbb{Z}^n$.
\end{proof}

\begin{Exa}\ 
\label{exa_Pell_A2}
In rank $2$ and $3$, we have in fact a slight simplification.
\begin{enumerate}
\item 
Let $n=2$.
From \Cref{gaussA}, we obtain
\begin{align*}
    36N + 12 
    &= 3\big(6q_1 + 3q_2 -2 \big)^2 + \big( 9q_2 \big)^2 \\
    &= 3\big(6q_1 + 3q_2 -2 \big)^2 + 9\big( 3q_2 \big)^2 
\end{align*}
Dividing both sides by $3$ yields
$$ 12N + 4 = (6q_1 + 3q_2-2)^2 + 3(3q_2)^2.$$
Therefore, we get two different Diophantine equations
$$ 3 x^2+ y^2 = 36N+12 \quad\mand\quad x^2+ 3y^2 = 12N+4,$$
which are in fact equivalent.
More precisely, if $(x,y)\in\cU_{\bd}(36N+12)$,
then $y^2=0\mod 3$ and thus $y=0\mod 3$, and we obtain a map
$$
\begin{array}{rrcl}
\al: &  \cU_{(3,1)}(36N+12) & \longmapsto &  \cU_{(1,3)}(12N+4)
\\
&(x,y) &\mapsto &(x,y/3),
\end{array}$$
which is a bijection by the same calculation as above.
By \Cref{new_sols_A}, we have a map
$$\begin{array}{rrcl}
\varphi: &  \widehat \cB(N) & \longmapsto &  \cU_{(3,1)}(36N+12)
\\
& (q_1,q_2)&  \longmapsto& (6q_1 + 3q_2-2, 9q_2).
\end{array}
$$
Composing with $\al$, we obtain the map
$$
\begin{array}{rrcl}
\al\circ\varphi: &  \widehat \cB(N) & \longmapsto &  \cU_{(1,3)} (12N+4)
\\
& (q_1,q_2)&  \longmapsto & (6q_1 + 3q_2-2, 3q_2).
\end{array}
$$
In \Cref{sec_Pell_A2}, we will consider instead the simpler equation
$x^2+ 3y^2 = 12N+4$,
and for simplicity, we will actually rename $\al\circ \varphi$ to $\varphi$.
\item Let $n=3$. 
From \Cref{gaussA}, we obtain
\begin{align*}
    96N + 60 
    &= 6\big(8q_1 + 4q_2 + 4q_3 -3 \big)^2 + 2\big(12q_2 + 4q_3 - 1\big)^2 + \big(16q_3 + 2\big)^2 \\
    & = 6\big(8q_1 + 4q_2 + 4q_3 -3 \big)^2 + 2\big(12q_2 + 4q_3 - 1\big)^2 + 4\big(8q_3 + 1\big)^2.
\end{align*}
Dividing both sides by $2$ yields
$$48N +30 = (12q_2 + 4q_3 - 1)^2 + 2(8q_3+1)^2 + 3(8q_1 + 4q_2 + 4q_3 - 3)^2.$$
Therefore, we get two different Diophantine equations
$$ 6 x^2+ 2 y^2 + z^2 = 96N+60 \quad\mand\quad x^2+ 2y^2 + 3z^2 = 48N+30,$$
which are again equivalent.
More precisely, if $(x,y,z)\in\cU_{(6,2,1)} (96N+60)$,
then $z^2=0\mod 2$ and thus $y=0\mod 2$, and we obtain a map
$$
\begin{array}{rrcl}
\al: &  \cU_{(6,2,1)}(96N+60) & \longmapsto &  \cU_{(1,2,3)}(48 N+ 30)
\\
&(x,y,z) &\mapsto &(y, x, z/2),
\end{array}$$
which is a bijection by the same calculation as above.
By \Cref{new_sols_A}, we have a map
$$\begin{array}{rrcl}
\varphi: &  \widehat \cB(N) & \longmapsto &  \cU_{(6,2,1)}(96N+60)
\\
& (q_1,q_2)&  \longmapsto& (8q_1 + 4q_2 + 4q_3 -3 , 12q_2 + 4q_3 - 1, 16q_3 + 2).
\end{array}
$$
Composing with $\al$, we obtain the map
$$
\begin{array}{rrcl}
\al\circ\varphi: &  \widehat \cB(N) & \longmapsto &  \cU_{(1,2,3)} (48N+30)
\\
& (q_1,q_2)&  \longmapsto & (12q_2 + 4q_3 - 1, 8q_1 + 4q_2 + 4q_3 -3 , 8q_3 + 1).
\end{array}
$$
In \Cref{sec_pell_A3}, we will this time consider instead the simpler equation
$x^2+ 2y^2 +3z^2= 48N+30$,
and we will again rename $\al\circ \varphi$ to $\varphi$.

\end{enumerate}
\end{Exa}


\section{Combinatorial models for generalised cores}
\label{sec_comb_models}

By the results of the previous section, it is natural to generalise the notion of core partitions by considering
affine Grassmannian elements,
or equivalently the elements of the lattice $M$.This is the point of view of \cite{STW2023}, and we shall use this approach in this paper.

\begin{Rem}
There is an even broader generalisation of cores.
Let $\fg$ denote the affine Kac-Moody algebra associated to 
the chosen generalised Cartan matrix.
It is well-known that we have a bijection
$$
\cO(\La_0)
\overset{1:1}{\longleftrightarrow}
W/W_0,$$
where $\cO(\La_0)$ is the orbit,  under the action of $W$, of the highest weight vertex in the $\fg$-crystal of the basic $\fg$-representation $V(\La_0)$.
Therefore, it is natural to consider a similar construction for an arbitrary dominant weight $\La$.
This is the approach of \cite{CG2022}, and also works for Kac-Moody algebras of finite type.
However, we will not need this degree of generality in 
the rest of this paper.
\end{Rem}

Now, we want to have combinatorial models for these generalised cores.
In classical untwisted types and in type $G_2^{(1)}$,
this has been achieved in \cite[Section 5]{STW2023},
see also \cite{HanusaJones2012, CottonWilliams2021},
where $M$ is put in bijection with a certain subset of partitions (often self-conjugate
with additional symmetry properties)
Moreover, Stucky, Thiel and Williams \cite{STW2023} show that
the statistic $\sL_{\Lambda_0}$ (whose formula is given in \Cref{thm_LA}) can be computed by counting boxes of the corresponding partition,
but using an appropriate weighting for boxes of a given residue,
hence generalising the formula of \Cref{rem_GKS}.
For the remaining types, namely the twisted ones, similar
results have been established by Lecouvey and Wahiche \cite{LW2024}.
We quickly recall the different results in all classical (untwisted and twisted) types and
in types $G_2^{(1)}$ and $D_4^{(3)}$.

\medskip

For $d \in\Z_{\geq 2}$ 
and $N\geq 0$, 
denote $\SC_d(N)$ the set of self-conjugate $d$-cores of size $N$,
and $\SC_d=\bigsqcup_{N\in\N}\SC_d(N)$.
For $\la\in\cC_d$ and for all $0\leq i\leq d-1$, we denote by $|\la|_i$ the number of
$i$-boxes of $\la$, so that $|\la| = \sum_{i=0}^{d-1}|\la|_i$.

\subsection{Type $C_n^{(1)}$}
\label{sec_cores_C}

let $W$ be the Weyl group of type $C_n^{(1)}$, see \Cref{table}.
Using a folding of the type $A_{2n-1}^{(1)}$ root system,
see \cite[Theorem 5.2]{STW2023} or \cite[Section 4.2]{LW2024} for details,
it is established that
\begin{equation}
\label{cores_C}
W/W_0 \overset{\sim}{\longleftrightarrow} \SC_{2n}.
\end{equation}
Moreover, for $q = q_1\sqrt{2}e_1 + \cdots + q_n\sqrt{2}e_n \in {{M}}$
with $q_i\in\Z$, we have
$$\sL_{\La_0}(t_q) = |\la|,$$ 
where $\la$ is the self-conjugate $2n$-core corresponding to $\be$ under Bijection \Cref{cores_C}.

\subsection{Type $D_{n+1}^{(2)}$}
\label{sec_cores_Dt}

let $W$ be the Weyl group of type $D_{n+1}^{(2)}$, see \Cref{table},
which is the transpose of type  $C_n^{(1)}$.
Using a similar root system folding,
see  \cite[Section 4.3]{LW2024} for details,
we have
\begin{equation}
\label{cores_Dt}
W/W_0 \overset{\sim}{\longleftrightarrow} \SC_{2n}.
\end{equation}
Moreover, for $q = q_1\sqrt{2}e_1 + \cdots + q_n\sqrt{2}e_n \in {{M}}$
with $q_i\in\Z$, we have
$$\sL_{\La_0}(t_q) = \frac{1}{2}\left(|\la| - |\la|_0 + |\la|_n\right),$$ 
where $\la$ is the self-conjugate $2n$-core corresponding to $\be$ under Bijection \Cref{cores_Dt}.

\medskip

In this case, there is an alternative partition model for which the size is exactly given by the
atomic length, see \cite{LW2024}. Since we will need it in \Cref{Pell_D32}, we quickly recall this construction here.
A \textit{bar partition} is a partition with only distinct parts (these are also called \textit{distinct} or \textit{$2$-regular} partitions). The combinatorics of bar partitions and their applications in representation
have been extensively studied in the literature, see for instance the seminal work \cite{Morris1965}.
To each bar partition $\la$, we can associate a \textit{double distinct} partition $\widetilde{\la}$,
obtained by shifting row $k$ of the Young diagram of $\la$ $k$ steps to the right,
and completing by the transpose of this shifted diagram, see \cite[Chapter III]{Macdonald1998} for details.

\newcommand\gry{\Yfillcolour{black!20}}
\newcommand\white{\Yfillcolour{white}}
\begin{Exa}
The partition $(4,2,1)$ is a bar partition of size $7$.
Its Young diagram is 
\Yboxdim{6pt}
$$
\yng(4,2,1),
$$
and its corresponding double distinct partition is
$$
\young(
!\gry<>!\white<><><><>,!\gry<><>!\white<><>,!\gry<><><>!\white<>,!\gry<>).
$$
\end{Exa}

\newcommand{\cD}{\mathcal{D}}
The bar partition $\la$ is called an \textit{$(2n+2)$-core} if $\widetilde{\la}$ is an $(2n+2)$-core.
Let us denote
by $\cD_{2n+2}(N)$ the set of $(2n+2)$-core bar partitions of size $N$ with no part equal to $(n+1)$,
and $\cD_{2n+2} = \sqcup_{N\in\N}\cD_{2n+2}(N)$.
As showed in \cite{LW2024}, there is a bijection
\begin{equation}
\label{cores_Dt_alt}
W/W_0 \overset{\sim}{\longleftrightarrow} \cD_{2n+2}
\end{equation}
such that
we have
$$\sL_{\La_0}(t_q) = |\la|$$ 
where $\la$ is the bar $(2n+2)$-core corresponding to $\be$ under Bijection \Cref{cores_Dt_alt}.
In fact, this bijection is explicit, and is essentially the same as the bijection in type $A_n^{(1)}$ recalled in \Cref{sec_cores_alcoves}.
More precisely, starting with $q=q_1\sqrt{2}e_1+ \cdots + q_n\sqrt{2}e_n\in {{M}}$,
we first associate the vector 
$\widetilde{q}=(0,q_1,\ldots, q_n,0,-q_n, \ldots ,-q_1)$.
By the procedure illustrated in \Cref{exa_cores_A},
this element $\widetilde{q}$ defines a $(2n+2)$-core partition,
which turns out to lie in $\cD_{2n+2}$:
this is the partition $\widetilde{\la}$.
We will give an illustration of this construction later, namely in \Cref{exa_D32_2}.

\begin{Exa}\
\label{exa_D32_1}
\begin{enumerate}
\item Take $n=2$, so that $2n+2=6$. 
Let us compute $\cD_6(N)$ for the first values of $N$.
$$
\Yboxdim{6pt}
\begin{array}{@{}l@{\hskip 20pt} @{}l@{}}
\hline 
\text{size $N$} & \text{Elements of }\cD_6(N)
\\
\hline
0 & 
\emptyset
\\
1 &
\young(<>)
\\
2 &
\young(<><>)
\\
3 &
\young(<><>,<>)
\\
4 &
\young(<><><><>)
\\
5 &
\young(<><><><>,<>), 
\young(<><><><><>)
\\
6 &
\text{none}
\\
\hline
\end{array}
$$
\item One can show that 
$\cD_6(35)=\left\{(17,11,5,2), (16,10,5,4) , (13,10,7,4,1) \right\}.$
\end{enumerate}
\end{Exa}

\subsection{Type $A_{2n}^{(2)}$}
\label{sec_cores_Ateven}

let $W$ be the Weyl group of type $A_{2n}^{(2)}$, see \Cref{table}.
Using again a folding of the type $A_{2n-1}^{(1)}$ root system,
see \cite[Section 4.4]{LW2024} for details,
it is established that we once again have
\begin{equation}
\label{cores_Ateven}
W/W_0 \overset{\sim}{\longleftrightarrow} \SC_{2n}.
\end{equation}
Moreover, for $q = q_1e_1 + \cdots + q_ne_n \in {{M}}$
with $q_i\in\Z$, we have
$$\sL_{\La_0}(t_q) =  |\la| + |\la|_0,$$ 
where $\la$ is the self-conjugate $2n$-core corresponding to $\be$ under Bijection \Cref{cores_Ateven}.

\subsection{Type $B_n^{(1)}$}
\label{sec_cores_B}

let $W$ be the Weyl group of type $B_{n}^{(1)}$, see \Cref{table}.
This time, one needs to first consider an embedding into
the type $D_{n+1}^{(2)}$ root system, and then use the results of \Cref{sec_cores_Dt}.
We will need the set
$$\SC_{2n}^+=\left\{ \la\in\SC_{2n} \mid \text{ $\la$ has an even number of diagonal boxes } \right\}.$$
Then we have, by \cite[Theorem 5.3]{STW2023} or \cite[Section 4.6]{LW2024},
\begin{equation}
\label{cores_B}
W/W_0 \overset{\sim}{\longleftrightarrow} \SC_{2n}^+.
\end{equation}
Moreover, for $q = q_1e_1 + \cdots + q_ne_n \in {{M}}$
with $q_i\in\Z$, we have
$$\sL_{\La_0}(t_q) =  \frac{1}{2}\left(|\la| - |\la|_0 + |\la|_n \right),$$ 
where $\la\in\SC_{2n}^+$ corresponds to $\be$ under Bijection \Cref{cores_B}.

\subsection{Type $A_{2n-1}^{(2)}$}
\label{sec_cores_Atodd}

let $W$ be the Weyl group of type $A_{2n-1}^{(2)}$, see \Cref{table}.
This time, one needs to first consider an embedding into
the type $A_{2n}^{(2)}$ root system, and then use the results of \Cref{sec_cores_Ateven}.
By \cite[Section 4.7]{LW2024}, we have a bijection
\begin{equation}
\label{cores_Atodd}
W/W_0 \overset{\sim}{\longleftrightarrow} \SC_{2n}^+.
\end{equation}
Moreover, for $q = q_1e_1 + \cdots + q_ne_n \in {{M}}$
with $q_i\in\Z$, we have
$$\sL_{\La_0}(t_q) =  \frac{1}{2}\left(|\la| - |\la|_0 \right),$$ 
where $\la\in\SC_{2n}^+$ corresponds to $\be$ under Bijection \Cref{cores_Atodd}.

\subsection{Type $D_n^{(1)}$}
\label{sec_cores_D}

let $W$ be the Weyl group of type $D_{n}^{(1)}$, see \Cref{table}.
One considers again an embedding into
the type $D_{n+1}^{(2)}$ root system.
We again get the following bijection, see \cite[Theorem 5.4]{STW2023} or \cite[Section 4.8]{LW2024},
\begin{equation}
\label{cores_D}
W/W_0 \overset{\sim}{\longleftrightarrow} \SC_{2n}^+.
\end{equation}
Moreover, for $q = q_1\sqrt{2}e_1 + \cdots + q_n\sqrt{2}e_n \in {{M}}$
with $q_i\in\Z$, we have
$$\sL_{\La_0}(t_q) =  \frac{1}{2}\left(|\la| - |\la|_0 - |\la|_n \right),$$ 
where $\la\in\SC_{2n}^+$ corresponds to $\be$ under Bijection \Cref{cores_D}.

\subsection{Type $G_2^{(1)}$}
\label{sec_cores_G}

let $W$ be the Weyl group of type $D_{n}^{(1)}$, see \Cref{table}.
By a more complicated construction, one obtains \cite[Section 4.11]{LW2024},
\begin{equation}
\label{cores_G21}
W/W_0 \overset{\sim}{\longleftrightarrow} \SC_{6}^\flat,
\end{equation}
where $\SC_{6}^\flat$ is the subset of self-conjugate $6$-cores
whose corresponding $6$-charge is of the form $(c_1,c_2,c_1-c_2,c_2-c_1,-c_2,-c_1)$.
Moreover, for $q = q_1e_1 + q_2e_2 + q_3e_3 \in {{M}}$
with $q_i\in\Z$ and $q_3=-q_1-q_2$, we have
$$\sL_{\La_0}(t_q) =  \frac{1}{2}\left(|\la| - |\la|_0 + |\la|_3 \right),$$ 
where $\la\in\SC_{6}^\flat$ corresponds to $\be$ under Bijection \Cref{cores_G21}.

\medskip

Alternatively, note that the $3$-cores also provide a combinatorial model for the affine Grassmannian elements 
in this type, see \cite[Section 5.4]{STW2023} and \cite{CottonWilliams2021}.

\subsection{Type $D_4^{(3)}$}
\label{sec_cores_D43}

let $W$ be the Weyl group of type $D_{n}^{(1)}$, see \Cref{table}.
By a similar construction as for $G_2^{(1)}$, one obtains \cite[Section 4.10]{LW2024},
\begin{equation}
\label{cores_D43}
W/W_0 \overset{\sim}{\longleftrightarrow} \SC_{6}^\flat,
\end{equation}
where $\SC_{6}^\flat$ is the subset of self-conjugate $6$-cores
whose corresponding $6$-charge is of the form $(c_1,c_2,c_1-c_2,c_2-c_1,-c_2,-c_1)$.
Moreover, for $q = q_1e_1 + q_2e_2 + q_3e_3 \in {{M}}$
with $q_i\in\Z$ and $q_3=-q_1-q_2$, we have
$$\sL_{\La_0}(t_q) =  \frac{1}{2}\left(|\la| - |\la|_0 - |\la|_3 \right),$$ 
where $\la\in\SC_{6}^\flat$ corresponds to $\be$ under Bijection \Cref{cores_D43}.

\medskip

In this case, we will also use the alternative partition model of \cite[Section 4.10]{LW2024},
which we denote here $\cD_{4}^\flat(N)$, in which the size of partitions is given by the atomic length.
In other terms, we have
\begin{equation}
\label{cores_D43_alt}
W/W_0 \overset{\sim}{\longleftrightarrow} \cD_{4}^\flat,
\end{equation}
and,
$$\sL_{\La_0}(t_q) = |\la|$$ 
where $\la$ is the partition corresponding to $\be$ under Bijection \Cref{cores_D43_alt}.
If we call  \textit{$i$-part} a non-zero part of $\la$ which is congruent to $i$ modulo $4$,
\cite{LW2024} show that the set $\cD_{4}^\flat$
consists of bar partitions $\la$ of size $N$ such that
$\widetilde{\la}$ has no $4$-hook below the diagonal\footnote{
Note that this does not quite mean that $\widetilde{\la}$ is a $4$-core, as
there might be boxes above or on the diagonal with hook-length $4$.
},
and with an extra technical condition
relating the number of $i$-parts of $\la$ (for $i=0,\pm,2$), see Formula \Cref{formule_i-parts} below.
Instead of recalling this relationship precisely here, 
we prefer to give the following construction, from which the exact formula 
can be easily derived.

\medskip

The combinatorial recipe to construct $\la\in\cD_4^\flat$ from $q\in {{M}}$
is as follows.
Write $q=q_1e_1+q_2e_2+q_3e_3\in {{M}}$. Recall that $q_3=-q_1-q_2$,
so the knowledge of $(q_1,q_2)$ is enough.
Then, define $m_0,m_{\pm1}, m_2$ by
\begin{equation}
\label{formules_david}
m_2 = |q_2|
\quad,\quad\quad
(m_1,m_{-1})=
\left\{\begin{array}{ll} 
(|q_2|,0) & \text{ if } q_2\leq 0 
\\
(0,|q_2|) & \text{ if } q_2 > 0 
\end{array}
\right.
\mand\quad
m_0=
\left\{\begin{array}{ll}
q_1+q_2 & \text{  if } q_1+q_2\geq 0 
\\
-(q_1+q_2)-1 & \text{  if } q_1+q_2 < 0.
\end{array}
\right.
\end{equation}
Then define $\la$ to be the only partition such that
\begin{itemize}
    \item $m_i$ is the number of $i$-parts of $\la$,
    \item $k$ is a part of $\la$ $\Rightarrow$ $k-4$ is a part of $\la$ (if $k-4>0$).
\end{itemize}
Note that the relation defining $m_0$ induces a relation between the $i$-parts of $\la$ of the form
\begin{equation}
\label{formule_i-parts}
m_0=\pm m_{\pm1} \pm m_2 -\delta \text{ \quad with }
\delta\in\{0,1\}.
\end{equation}

\begin{Exa}\
\label{exa_cores_D43}
\begin{enumerate}
\item
Let $q=-3e_1+e_2+e_3$, so that $(q_1,q_2)=(-3,1)$.
Formulas \Cref{formules_david} yield
\begin{itemize}
    \item $m_2=3$, so $\la$ contains the parts $2,6,10$,
    \item $m_{-1}=1$, so $\la$ contains the part $3$,
    \item $m_0=1$, so $\la$ contains the part $4$.
\end{itemize}
\Yboxdim{6pt}
We get $\la=(10,6,4,3,2)=\yng(10,6,4,3,2)$.
\item 
Let us compute $\cD_4^\flat(N)$ for the first values of $N$.
$$
\Yboxdim{6pt}
\begin{array}{@{}l@{\hskip 20pt} @{}l@{}}
\hline 
\text{size $N$} & \text{Elements of }\cD_4^\flat(N)
\\
\hline
0 & 
\emptyset
\\
1 &
\young(<>)
\\
2 &
\young(<><>)
\\
3 &
\young(<><>,<>)
\\
4 &
\text{none}
\\
5 &
\young(<><><>,<><>)
\\
6 &
\young(<><><><>,<><>)
\\
7 &
\young(<><><><>,<><><>),\young(<><><><>,<><>,<>)
\\
8 &
\young(<><><><><>,<><>,<>)
\\
9 &
\text{none}
\\
10 &
\young(<><><><><>,<><>,<>)
\\
\hline
\end{array}
$$
\end{enumerate}
\end{Exa}


\section{Pell-type equations and group actions}

\label{sec_pell}

In this short section, we investigate some $G$-set structures 
on the solution sets of certain Pell-like Diophantine equations,
where $G$ is a well-chosen finite group.
We will see in \Cref{sec_param} that these equations
will appear naturally in the study of generalised cores.

\subsection{The Diophantine equation $x^2 +y^2 = k$}
\label{sec_x2y2k}

Let $k\geq 2$ be an integer. We consider the equation
\begin{equation}
\label{eq_1}
x^2+y^2 =k
\end{equation}
in the variables $x,y$.
We denote by $\cU(k)$ the set of integer solutions of \Cref{eq_1}.

\medskip

First, we recall some generalities about \Cref{eq_1}. For that, we introduce 
some notation.
We denote by $\mathcal P_1(k)$ and $\mathcal P_3(k)$ the sets of prime
factors of $k$ with residue $1$ and $3$ modulo $4$.
Furthermore, we identify $\Z^2$ with the ring of Gaussian integers $\Z[i]$ by
$(x,y)\in\Z\mapsto x+iy\in\Z[i]$.
By abuse of notation, we will continue to write $\mathcal U(k)$ for
the set of elements $x+iy\in\Z[i]$ with $(x,y)$ solution of \Cref{eq_1}.

\medskip

We also recall that $\Z[i]$ is a principal ideal domain with group of
units $(\Z[i])^{\times}=\{1,-1,i,-i\}$, whose irreducible elements (up to
association) are described as follows.
\begin{itemize}
\item The element $(1+i)$.
\item The set of prime integers with residue $3$ modulo $4$.
\item The set of elements of the form $a+ib$ and $a-ib$ whose
product is a prime integer with residue $1$ modulo $4$.
\end{itemize}

Let $p$ be a prime. We denote by $\nu_p(k)$ the $p$-valuation of $k$, that is $k = p^{\nu_p(k)}q$ with $p \nmid q$.
Let $\alpha\in\N$ and $m$ be an odd integer such that $k=2^{\alpha}m$.
By the famous \emph{sum of two squares theorem}, we know that
\Cref{eq_1} has solutions if and only if $\nu_p(k)$ is even for
all $p\in\mathcal P_3(k)$. In particular, $\mathcal U(k)$ is empty if
and only if $\mathcal U(m)$ is.

\medskip

From now on, we assume that Equation \Cref{eq_1} has integer solutions.
Then the product (with multiplicity) of the prime factors with residue $3$
modulo $4$ is a square denoted by $c^2$.

\begin{Lem} 
\label{lembij}
The map
$$\Phi:\mathcal U(m) \longrightarrow \mathcal U(k),\ z\longmapsto
(1+i)^\alpha cz,$$
is a bijection.
\end{Lem}

\begin{proof}
First, we remark that if $z'= (1+i)^{\alpha} cz$ with
$z\in\mathcal U(m)$ (that is $z\overline z=m$), then
$$z'\overline z'=(1+i)^{\alpha}(1-i)^{\alpha}c^2z\overline
z=2^{\alpha}c^2m=k.$$
Hence, $\Phi$ is well-defined.
Now, assume that $z'\in\mathcal U(k)$. 
For any $p\in\mathcal P_1(k)$, we fix $z_p\in \Z[j]$ such that
$p=z_p\overline z_p$.
Since $(1-i)=-i(1+i)$, the $\Z[i]$-decomposition of $k$ into irreducible elements
is
\begin{equation}
\label{eq:dec}
k=(-i)^{\alpha}(1+i)^{2\alpha}\prod_{p\in\mathcal
P_3}p^{\nu_p(k)}\prod_{p\in\mathcal P_1}z_p^{\nu_p(k)}\overline
z_p^{\nu_p(k)}.
\end{equation}

Then by uniqueness of the $\Z[i]$-decomposition into irreducible factors
of $k$, we deduce from $k=z'\overline z'$ that 
\begin{itemize}
\item there are $z_1\in\Z[i]$ with only irreducible factor $(1+i)$ and 
$z\in\Z[i]$ whose irreducible factors are some $z_p$ for
$p\in\mathcal P_3(k)$, such that $z'=z_1cz$,
\item we must have $z_1\overline z_1=2^{\alpha}$ and $z\overline z=m$.
\end{itemize}
Since $(1+i)$ is the only irreducible factor of $z_1$, there is a unit
element $\xi\in\Z[i]$ and a positive integer $u$ such that $z_1=\xi
(1+i)^u$. However, $z_1\overline z_1=(-i)^u(1+i)^{2u}$, and the uniqueness of
the decomposition of $2^{\alpha}$ into irreducibles implies that
$u=\alpha$. Hence $z_1=\xi(1+i)^{\alpha}$ and $\Phi(\xi z)=z'$, with $\xi
z\in\mathcal U(m)$, proving
that $\Phi$ is surjective. The injectivity of $\Phi$ is immediate.
\end{proof}

Lemma~\Cref{lembij} reduces the study of \Cref{eq_1} to the case of odd
integers $k$. On the other hand, considering the quadratic residues modulo
$4$, if $\mathcal U(k)\neq \emptyset$, then $k= 1\mod 4$. Note that this last
condition is equivalent to 
$$k= 1\mod 8\quad\text{or}\quad k=
5\mod 8.$$ 
In the following, we will put
natural families of solutions of $\mathcal U(k)$ ($k$ odd)
in bijective correspondence with (extended)
Grassmannian elements with respect to the type $C_2^{(1)}$ in \Cref{sec_Pell_C2}.
More precisely,
\begin{itemize}
\item using the affine
Grassmannian elements in type $C_2^{(1)}$ of atomic length equal to $N$
(that are naturally in bijection with the elements of $\SC_4(N)$), we
obtain a description of the solutions of the Diophantine equation
$x^2+y^2=8N+5$,
\item similarly, but using the atomic length associated with the
fundamental weight $\Lambda_1$ in type $C_2^{(1)}$, we obtain a
description of the solutions of the Diophantine equation $x^2+y^2=8N+1$ in
terms of the extended Grassmannian elements.
\end{itemize}

\subsection{Group actions on the solution set of $x^2+y^2=k$}

Let $D_8 = \langle r, s \mid r^4=s^2=(rs)^2 = 1\rangle$ be the dihedral group of order $8$.
Clearly, the formulas 
\begin{equation}
\label{actionD8}
\left\{
\begin{array}{l}
r(x,y) = (-y,x) \\
s(x,y) = (y,x)
\end{array}
\right.
 \text{ \quad \quad for all } (x,y)\in\cU(k)
\end{equation}
yield an action of $D_8$ on $\cU(k)$ when $\cU(k)$ is non-empty.
The $D_8$-orbit of $(x,y)$ in $\cU(k)$ is
\begin{equation}
\label{orbitD8}
D_8(x,y) = \left\{ (x,y), (x, -y), (-x,y), (-x,-y), (y,x), (y, -x) , (-y,x),  (-y, -x)\right\}.
\end{equation}

\begin{Rem}
\label{idaction}
Through the identification $(x,y)\mapsto x+iy\in\Z[i]$ used in \Cref{sec_x2y2k}, 
observe that the elements $r,\, s\in
 D_8$ given in \Cref{actionD8} act on $\mathcal U(k)$ respectively
 by multiplication by $i$ and $i \emph{conj}$, where $\emph{conj}$ is the complex conjugation.  
\end{Rem}

Now, all the subgroups of $ D_8$ also act on $\cU(k)$.
In \Cref{sec_param}, we will need the following subgroup actions.
\begin{itemize}
\item 
The action of $\left\langle r \right\rangle \simeq C_4$ (the cyclic group of order $4$)
\begin{equation}
\label{actionC4}
\left.
\begin{array}{l}
r(x,y) = (-y,x)
\end{array}
\right.
 \text{ \quad \quad for all } (x,y)\in\cU(k).
\end{equation}
The orbit of $(x,y)$ under this $C_4$-action is
\begin{equation}
\label{orbitC4}
C_4(x,y) = \left\{ (x,y), (-y,x), (-x,-y), (y,-x)\right\}.
\end{equation}
\item 
The action of $\left\langle rs , sr \right\rangle \simeq V_4$ (the Klein four-group)
\begin{equation}
\label{actionV4}
\left\{
\begin{array}{l}
rs(x,y) = (-x,y) \\
sr(x,y) = (x,-y)
\end{array}
\right.
 \text{ \quad \quad for all } (x,y)\in\cU(k).
\end{equation}
The orbit of $(x,y)$ under this $V_4$-action is
\begin{equation}
\label{orbitV4}
V_4(x,y) = \left\{ (x,y), (x, -y), (-x,y), (-x,-y)\right\}.
\end{equation}
\end{itemize}

\begin{Lem}
\label{lem_free_1}
Let $(x,y)\in\cU(k)\setminus\{(0,0)\}$. 
\begin{enumerate}
\item The $C_4$-orbit of $(x,y)$ described in \Cref{orbitC4} has size $4$.
\item The $D_8$-orbit of $(x,y)$ described in \Cref{orbitD8} has size $8$ if and only if 
$x,y\neq 0$ and $x\neq \pm y$.
\end{enumerate}
\end{Lem}

\begin{proof}
Point (1) is clear by the explicit description of the orbit \Cref{orbitC4}.
Now, if $x,y\neq 0$ and $x\neq \pm y$, all the elements listed in \Cref{orbitD8}
are pairwise distinct. Conversely, if $x=0$ or $y=0$ or $x=\pm y$, 
then the action of $\langle rs \rangle$, respectively  $\langle sr\rangle$,
respectively $\langle s\rangle$, is not free.
\end{proof}

\begin{Prop}
\label{charac_free}
Let $k\geq 2$. Assume that $\mathcal U(k)$ is non-empty. Then the $D_8$-action on $\mathcal U(k)$ is free if and only if $k$ is neither a perfect square nor twice a perfect square.
\end{Prop}

\begin{proof}
First, suppose that there is an integer $\ell$ such that $k=\ell^2$. Then $(\ell,0)\in \cU(k)$, and the $D_8$-orbit of $(\ell,0)$ has size $4$ by 
\Cref{lem_free_1}. Similarly, if $k=2\ell^2$, then $(\ell,\ell)\in\cU(k)$ and again, we conclude with \Cref{lem_free_1}.
It follows that $D_8$ does not act freely on $\mathcal U(k)$ in these cases.
\medskip 

Now, assume that $k$ is not a square or twice a square. Then there is an odd $m$ which is not a square, and a non negative integer $\alpha$, such that $k=2^{\alpha}m$.
As above, we identify $\Z^2$ with $\Z[i]$. As noted in \Cref{idaction}, 
the group $D_8$ is generated by the multiplication by $i$ and complex 
conjugation.  
Let $(x,y)\in \mathcal U(k)$ and set $z'=x+iy$. By \Cref{lembij} there are an integer $c$ with prime factors in $\mathcal P_3$ and $z\in \mathcal U_m$ such that $z'=(1+i)^{\alpha}c z$. 
Suppose that $\overline z' \in \langle i\rangle \cdot z'$, that is, there is
$\beta\in\N$ such that $\overline z'=i^{\beta}z'$. Since
$$\overline z'=(1-i)^{\alpha}c\overline
z=i^{3\alpha}(1+i)^{\alpha}c\overline{z},$$
we deduce that $$\overline{z}=i^{\beta-3\alpha}z.$$
Write
$z=u+iv$ with $u,\,v\in\Z$. We remark that $u\neq 0$ and $v\neq 0$ (otherwise $m=u^2$ or $m=v^2$) and if $u=\pm v$, that is $z=u(1\pm i)$, then
$$m=|i^{\beta-3\alpha}z|^2=2u^2$$ is even, which is a contradiction. Thus, by \Cref{lem_free_1}, the $D_8$-orbit of $(x,y)$ has size $8$, as required.
\end{proof}

\begin{Rem}
In fact, in the non-free case, let us show that the $D_8$-orbits of $\mathcal{U}(k)$ all have size $8$ except exactly one with size $4$. 
Assume $\mathcal U(k)\neq\emptyset$ where $k$ is a square or twice a square. 
Let $(x,y)\in \cU(k)$ and assume that the $D_8$-orbit of $(x,y)$ has size $4$.
\begin{itemize}
\item Suppose that $k=\ell^2$ is a square. If $x=\pm y$, then $x^2+y^2=2x^2$, and $\ell^2=2x^2$ which is impossible. Hence, by \Cref{lem_free_1}, $x$ or $y$ is equal to $0$. Say $y=0$ (the case $x=0$ is analogue). Then $x^2=k=\ell^2$ and $x=\pm \ell$. It follows that $(x,y)$ and $(\ell,0)$ have the same $D_8$-orbit.
\item Suppose now that $k=2\ell^2$ for some integer $\ell$. If $y=0$ (the case $x=0$ is similar) then $x^2=k=2\ell^2$ which is impossible. Hence, by \Cref{lem_free_1}, $x=\mp y$ and $2x^2=2\ell^2$ implying $x=\pm \ell$. The $D_8$-orbits of $(x,y)$ and $(\ell,\ell)$ are the same.
\end{itemize} 
\end{Rem}

We will investigate the set $\cU(k)$ in the case where $k$ is of the form
$$k=aN+b$$
for some $a,b\in\N$, where $N\in\N$ is fixed.

\begin{Cor}
\label{cor_free} Let $k$ be as above.
Assume that $b$ is neither a quadratic residue nor a twice quadratic residue modulo $a$.
Then the action \Cref{actionD8} is free.
\end{Cor}

\begin{proof}
By contraposition, suppose that the $D_8$-action is not free. 
Then by \Cref{charac_free}, $k$ is either a square or twice a square.
Since $k=aN+b$, reducing modulo $a$ implies that $b$ is either a quadratic residue or twice a quadratic residue.
\end{proof}

\begin{Rem}
In higher rank, that is for the Diophantine equation
$x_1^2+\cdots+x_n^2 =k$, there is an analogous action of the hyperoctahedral group $H_n$.
We will see in \Cref{sec_Pell_hyperoct} that these equations
naturally arise in affine types that have an underlying finite Weyl group of type $B_n$ or $C_n$ 
but that there is no hope of having a free action in general.
\end{Rem}

\subsection{The Diophantine equation $x^2 +3y^2 = k$}
\label{pell_action_eq2}

Let $k\in\N$ and consider the equation
\begin{equation}
\label{eq_2}
x^2+3y^2 =k
\end{equation}
in the variables $x,y$.
We denote again by $\cU(k)$ the set of integer solutions of \Cref{eq_2}.

\medskip 

In this situation, we do not have an action of $D_8$, but we can will consider the following actions on $\cU(k)$.
\begin{itemize}
\item The action of $V_4$ defined as in \Cref{actionV4}.
\item If $k$ is even, the action of the cyclic group $C_6=\langle t \rangle$ defined by 
\begin{equation}
    \label{action_C6}
t(x,y) = (x',y') \text{ where } x'+jy'= -j(\frac{1}{2}x+\frac{i\sqrt{3}}{2}y),
\end{equation} 
that is, $(x,y)$ is identified with the Eisenstein integer $\frac{x}{2}+i\sqrt{3}\frac{y}{2}\in\Z[j]$.
Indeed, by the same argument as \cite{BrunatNath2022}, this identification makes sense 
and this group action is well-defined if and only if $k$ is even.
\end{itemize}

\medskip

We state the following lemma, analogous to \Cref{lem_free_1}, and whose proof
is exactly the same.

\begin{Lem}
\label{lem_V4}
Let $(x,y)\in\cU(k)$.
The $V_4$-orbit of $(x,y)$ described in \Cref{orbitV4} has size $4$ if and only if $x,y\neq 0$.
In particular, if this condition holds for all $(x,y)\in\cU(k)$, $V_4$ acts freely on $\cU(k)$.
\end{Lem}

\section{Solution sets via generalised cores : rank two}
\label{sec_sols_rank2}

In this section, we study the solutions of the Pell-type equations naturally arising
from the atomic length in different Dynkin types of rank two.
This new approach is based on the study of the generalised cores (recall that these are the 
(extended) affine Grassmannian elements).
We will start by recalling the results of \cite{BrunatNath2022} (involving usual $3$-cores), reformulating them in our context, and we will use them to give a neater description of the set
\begin{equation}\label{U(12N+4)}
    \mathcal{U}(12N+4) = \{(x,y) \in \mathbb{Z}^2~|~x^2 + 3y^2 = 12N+4\}.
\end{equation}

Then, we will see that our approach gives a simple combinatorial proof of 
\cite[Theorem 7]{Alpoge2014} (involving self-conjugate $4$-cores).
The results of this section provide far-reaching extensions of these two papers, see also \cite{OnoSze1997}
for connected results on $3$ and $4$-cores.

\medskip

Surprisingly, certain Dynkin types (notably $A_n^{(1)}$, $C_n^{(1)}$, $D_{n+1}^{(2)}$ and $D_4^{(3)}$)
will provide very fruitful results, while others will not.
The wide variety of results in this section illustrate the relevance of working with both untwisted and twisted Dynkin types.

\medskip

Moreover, in type $C_2^{(1)}$, we will unsparingly use the theory developed in \Cref{sec_AL} and \Cref{sec_ext_AL}
as the $\La_1$-atomic length and the extended atomic length will play a crucial role, see \Cref{sec_Pell_C2}.
The extended atomic length will also be heavily reinvested in \Cref{same orbit A2} and \Cref{sec_variety} 
to refine the parametrisation of the solution sets investigated.

\label{sec_param}

\subsection{Type $A_2^{(1)}$}
\label{sec_Pell_A2}

In this section, 
we use the theory of extended affine Grassmannian elements to obtain a refinement of the results \cite{BrunatNath2022}.
We let $W$ be the Weyl group of type $A_2^{(1)}$.

\subsubsection{Extended affine Grassmannian elements}

Recall from \Cref{sec_ext_AG} that extended affine Grassmannian elements are in bijection with
\begin{equation}\label{dec L A2}
L = \Z\varpi_1 \oplus \Z\varpi_2 = {{M}} \sqcup (\varpi_1+ w_{0,1}w_0({{M}}) ) \sqcup (\varpi_2+ w_{0,2}w_0({{M}}) ),
\end{equation}
where the second equality follows from \Cref{L in terms of M}.
We also recall that
$\widehat W=T(L)\rtimes W_0=\Sigma\ltimes W$,
therefore any element $\widehat{w} \in \widehat{W}$ decomposes both as 
$\widehat{w} = t_r w$ for some $r \in L$ and $w \in W_0$; and as $\widehat{w} = \sigma_i  v$ for some $i = 0,1,2$ and $v \in W$. 

\medskip

It will be convenient to use matrix representation for elements of $\widehat W$,
and from now on we will always work with the basis $\cC$ of ${{M}}$ introduced in \Cref{basis_eps}.
From \Cref{fund_weights_A}, we see that the fundamental weights write
\begin{equation}
\label{eq:poidsA2}
\omega_1=\varpi_1=
\frac 1 3
\begin{pmatrix}
2 \\ -1
\end{pmatrix}
\quad \text{and}\quad
\omega_2=\varpi_2=
\frac 1 3
\begin{pmatrix}
1 \\ 1
\end{pmatrix}.
\end{equation}
From \Cref{w01w0_can_A},  we deduce that the matrix of $w_{0,1}w_0$ (viewed as a linear map on $V_0$) is
\begin{equation}\label{matrix w01w0 def}
B= \text{Mat}_{\cC}(w_{0,1}w_{0}) =
 \begin{pmatrix}
-1 & -1  \\
1 & 0    
\end{pmatrix}.
\end{equation}
Moreover, as we already explained in \Cref{section Sigma type A}, 
$\mathrm{Mat}_{\cC} (w_{0,2}w_0) = B^2$.
Writing $q=(q_1,q_2)$ one computes
\begin{equation}
\label{eq:sigmatq}
\varpi_1 + w_{0,1}w_0(q) =\frac 1 3
\begin{pmatrix}
 -3q_1-3q_2+2 \\3q_1-1
\end{pmatrix}
\quad\text{and}\quad 
\varpi_2 + w_{0,2}w_0(q) =\frac 1 3
\begin{pmatrix}
 3q_2+1 \\-3q_1-3q_2+1
\end{pmatrix},
\end{equation}
which constitute, together with $(q_1,q_2)$, the typical elements of $L$ by \Cref{dec L A2}.

\subsubsection{From $3$-cores to integer solutions and back}\label{crank and 3-cores}

Let $N\in\N$. 
Recall that we have establised in \Cref{exa_Pell_A2} (1) the following identity, for all $q = q_1\eps_1 + q_2\eps_2 \in {{M}}$.
\begin{Lem}\label{lem_A2}
For all $q\in {{M}}$, we have
$$ 12 \sL_{\La_0}(t_q)  + 4 = (6q_1 + 3q_2-2)^2 + 3(3q_2)^2.$$
\end{Lem}
Therefore, we consider the Diophantine equation 
\begin{equation}\label{eq_A2}
x^2+3y^2 = 12N+4
\end{equation}
where $N\in\N$ is fixed, and denote by $\cU(12N+4)$, as in \Cref{section
gaussian equation gen}, the set of integer solutions to \Cref{eq_A2}.
We have also explained in \Cref{exa_Pell_A2} (1) why the map
\begin{equation}\label{map p A2}
\begin{array}{ccccc}
\varphi & : & \widehat \cB(N) & \longrightarrow & \cU(12N+4) \\
   & & q &\longmapsto &(6q_1 + 3q_2-2, 3q_2).
   \end{array}
\end{equation}
is well-defined,
that is, we can attach an integral solution of \Cref{eq_A2} to any extended affine Grassmannian element.

\medskip

Let us express \Cref{map p A2} in matrix form, with respect to
the basis $\cC$ of $V_0$ and the standard basis of $\mathbb R^2$, which gives
\begin{equation}
\label{eq:matricephi_rk2}
\varphi(q)=Q_0q+u_0,
\end{equation}
where
\begin{equation}\label{matrix B def rk2}
Q_0 = 
 \begin{pmatrix}
6 & 3  \\
0 & 3 \\ 
\end{pmatrix}
\quad \text{and}\quad 
u_0 = 
 \begin{pmatrix}
 -2\\0
\end{pmatrix}.
\end{equation}

Now, in \cite[Theorem 2.1, Theorem 3.2]{BrunatNath2022}, the authors proved the following theorem.
Let
\begin{equation}\label{matrix R BN}
    R = \frac{1}{2}
 \begin{pmatrix}
1 & -3 \\
1 & 1 \\
\end{pmatrix}.
\end{equation}
Then the group $G = \langle R \rangle \subset \text{GL}(\mathbb{R}^2)$ 
is isomorphic to the cyclic group $C_6$ introduced in \Cref{pell_action_eq2}, via $R\mapsto t$.
Therefore, $G$ acts on $\mathcal{U}(12N+4)$.

\begin{Th}
\label{thm_A2}
Let $N\in\N$.
\begin{enumerate}
\item The action of $G$ on $\cU(12N+4)$ is free.
\item The set $\Im(\varphi)$ is a complete set of representatives of the $G$-orbits of $\cU(12N+4)$.
\end{enumerate}
\end{Th}

Therefore, any orbit under this action provides 6 points where only one of them comes from a  
$3$-core.
Moreover, if we denote by $\mathcal{O}_G(\varphi(q))$ the orbit of $q\in\cB(N)$ under the action $G$,
we have the following decomposition 
\begin{align}
 \mathcal{U}(12N+4) 
& = \bigsqcup\limits_{q \in \mathcal{B}(N)} \mathcal{O}_G(\varphi(q)). \label{decomposition orbite A}
\end{align}

Recall that, by \Cref{sec_cores}, that we have a bijection
\begin{equation}\label{bij_3cores}
 \cB(N) \overset{\sim}{\longleftrightarrow} \cC_3(N)
\end{equation}
where $\cC_3(N)$ is the set of $3$-core partitions of size $N$.
In particular, the following  identity is straightforward.

\begin{Cor}\label{cor_A2}
For all $N\in\N$, we have 
$$ \left| \cC_3(N) \right| = \frac{1}{6} \left| \cU(12N+4) \right|.$$
\end{Cor}

\subsubsection{The decomposition theorem for $A_2^{(1)}$}\label{sec_pig A2}

In this subsection, we will refine the decomposition \Cref{decomposition orbite A} using extended affine Grassmannian elements.
 The following theorem shows that the subgroup $\langle R^2 \rangle$ (of order $3$)
is isomorphic to the fundamental group $\Sigma$ as permutation group.

\begin{Th}\label{same orbit A2}
Let $q\in L$. We have
$\varphi(\sigma_1 t_q)=R^2\varphi(q)$.
\end{Th}

\begin{proof}
We perform a simple check using the matrix formulas previously established. On the one hand,
$$
\varphi(\sigma_1 t_q)
\underset{\substack{(\ref{eq:sigmatq})\\(\ref{map p A2})}}{=}
\begin{pmatrix} 6&3 \\ 0&3 \end{pmatrix}
\frac{1}{3}\begin{pmatrix} -3q_1-3q_2+2 \\ 3q_1-1 \end{pmatrix}
+\begin{pmatrix} -2 \\ 0 \end{pmatrix}
=
\begin{pmatrix} -3q_1-6q_2+1 \\ 3q_1-1\end{pmatrix}.
$$
On the other hand, we have
$$R^2 \varphi(q)
=\frac 1 2\begin{pmatrix}-1& -3\\1& -1\end{pmatrix}
\begin{pmatrix} 3q_1+6q_2-1 \\ 3q_1-1 \end{pmatrix}
= 
\begin{pmatrix} -3q_1-6q_2+1 \\ 3q_1-1\end{pmatrix}.
$$.
\end{proof}

\begin{Cor}\label{prop_phi_sols}
Let $q\in {{M}}$, $j \in \{1,2\}$ and $N = \sL_{\Lambda_0}(t_q)$. We have
$\varphi(\sigma_jt_q) \in\cU(12N +4)$.
\end{Cor}

\begin{proof} 
By the results of \Cref{crank and 3-cores}, if $q\in {{M}}$ then $\varphi(q)\in\cU(12N+4)$,
and we know that $G$ acts on $\cU(12N+4)$.
Therefore, $R^2\varphi(q), R^4\varphi(q)\in\cU(12N+4)$.
We conclude by \Cref{same orbit A2}.
\end{proof}

\begin{Rem}
Note that \Cref{prop_phi_sols} could also be derived from the results of \Cref{section gaussian equation gen}.
Indeed, it is straightforward from the formulas \Cref{eq:sigmatq} that $\varphi(\sigma_jt_q) \in\Z^2$ for all $j \in \{1,2\}$.
The conclusion follows  from \Cref{new_sols_ext}.
In fact, this argument will be used in \Cref{sec_A3_phi} to establish a similar result in higher rank.
\end{Rem}

\medskip

\begin{Exa}\label{example A2 p}\
We represent in \Cref{ellipses} the sets $\cU(12N+4)$ in the cases $N=0,1,6$
and give the corresponding extended affine Grassmannian elements of atomic length $N$ (or equivalently, $3$-cores of size $N$).
Let us detail the computations.
\begin{enumerate}
\item Let $N=0$.
 The empty partition is the only $3$-core of size $0$, and corresponds to the affine Grassmannian element
 $w = 1$.  
 The corresponding element in $M$ is $q = 0 = (0,0)$ (written in the basis $\cC$). 
 The associated extended affine Grassmannian elements are
 $\sigma_1 = t_{\varpi_1}s_1s_2$ and $\sigma_2 = t_{\varpi_2}s_2s_1$. 
Using \Cref{same orbit A2}, we get three integer solutions 
$$\varphi((0,0)) =(-2,0) =: X \quad, \quad 
\varphi(\varpi_1) = (1,-1) =: X^{[1]} \mand
\varphi(\varpi_2) = (1,1) =: X^{[2]}.$$
Now, acting by $\langle R^3\rangle$ on $X$, $X^{[1]}$ and $X^{[2]}$ we obtain the full solution set
$$\mathcal U(4) =\{(2,0), (1,-1), (1,1), (2,0),(-1,1),  (-1,-1), \}.$$ 
 \item Let $N=1$. 
 The partition  $\Yboxdim{5pt}\yng(1)$ is the only $3$-core of size $1$, and corresponds to the affine Grassmannian element
 $w = s_0 = t_{\alpha_1+\alpha_2} s_1s_2s_1$.  
 The corresponding element in ${{M}}$ is $q = \alpha_1 + \alpha_2 = (1,0)$ (written in the basis $\cC$). 
 The associated extended affine Grassmannian elements are
 $\sigma w = t_{\varpi_1 + w_{0,1}w_0(q)}s_1s_2 s_1s_2s_1 = t_{\varpi_1 +  w_{0,1}w_0(q)}s_2$ and $\sigma_2 w = t_{\varpi_2 +  w_{0,2}w_0(q)}s_2s_1 s_1s_2s_1 = t_{\varpi_2 +  w_{0,2}w_0(q)}s_1$.  
Using \Cref{same orbit A2}, we get three integer solutions 
$$\varphi(q) =(4,0) =: Y \quad, \quad 
\varphi(\varpi_1 + w_{0,1}w_0(q)) = (-2,2) =: Y^{[1]} \mand
\varphi(\varpi_2 + w_{0,2}w_0(q)) = (-2,-2) =: Y^{[2]}.$$
Now, acting by $\langle R^3\rangle$ on $Y$, $Y^{[1]}$ and $Y^{[2]}$ we obtain the full solution set
$$\mathcal{U}(16)=  \{(4,0),(-2,2),(-2,-2),(-4,0),(2,-2),(2,2)\}.$$
\item Let $N=6$. 
There are two  $3$-cores of size $6$, namely $\Yboxdim{4pt}\yng(4,2)$ and its conjugate $\Yboxdim{4pt}\yng(2,2,1,1)$.
These correspond to the affine Grassmannian elements
$w = t_{2\alpha_1 + \alpha_2}s_1s_2$ and $v= t_{\alpha_1 + 2\alpha_2}s_2s_1$  respectively.
The corresponding elements of ${{M}}$ are $q = 2\alpha_1 + \alpha_2 = (2,-1)$ and $q' = \alpha_1 + 2\alpha_2 = (1,1)$ respectively (written in the basis $\cC$). 
The associated extended affine Grassmannian elements are
 $\sigma w = t_{\varpi_1 + w_{0,1}w_0(q)}s_1s_2 s_1s_2s_1 = t_{\varpi_1 +  w_{0,1}w_0(q)}s_2$ and $\sigma_2 w = t_{\varpi_2 +  w_{0,2}w_0(q)}s_2s_1 s_1s_2s_1 = t_{\varpi_2 +  w_{0,2}w_0(q)}s_1$,
 and similarly for $q'$.
Using \Cref{same orbit A2}, we get three integer solutions from $q$ and three from $q'$:
$$\varphi(q) =(7,-3) =: Z \quad, \quad 
\varphi(\varpi_1 + w_{0,1}w_0(q)) = (1,5) =: Z^{[1]} \mand
\varphi(\varpi_2 + w_{0,2}w_0(q)) = (-8,-2) =: Z^{[2]}.$$
$$\varphi(q') =(7,3) =: Z' \quad, \quad 
\varphi(\varpi_1 + w_{0,1}w_0(q')) = (-8,2) =: Z'^{[1]} \mand
\varphi(\varpi_2 + w_{0,2}w_0(q')) = (1,-5) =: Z'^{[2]}.$$
Now, acting by $\langle R^3\rangle$ on $Z$, $Z^{[1]}$ and $Z^{[2]}$, and on $Z'$, $Z'^{[1]}$ and $Z'^{[2]}$,  we obtain the full solution set
as disjoint union of two orbits:
$$
\begin{array}{l}
\mathcal{U}(76)=  \{(7,-3),(1,5),(-8,-2), (-7,3), (-1,-5), (8,2)\} 
\\
\hspace{3cm}
\sqcup \, \{ (7,3), (-8,2), (1,-5), (-7,-3), (8,-2), (-1,5) \}.
\end{array}
$$
\begin{figure}[h!]
\centering 
\begin{tikzpicture}[scale=0.6]
    \draw[->] (-10,0) -- (10,0);
    \foreach \x in {-10,-9,...,10}
        \draw (\x,0.2) -- (\x,-0.2);
    \draw[->] (0,-7) -- (0,7);
    \foreach \y in {-7,-6,...,7}
        \draw (0.2,\y) -- (-0.2,\y);
    
    \draw[black] (0,0) ellipse ({sqrt(4)} and {sqrt(4/3)});
    \draw[black] (0,0) ellipse ({sqrt(16)} and {sqrt(16/3)});
    \draw[black] (0,0) ellipse ({sqrt(76)} and {sqrt(76/3)});
    
    \fill[black] (2,0) circle (4pt) node[above right, scale=1] {};
    \fill[blue] (1,1) circle (4pt) node[above right, scale=1] {$X^{[2]}$};
    \fill[green!60!black] (-1,1) circle (4pt) node[above left, scale=1] {};
    \fill[black] (-2,0) circle (4pt) node[above left, scale=1] {$X$};
    \fill[blue] (-1,-1) circle (4pt) node[below left, scale=1] {};
    \fill[green!60!black] (1,-1) circle (4pt) node[below right, scale=1] {$X^{[1]}$};
    
    \fill[black] (4,0) circle (4pt) node[above right, scale=1] {$Y$};
    \fill[blue] (2,2) circle (4pt) node[above right, scale=1] {};
    \fill[green!60!black] (-2,2) circle (4pt) node[above left, scale=1] {$Y^{[1]}$};
    \fill[black] (-4,0) circle (4pt) node[above left, scale=1] {};
    \fill[blue] (-2,-2) circle (4pt) node[below left, scale=1] {$Y^{[2]}$};
    \fill[green!60!black] (2,-2) circle (4pt) node[below right, scale=1] {};
    
    \fill[blue] (8,2) circle (4pt) ;
    \fill[black] (7,3) circle (4pt) node[anchor=south west]{${Z'}$};
    \fill[green!60!black] (1,5)  circle (4pt) node[anchor=south west]{$Z^{[1]}$};
    \fill[blue] (-1,5) circle (4pt)  {};
    \fill[black] (-7,3)  circle (4pt) {};
    \fill[green!60!black] (-8,2) circle (4pt)  node[anchor=south east]{$Z'^{[1]}$};
    \fill[blue] (-8,-2)  circle (4pt)  node[anchor=north east]{$Z^{[2]}$};
    \fill[black] (-7,-3) circle (4pt)  {};
    \fill[green!60!black] (-1,-5)  circle (4pt)  {};
    \fill[blue] (1,-5) circle (4pt)  node[anchor=north west]{$Z'^{[2]}$};
    \fill[black] (7,-3)  circle (4pt)  node[anchor=north west]{$Z$};
    \fill[green!60!black] (8,-2) circle (4pt) {};
\end{tikzpicture}
\caption{
The $G$-orbit of the $3$-cores for $N=0$ (on the small ellipse), $N=1$ (medium) and $N=6$ (large).
The different colors represent the three layers of $\widehat{W}$ from which the integer points arise (black from $W^0$, green from $\sigma_1W^0$ and blue from $\sigma_2W^0$).
}
\label{ellipses}
\end{figure}
\end{enumerate}
\end{Exa}

\subsection{Type $C_2^{(1)}$} \label{sec_Pell_C2}

We now turn to the simplest non simply-laced Dynkin type.
We obtain here two astonishing results. 
\Cref{thm_C2} simply uses the $\La_0$-atomic length, 
and gives a bijective proof of Alpoge's number-theoretic result \cite[Theorem 7]{Alpoge2014}
on the solutions of $x^2+y^2=8N+5$, which we recover in \Cref{cor_C2}.
\Cref{thm_C2_1} is both truly important, as it treats the complementary equation $x^2+y^2=8N+1$,
and very interesting, as it involves the $\La_1$-atomic length, in its extended version.
As such, it relies on virtually all the preceding results of this paper.

\subsubsection{Pell-type equation from the  basic atomic length}

Let $W$ be the Weyl group of type $C_2^{(1)}$.
By \Cref{table}, elements in ${{M}}$
write $q=q_1\sqrt{2}e_1 + q_2\sqrt{2}e_2$ with $(q_1,q_2)\in\Z^2$.
Moreover, denote $\cB(N) = \left\{ \be\in {{M}} \mid \sL_{\La_0}(t_q) = N \right\}$.
By \Cref{sec_cores_C}, we have a bijection
\begin{equation}\label{bij_sc4cores}
 \cB(N) \overset{\sim}{\longleftrightarrow} \SC_4(N)
\end{equation}
where $\SC_4(N)$ is the set of self-conjugate $4$-core partitions of size $N$.
Using the expression of the atomic length established in \Cref{thm_LA}, we find
\begin{align}
\label{LA_C2}
\sL_{\La_0}(t_q) 
&  =  4q_1^2+4q_2^2 - 3q_1-q_2.
\end{align}

We now use a trick specific to this rank 2 type.
Let $u:\R^2\to\R^2$ be the isomorphism whose matrix in the standard basis $\{e_1, e_2\}$ is
\begin{equation}
\label{iso_u}
\frac{1}{\sqrt{2}}\begin{pmatrix} 1&1\\1&-1 \end{pmatrix}.
\end{equation}
Recall from \Cref{table} that the simple roots are $\al_1=\frac{1}{\sqrt{2}}(e_1-e_2)$
and $\al_2=\sqrt{2}e_2$.
One can check that $u$ is an involutive isometry and that $u(\al_1)=e_2$, $u(\al_2)= e_1-e_2$.
In particular, $\{ u(\al_2), u(\al_1) \}$ forms a simple root system of type $B_2$.
The image lattice $M'=u({{M}})$ is 
$$M' = \{ q_1e_1+q_2e_2 \,;\, q_1,q_2\in\Z \text{ such that } q_1+q_2 \text{ is even} \},$$
and we can consider the statistic $\sL'_{\La_0}:T(M')\to\N, t_q \mapsto \sL_{\La_0}(t_{u(q)})$.
The map $u$ restricts to a bijection 
\begin{equation}
\label{bij_sc4cores_bis}
\cB(N)  \overset{\sim}{\longrightarrow} \cB'(N)
\end{equation} where
$$\cB'(N) = u(\cB(N)) =\{ q\in M' \mid \sL'_{\La_0}(t_{q}) = N \}.$$
Moreover, by \Cref{LA_C2},
we obtain for all $q=q_1e_1+q_2e_2\in M'$
$$\sL'_{\La_0}(t_{q}) = 2q_1^2 + 2q_2^2 - 2q_1-q_2$$
We can rewrite this as
$$\sL'_{\La_0}(t_{q}) = \frac{1}{8}(4q_1-2)^2 + \frac{1}{8}(4q_2 - 1)^2 -\frac{5}{8},$$
and we obtain the following equivalent identity
\begin{Lem}\label{lem_C2}
For all $q\in M'$, we have
$$8 \sL'_{\La_0}(t_{q}) + 5 = (4q_1-2)^2 + (4q_2 - 1)^2.$$
\end{Lem}

Therefore, we consider Equation \Cref{eq_1} for $a=8$ and $b=5$, that is
\begin{equation}
\label{eq_C2}
x^2+y^2 = 8N+5
\end{equation}
and the corresponding solution set $\cU(8N+5)$.
\Cref{lem_C2} ensures that the following map is well-defined
\begin{equation}
\begin{array}{cccc}
\varphi: & \cB'(N) & \longrightarrow & \cU(8N+5)
\\
& (q_1,q_2) & \longmapsto & (4q_1-2, 4q_2-1)
\end{array}
\end{equation}

\begin{Th}
\label{thm_C2}
Let $N\in\N$.
\begin{enumerate}
\item The action of $D_8$ on $\cU(8N+5)$ given in \Cref{actionD8} is free.
\item The set $\Im(\varphi)$ is a complete set of representatives of the $D_8$-orbits of $\cU(8N+5)$.
\end{enumerate}
\end{Th}

\begin{proof} Let $N\in\N$.
\begin{enumerate} 
\item By \Cref{cor_free}, it suffices to show that $5$ is neither a quadratic residue 
nor twice a quadratic residue modulo $8$. We have $5=-3 \mod 8$,
and the quadratic residues modulo $8$ are given in the following table.
$$
\begin{array}{@{}l@{\hskip 20pt} @{}l@{\hskip 20pt} @{}l@{\hskip 20pt} @{}l@{\hskip 20pt} @{}l@{\hskip 20pt} @{}l@{}}
\hline 
x & 0&\pm 1&\pm 2&\pm 3 & 4
\\
\hline
x^2 &0 & 1& 4 & 1 & 0
\\
\hline
\end{array}
$$
The result follows.
\item We need to show that each orbit contains exactly one element of $\Im(\varphi)$.
Let us prove the existence first.
Let $(x,y)\in\cU(8N+5)$ and consider the associated orbit $\cO$.
By the above table, one of $x,y$ must be $\pm 2\mod 8$, so it must be $\pm2\mod 4$, 
and  the other must be $\pm 1$ or $\pm3$, so it must be $\pm 1\mod 4$.
Therefore, 
$\cO$ contains an element $(x,y)$ verifying
$x=-2\mod 4$ and $y=-1\mod 4$. Write $x=4r_1-2$ and $y=4r_2-1$.
\begin{itemize}
\item 
If $r_1+r_2$ is even, then $(r_1,r_2)\in\cB'(N)$ and
$$(x,y) = \varphi(r_1,r_2).$$
\item 
If $r_1+r_2$ is odd, then $(1-r_1, r_2)\in\cB'(N)$.
Moreover, the pair $(-x,y)$ also belongs to $\cO$ and we have
$$(-x,y)=(4(1-r_1)-2, 4r_2-1) = \varphi(1-r_1,r_2).$$
\end{itemize}
Finally, let us prove uniqueness.
Let $(r_1,r_2), (t_1,t_2)\in\cB'(N)$ be distinct. 
Suppose that $\varphi(r_1,r_2)$ and $\varphi(t_1,t_2)$ are in the same orbit.
Then we must have $4r_2-1 = 4t_2-1$ and $4r_1-2=-4t_1+2$, therefore
$r_2=t_2$ and $r_1=1-t_1$. But then $r_1+r_2 = 1-t_1+t_2$, so $r_1+r_2$ is odd 
since $t_1+t_2$ is even, which is a contradiction.
\end{enumerate}
\end{proof}

\Cref{thm_C2} implies that
$$ \left| \cB'(N) \right| = \frac{1}{8} \left| \cU(8N+5) \right|.$$
Combining this with Bijections \Cref{bij_sc4cores} and \Cref{bij_sc4cores_bis}, we get the following result. 
This recovers \cite[Theorem 7]{Alpoge2014}.

\begin{Cor}\label{cor_C2}
For all $N\in\N$, we have 
$$ \left| \SC_4(N) \right| = \frac{1}{8} \left| \cU(8N+5) \right|.$$
\end{Cor}

\begin{Exa}
\label{exa_C2}
Take $N=40$.
One can check that there
are $3$ self-conjugate $4$-cores of size $40$,
namely 
\Yboxdim{5pt}
\Yvcentermath{0}
$$ \yng(11,8,5,4,3,2,2,2,1,1,1) \quad, \quad
\yng(10,7,6,5,4,3,2,1,1,1) \quad, \quad
\yng(8,8,6,6,4,4,2,2) \quad,
$$
which correspond respectively to the following elements $\sqrt{2}(q_1,q_2)$ of ${{M}}$:
$$ \sqrt{2}(1,-3) \quad,\quad \sqrt{2}(-1,3) \quad,\quad \sqrt{2}(-2,-2).$$
Applying $u$, we obtain the following corresponding elements  $(q_1,q_2)$ of $M'$:
$$ (-2,4) \quad,\quad (2,-4) \quad,\quad (-4,0).$$
Finally, applying $\varphi$ yields respectively the following pairs $(x,y)$:
$$ (-10,15) \quad,\quad (6,-17) \quad,\quad (-18,-1).$$
These are the representatives in $\Im(\varphi)$ of the three $D_8$-orbits
on the set of solutions of $x^2+y^2 = 325$ (there are $24$ integer solutions in total).
\end{Exa}

\subsubsection{Type $C_2^{(1)}$: a new Diophantine equation using the $\La_1$-atomic length}

\newcommand{\sH}{\mathscr{H}}
\newcommand{\cE}{\mathcal{E}}

Let us now study $\sL_{\La_1}$, that is, the $\La_1$-atomic length.
We will consider its extended version, but restricted to $T(L)$, that is, 
elements of the form $t_q$ for $q\in L$.
One can check that
$\varpi_1 = \sqrt{2}e_1$ and $\varpi_2=\frac{1}{\sqrt{2}}(e_1+e_2)$,
which yields
$$L = \left\{  \sqrt{2}q_1e_1 + \sqrt{2}q_2e_2 ~|~ q_1,q_2\in\Z \text{ such that } q_1+q_2\in\Z \right\}.$$
Recall that $\La_1=\La_0+\om_1$ where $\om_1=\frac{1}{\sqrt{2}}e_1$.
Combining \Cref{thm_LAi} and \Cref{atomic extended} (for $\La=\La_1$ and $x=q$), we have
\begin{align}
\label{LA_C2_1}
\sL_{\La_1}(t_q) 
&  =  4q_1^2+4q_2^2 +q_1-q_2.
\end{align}
In particular, we get a statistic
$$\sL_{\La_1} : T(L) \longrightarrow \N.$$

Similarly to the previous section,
we consider $\sL'_{\La_1}:T(L')\to\N, t_q\mapsto \sL_{\La_1}(t_{u(q)})$,
where 
$$L'=u(L)=\left\{ q_1e_1 + q_2e_2 \mid q_1,q_2\in\Z\right\}.$$
Note that $u$ is compatible with the decomposition $L={{M}}\sqcup (L\setminus {{M}})$, that is, 
$$u(L) =u({{M}})\sqcup u(L\setminus {{M}}).$$
Since $M'=u({{M}})$ is again the lattice of integer pairs of even sum,
$L'\setminus M' = u(L\setminus {{M}})$ is the set of integer pairs of odd sum.

\medskip

Similarly to the previous section, and in accordance with the notation of \Cref{definition B hat}, we let
$$\widehat{\cB_1}(N)= \{ \be\in L \mid \sL_{\La_1}(t_q) = N \}
\mand
\widehat{\cB_1'}(N) = u(\widehat{\cB_1}(N)) =\{ q\in L' \mid \sL'_{\La_1}(t_{q}) = N \}.$$
The map $u$ restricts to a bijection 
\begin{equation}
\label{bij_sc4cores_alt}
\widehat{\cB_1}(N)  \overset{\sim}{\longrightarrow} \widehat{\cB_1'}(N)
\end{equation}
Moreover, by \Cref{LA_C2_1},
we obtain for all $q\in L'$
\begin{equation}
\label{LA_C2_1_prime}
\sL'_{\La_1}(t_{q}) = 2q_1^2 + 2q_2^2 + q_2.
\end{equation}
We can rewrite this as
$$\sL'_{\La_1}(t_{q}) = 2q_1^2 + 2(q_2 + \frac{1}{4})^2 -\frac{1}{8},$$
and we obtain the following equivalent identity.
\begin{Lem}\label{lem_C2_1}
For all $q\in L'$, we have
$$8 \sL_{\La_1}'(t_{q}) + 1 = (4q_1)^2 + (4q_2 + 1)^2.$$
\end{Lem}

Therefore, we consider Equation \Cref{eq_1} for $a=8$ and $b=1$, that is
\begin{equation}
\label{eq_C2_1}
x^2+y^2 = 8N+1
\end{equation}
and the corresponding solution set $\cU(8N+1)$.
\Cref{lem_C2_1} ensures that the following map is well-defined
\begin{equation}
\label{eq8N+1a}
\begin{array}{cccc}
\varphi: & \widehat{\cB_1'}(N) & \longrightarrow & \cU(8N+1)
\\
& (q_1,q_2) & \longmapsto & (4q_1, 4q_2+1)
\end{array}
\end{equation}

\begin{Th}
\label{thm_C2_1}
Let $N\in\N$.
\begin{enumerate}
\item The action of the cyclic group $C_4$ on $\cU(8N+1)$ given in \Cref{actionC4} is free.
\item The sets  $\varphi(M')$ and $ \varphi(L'\setminus M')$ are disjoint. 
Moreover, $\varphi(L')$ is a complete set of representatives of the $C_4$-orbits of $\cU(8N+1)$.
\end{enumerate}
\end{Th}

\begin{proof}\
\begin{enumerate}
\item Let $(x,y)\in\cU(8N+1)$. 
In particular, $x^2+y^2 = 1 \mod 4$.
Recall the  quadratic residues modulo $8$ given in the proof of \Cref{thm_C2}.
Up to permutation, we must have, modulo $8$, $x\in\{0,4\}$ and $y\in\{\pm 1, \pm3\}$.
In particular, $x=0\mod 4$ and $y=\pm1\mod 4$.
Using \Cref{orbitC4}, we see that the $C_4$-orbit of $(0,1)$ and $(0,-1)$ are both equal to
$$\left\{ (0,1), (-1,0), (0,-1), (1,0)\right\},$$
and have size $4$.
A fortiori, the $C_4$-orbit of $(x,y)$ also has size $4$, that is, the action is free.
\item If $(x,y)\in \varphi(M')$ then $(x,y) = (4q_1,4q_2+1)$ for some $q_1,q_2$ with $q_1+q_2$ even. Then $x+y = 1\mod 8$. On the contrary, if $(x,y)\in \varphi(L'\setminus M')$ then $(x,y) = (4q_1,4q_2+1)$ for some $q_1,q_2$ with $q_1+q_2$ odd, and in this case $x+y = 5\mod 8$.
Therefore $\varphi(M')$ and $ \varphi(L'\setminus M')$ are disjoint.

Let us prove the second statement, by showing that each orbit contains a unique
$(x,y)$ such that $\varphi(q)=(x,y)$ for some $q\in L'$.
Let $(x,y)\in\cU(8N+1)$ and consider its $C_4$-orbit $\cO$.
By the same consideration as in (1), there is a unique $(x,y)\in\cO$ such that
$x=0\mod 4$ and $y= 1\mod 4$. 
\begin{itemize}
\item If $x=0\mod 8$ and $y=1\mod 8$, then $x=8k_1$ and $y=8k_2+1$.
Then $$(x,y) = (4(2k_1), 4(2k_2)+1)= \varphi (2k_1,2k_2)\in\varphi(M').$$
\item If $x=0\mod 8$ and $y=5\mod 8$, then $x=8k_1$ and $y=8k_2+5$.
Then $$(x,y) = (4(2k_1), 4(2k_2+1)+1)= \varphi (2k_1,2k_2+1)\in\varphi(L'\setminus M').$$
\item If $x=4\mod 8$ and $y=1\mod 8$, then $x=8k_1+4$ and $y=8k_2+1$.
Then $$(x,y) = (4(2k_1+1), 4(2k_2)+1)= \varphi (2k_1+1,2k_2)\in\varphi(L'\setminus M').$$
\item If $x=4\mod 8$ and $y=5\mod 8$, then $x=8k_1+4$ and $y=8k_2+5$.
Then $$(x,y) = (4(2k_1+1), 4(2k_2+1)+1)= \varphi (2k_1+1,2k_2+1)\in\varphi(M').$$
\end{itemize}
\end{enumerate}
\end{proof}

Using \Cref{bij_sc4cores_alt}, we obtain an analogue of \Cref{cor_C2} for the Diophantine equation $x^2+y^2=8N+1$.

\begin{Cor}
\label{cor_C2_1}
For all $N\in\N$, we have
$| \widehat{\cB_1}(N) | = \frac{1}{4}|\cU(8N+1)|.$
\end{Cor}

\begin{Exa} 
We illustrate in \Cref{table_8N+1}
how the solutions of $x^2+y^2 = 8N+1$ correspond to elements of $M'$
and $L'\setminus M'$ taken in combination.
\begin{figure}
\scriptsize
\rowcolors{2}{gray!25}{white}
\renewcommand\arraystretch{0.9}
$$
\begin{array}{@{}l@{\hskip 20pt} @{}l@{\hskip 20pt}  @{}l@{\hskip 20pt}  @{}l@{\hskip 20pt} @{}l@{\hskip 20pt} @{}l@{}}
\hline
N & q\in M' \mid \sL'_{\La_1}(t_{q})=N & \varphi(q) 
& q\in L'\setminus M' \mid \sL'_{\La_1}(t_q)=N &  \varphi(q) & \text{Solutions of } x^2+y^2 =  8N+1
\\
\hline
0 
    &
    \begin{array}[t]{ll}
(0, 0)
\end{array}
    &
    \begin{array}[t]{ll}
(0, 1)
\end{array}
    &
    \begin{array}[t]{ll}

\end{array}
    &
    \begin{array}[t]{ll}

\end{array}
    &
    \begin{array}[t]{ll}
(-1, 0)\\(0, -1)\\(0, 1)\\(1, 0)
\end{array}
    \\
1 
    &
    \begin{array}[t]{ll}

\end{array}
    &
    \begin{array}[t]{ll}

\end{array}
    &
    \begin{array}[t]{ll}
(0, -1)
\end{array}
    &
    \begin{array}[t]{ll}
(0, -3)
\end{array}
    &
    \begin{array}[t]{ll}
(-3, 0)\\(0, -3)\\(0, 3)\\(3, 0)
\end{array}
    \\
2 
    &
    \begin{array}[t]{ll}

\end{array}
    &
    \begin{array}[t]{ll}

\end{array}
    &
    \begin{array}[t]{ll}
(-1, 0)\\(1, 0)
\end{array}
    &
    \begin{array}[t]{ll}
(-4, 1)\\(4, 1)
\end{array}
    &
    \begin{array}[t]{ll}
(-4, -1)\\(-4, 1)\\(-1, -4)\\(-1, 4)\\(1, -4)\\(1, 4)\\(4, -1)\\(4, 1)
\end{array}
    \\
3 
    &
    \begin{array}[t]{ll}
(-1, -1)\\(1, -1)
\end{array}
    &
    \begin{array}[t]{ll}
(-4, -3)\\(4, -3)
\end{array}
    &
    \begin{array}[t]{ll}
(0, 1)
\end{array}
    &
    \begin{array}[t]{ll}
(0, 5)
\end{array}
    &
    \begin{array}[t]{ll}
(-5, 0)\\(-4, -3)\\(-4, 3)\\(-3, -4)\\(-3, 4)\\(0, -5)\\(0, 5)\\(3, -4)\\(3, 4)\\(4, -3)\\(4, 3)\\(5, 0)
\end{array}
    \\
4 
    &
    \begin{array}[t]{ll}

\end{array}
    &
    \begin{array}[t]{ll}

\end{array}
    &
    \begin{array}[t]{ll}

\end{array}
    &
    \begin{array}[t]{ll}

\end{array}
    &
    \begin{array}[t]{ll}

\end{array}
    \\
5 
    &
    \begin{array}[t]{ll}
(-1, 1)\\(1, 1)
\end{array}
    &
    \begin{array}[t]{ll}
(-4, 5)\\(4, 5)
\end{array}
    &
    \begin{array}[t]{ll}

\end{array}
    &
    \begin{array}[t]{ll}

\end{array}
    &
    \begin{array}[t]{ll}
(-5, -4)\\(-5, 4)\\(-4, -5)\\(-4, 5)\\(4, -5)\\(4, 5)\\(5, -4)\\(5, 4)
\end{array}
    \\
6 
    &
    \begin{array}[t]{ll}
(0, -2)
\end{array}
    &
    \begin{array}[t]{ll}
(0, -7)
\end{array}
    &
    \begin{array}[t]{ll}

\end{array}
    &
    \begin{array}[t]{ll}

\end{array}
    &
    \begin{array}[t]{ll}
(-7, 0)\\(0, -7)\\(0, 7)\\(7, 0)
\end{array}
    \\
7 
    &
    \begin{array}[t]{ll}

\end{array}
    &
    \begin{array}[t]{ll}

\end{array}
    &
    \begin{array}[t]{ll}

\end{array}
    &
    \begin{array}[t]{ll}

\end{array}
    &
    \begin{array}[t]{ll}

\end{array}
    \\
8 
    &
    \begin{array}[t]{ll}
(-2, 0)\\(2, 0)
\end{array}
    &
    \begin{array}[t]{ll}
(-8, 1)\\(8, 1)
\end{array}
    &
    \begin{array}[t]{ll}
(-1, -2)\\(1, -2)
\end{array}
    &
    \begin{array}[t]{ll}
(-4, -7)\\(4, -7)
\end{array}
    &
    \begin{array}[t]{ll}
(-8, -1)\\(-8, 1)\\(-7, -4)\\(-7, 4)\\(-4, -7)\\(-4, 7)\\(-1, -8)\\(-1, 8)\\(1, -8)\\(1, 8)\\(4, -7)\\(4, 7)\\(7, -4)\\(7, 4)\\(8, -1)\\(8, 1)
\end{array}
    \\
9 
    &
    \begin{array}[t]{ll}

\end{array}
    &
    \begin{array}[t]{ll}

\end{array}
    &
    \begin{array}[t]{ll}
(-2, -1)\\(2, -1)
\end{array}
    &
    \begin{array}[t]{ll}
(-8, -3)\\(8, -3)
\end{array}
    &
    \begin{array}[t]{ll}
(-8, -3)\\(-8, 3)\\(-3, -8)\\(-3, 8)\\(3, -8)\\(3, 8)\\(8, -3)\\(8, 3)
\end{array}
    \\
10 
    &
    \begin{array}[t]{ll}
(0, 2)
\end{array}
    &
    \begin{array}[t]{ll}
(0, 9)
\end{array}
    &
    \begin{array}[t]{ll}

\end{array}
    &
    \begin{array}[t]{ll}

\end{array}
    &
    \begin{array}[t]{ll}
(-9, 0)\\(0, -9)\\(0, 9)\\(9, 0)
\end{array}
    \\
11 
    &
    \begin{array}[t]{ll}

\end{array}
    &
    \begin{array}[t]{ll}

\end{array}
    &
    \begin{array}[t]{ll}
(-2, 1)\\(2, 1)
\end{array}
    &
    \begin{array}[t]{ll}
(-8, 5)\\(8, 5)
\end{array}
    &
    \begin{array}[t]{ll}
(-8, -5)\\(-8, 5)\\(-5, -8)\\(-5, 8)\\(5, -8)\\(5, 8)\\(8, -5)\\(8, 5)
\end{array}
    \\
12 
    &
    \begin{array}[t]{ll}

\end{array}
    &
    \begin{array}[t]{ll}

\end{array}
    &
    \begin{array}[t]{ll}
(-1, 2)\\(1, 2)
\end{array}
    &
    \begin{array}[t]{ll}
(-4, 9)\\(4, 9)
\end{array}
    &
    \begin{array}[t]{ll}
(-9, -4)\\(-9, 4)\\(-4, -9)\\(-4, 9)\\(4, -9)\\(4, 9)\\(9, -4)\\(9, 4)
\end{array}
    \\
13 
    &
    \begin{array}[t]{ll}

\end{array}
    &
    \begin{array}[t]{ll}

\end{array}
    &
    \begin{array}[t]{ll}

\end{array}
    &
    \begin{array}[t]{ll}

\end{array}
    &
    \begin{array}[t]{ll}

\end{array}
    \\
14 
    &
    \begin{array}[t]{ll}
(-2, -2)\\(2, -2)
\end{array}
    &
    \begin{array}[t]{ll}
(-8, -7)\\(8, -7)
\end{array}
    &
    \begin{array}[t]{ll}

\end{array}
    &
    \begin{array}[t]{ll}

\end{array}
    &
    \begin{array}[t]{ll}
(-8, -7)\\(-8, 7)\\(-7, -8)\\(-7, 8)\\(7, -8)\\(7, 8)\\(8, -7)\\(8, 7)
\end{array}
    \\
\hline
\end{array}
$$
\caption{Solutions of $x^2+y^2 = 8N+1$ via affine Grassmannian elements of type $C_2^{(1)}$.}
\label{table_8N+1}
\end{figure}
\end{Exa}

\begin{Rem}\
\begin{enumerate}
\item Let $\sH:T({{M}})\to \N$ defined by the formula
$$\sH(t_q)=4q_1^2+4q_2^2-3q_1+3q_2$$
for all $q=q_1\sqrt{2}e_1+ q_2\sqrt{2}e_2\in {{M}}.$
Again, using the change of variables
induced by $u$, we get the statistic $\sH':T(M')\to \N, t_q\mapsto \sH(t_{u(q)})$ whose formula is
$$\sH'(t_{q})=2q_1^2+2q_2^2-3q_2.$$
This can be rewritten as
$$\sH'(t_{q}) = 2q_1^2 + 2(q_2 - \frac{3}{4})^2 -\frac{9}{8},$$
and we obtain the following equivalent identity:
$$8 (\sH'(t_{q})+1) + 1 = (4q_1)^2 + (4q_2 -3)^2$$
for all $q\in M'$.
Now, observe that
$$\sH'(t_{q})= 2q_1^2+2q_2^2-3q_2 = 2q_1^2+2(q_2-1)^2+(q_2-1) + 1 =
\sL'_{\La_1}(t_{(q_1,q_2-1)}) + 1.$$ 
Therefore, it is equivalent to consider the combination of $\sL_{\La_1}$ and $\sH$ on the lattice $T({{M}})$
or only $\sL_{\La_1}$ on the finer lattice $T(L)$.
\item 
We leave open the problem of finding a combinatorial model as in \Cref{sec_comb_models} in which the size of partitions is given by the statistics $\sL_{\La_1}$ and $\sH$ on $T({{M}})$. 
\end{enumerate}

\end{Rem}

\subsection{Type $D_3^{(2)}$}
\label{Pell_D32}

Let $W$ be the Weyl group of type $D_3^{(2)}$.
By \Cref{table}, elements in ${{M}}$
write $q=q_1\sqrt{2}e_1 + q_2\sqrt{2}e_2$ with $(q_1,q_2)\in\Z^2$
Moreover, denote $\cB(N) = \left\{ \be\in {{M}} \mid \sL_{\La_0}(t_q) = N \right\}$.
By \Cref{sec_cores_Dt}, we have a bijection
\begin{equation}\label{bij_cores_D3t}
 \cB(N) \overset{\sim}{\longleftrightarrow} \cD_6(N).
\end{equation}

Using the expression of the atomic length established in \Cref{thm_LA}, we find
\begin{align*}
\sL_{\La_0}(t_q) 
&  =  3q_1^2+3q_2^2 - 2q_1-q_2
\\
&  = 3(q_1 -\frac{2}{3})^2 + 3(q_2 - \frac{1}{3})^2 -\frac{5}{12} 
\end{align*}

We can get rid of the denominators and obtain the following equivalent identity.

\begin{Lem}\label{lem_D3t}
For all $\be\in {{M}}$, we have 
$$12 \sL_{\La_0}(t_q)  + 5 = (6q_1 -2)^2 + (6q_2 - 1)^2.$$
\end{Lem}

Therefore, we consider Equation \Cref{eq_1} for $a=12$ and $b=5$, that is
\begin{equation}
\label{eq_D3t}
x^2+y^2 = 12N+5
\end{equation}
and the corresponding solution set $\cU(12N+5)$.
\Cref{lem_D3t} ensures that the following map is well-defined
\begin{equation}
\begin{array}{cccc}
\varphi: & \cB(N) & \longrightarrow & \cU(12N+5)
\\
& (q_1,q_2) & \longmapsto & (6q_1-2, 6q_2-1)
\end{array}
\end{equation}

\begin{Th}
\label{thm_D3t}
Let $N\in\N$.
\begin{enumerate}
\item The action of $D_8$ on $\cU(12N+5)$ given in \Cref{actionD8} is free.
\item The set $\Im(\varphi)$ is a complete set of representatives of the $D_8$-orbits of $\cU(12N+5)$.
\end{enumerate}
\end{Th}

\begin{proof} Let $N\in\N$.
\begin{enumerate} 
\item By \Cref{cor_free}, it suffices to show that $5$ is neither a quadratic residue 
nor twice a quadratic residue modulo $12$. 
The quadratic residues modulos $12$ are given in the following table.
$$
\begin{array}{@{}l@{\hskip 20pt} @{}l@{\hskip 20pt} @{}l@{\hskip 20pt} @{}l@{\hskip 20pt} @{}l@{\hskip 20pt} @{}l@{\hskip 20pt} @{}l@{\hskip 20pt} @{}l@{}}
\hline 
x & 0&\pm 1&\pm 2&\pm 3&\pm 4&\pm 5& 6
\\
\hline
x^2 &0&1&4&-3&4&1&0
\\
\hline
\end{array}
$$
and the result follows.
\item We need to show that each orbit contains exactly one element of $\Im(\varphi)$.
In order to do this, we look at residues modulo $6$.
Let $(x,y)\in\cU(12N+5)$, so that
$$x^2+y^2 = -1\mod 6.$$ 
The quadratic residues modulo $6$ are
$$
\begin{array}{@{}l@{\hskip 20pt} @{}l@{\hskip 20pt} @{}l@{\hskip 20pt} @{}l@{\hskip 20pt} @{}l@{}}
\hline 
x & 0&\pm 1&\pm 2& 3
\\
\hline
x^2 &0&1&-2&3
\\
\hline
\end{array}
$$

Therefore, we must have
 $(x\mod 6, y\mod 6)\in\{ (\pm 1, \pm 2), (\pm 1, \pm2) \}$.
We now observe (see \Cref{orbitD8}) that this set is precisely the orbit of $(x\mod 6, y\mod 6)$ under the action \Cref{actionD8} of $D_8$.
Therefore, the pair $(-2,-1)$ yields a unique element in the orbit of $(x,y)$, that writes $(6q_1-2,6q_2-1)$, and this element clearly lies
in $\Im(\varphi)$.
\end{enumerate}
\end{proof}

\Cref{thm_C2} implies that
$$ \left| \cB(N) \right| = \frac{1}{8} \left| \cU(12N+5) \right|.$$
Combining this with Bijection \Cref{bij_cores_D3t}, we get the following result.

\begin{Cor}\label{cor_D3t}
For all $N\in\N$, we have 
$$ \left| \cD_6(N) \right| = \frac{1}{8} \left| \cU(12N+5) \right|.$$
\end{Cor}

\begin{Exa} 
\label{exa_D32_2}
Recall \Cref{exa_D32_1} (2), where we had determined $\cD_6(35)$.
We found the following partitions

$$
\Yboxdim{6pt}
\yng(17,11,5,2) \quad,\quad
\yng(16,10,5,4) \quad,\quad
\yng(13,10,7,4,1),
$$

whose corresponding double distinct partitions are respectively
$$
\Yboxdim{6pt}
\young(
!\gry<>!\white<><><><><><><><><><><><><><><><><>,!\gry<><>!\white<><><><><><><><><><><>,!\gry<><><>!\white<><><><><>,!\gry<><><><>!\white<><>,!\gry<><><><>,!\gry<><><>,!\gry<><><>,!\gry<><>,!\gry<><>,!\gry<><>,!\gry<><>,!\gry<><>,!\gry<>,!\gry<>,!\gry<>,!\gry<>,!\gry<>)\quad,\quad
\young(
!\gry<>!\white<><><><><><><><><><><><><><><><>,!\gry<><>!\white<><><><><><><><><><>,!\gry<><><>!\white<><><><><>,!\gry<><><><>!\white<><><><>,!\gry<><><><>,!\gry<><><><>,!\gry<><><><>,!\gry<><>,!\gry<><>,!\gry<><>,!\gry<><>,!\gry<>,!\gry<>,!\gry<>,!\gry<>,!\gry<>)\quad,\quad
\young(
!\gry<>!\white<><><><><><><><><><><><><>,!\gry<><>!\white<><><><><><><><><><>,!\gry<><><>!\white<><><><><><><>,!\gry<><><><>!\white<><><><>,!\gry<><><><><>!\white<>,!\gry<><><><>,!\gry<><><><>,!\gry<><><>,!\gry<><><>,!\gry<><>,!\gry<><>,!\gry<>,!\gry<>).
$$
The corresponding vectors, computed by the abacus method illustrated in \Cref{exa_cores_A}
are respectively
$$(0,-3,1,0,-1,3)\quad,\quad (0,-1.-3,0,3,1)\quad, \quad (0,3,-2,0,2,-3).$$
Recall that by \Cref{sec_cores_Dt}, these are of the form $(0,q_1,q_2,0,-q_2,-q_1)$
and we can extract the associated elements $(q_1,q_2)\in {{M}}$,
which constitute the set $\cB(35)$:
$$(-3,1)\quad,\quad (-1.-3)\quad, \quad (3,-2).$$
Finally, we use the map $\varphi$ to get the corresponding solutions of $x^2+y^2 = 425$, which are
$$(-20, 5)\quad, \quad(-8, -19)\quad,\quad (16, -13).$$
The remaining solutions are obtained by letting $D_8$ act on these particular solutions, that is by 
allowing sign changes and component switching.
\end{Exa}

\subsection{Type $A_4^{(2)}$}
\label{sec_Pell_A42}

Let $W$ be the Weyl group of type $A_4^{(2)}$.
By \Cref{table}, elements in ${{M}}$
write $q=q_1e_1 + q_2e_2$ with $(q_1,q_2)\in\Z^2$
Moreover, denote $\cB(N) = \left\{ \be\in {{M}} \mid \sL_{\La_0}(t_q) = N \right\}$.
For $q=q_1e_1 + q_2e_2\in {{M}}$, 
$$ 
|\be|^2 = q_1^2+q_2^2 
$$
and
$$
\h(q) = \frac{3}{2}q_1 + \frac{1}{2}q_2.
$$

By \Cref{thm_LA}, we obtain
\begin{align*}
\sL_{\La_0}(t_q) &  =  \frac{5}{2}q_1^2 + \frac{5}{2}q_2^2  -\frac{3}{2} q_1 - \frac{1}{2}q_2
 =  \frac{5}{2}(q_1 - \frac{3}{10})^2 + 
\frac{5}{2}(q_2-\frac{1}{10})^2 - \frac{1}{4}
\end{align*}
which is equivalent to the following identity.

\begin{Lem}\label{lem_A42}
For all $\be\in {{M}}$, we have 
$$ 40 \sL_{\La_0}(t_q) + 10  = (10q_1-3)^2 + (10q_1-1)^2.$$
\end{Lem}

Therefore, we now consider Equation \Cref{eq_1} for $a=40$ and $b=10$, that is
\begin{equation}
\label{eq_A4t}
x^2+y^2 = 40N+10
\end{equation}
and the corresponding solution set $\cU(40N+10)$.
\Cref{lem_A42} ensures that the following map is well-defined
\begin{equation}
\begin{array}{cccc}
\varphi: & \cB(N) & \longrightarrow & \cU(40N+10)
\\
& (q_1,q_2) & \longmapsto & (10q_1-3, 10q_2-1)
\end{array}
\end{equation}

Unfortunately, here, no direct analogue of Theorems \ref{thm_A2}, \ref{thm_C2}, \ref{thm_C2_1} and \ref{thm_D3t}, as illustrated in the following example.
However, we are able to prove the following.
Recall the $D_8$ action on $\cU(40N+10)$ given in \Cref{actionD8}.

\begin{Prop}
\label{prop_A42}
Let $q\in {{M}}$. Then the $D_8$-orbit of $\varphi(q)$ has size $8$.
\end{Prop}

\begin{proof} Let $q\in {{M}}$ and  $(x,y)=\varphi(q)$.
We have $x=-3\mod 10$ and $y=-1\mod 10$.
Therefore, \Cref{lem_free_1} applies and the result follows.
\end{proof}

\begin{Exa}
We determine in \Cref{table_40N+10}
the affine Grassmannian elements of atomic length $N$, give the corresponding solution
of $x^2+y^2 = 40N+10$, and compare this to the remaining solutions of the equation.
We see that the action of $D_8$ is not necessarily free. 
For instance, for $N=1$, the orbit of $(5,5)$ has size only $4$.
However, the orbit of $\varphi(q)$ has size $8$ for all $q\in {{M}}$, which illustrates \Cref{prop_A42}.
It is unclear how to get the remaining solutions.
\begin{figure}
\renewcommand\arraystretch{0.9}
$$
\footnotesize
\rowcolors{2}{gray!25}{white}
\begin{array}{@{}l@{\hskip 20pt} @{}l@{\hskip 20pt} @{}l@{\hskip 20pt} @{}l@{}}
\hline
N & q\in\cB(N) & \varphi(q)
& \text{solutions of } x^2+y^2 =  40N+10
\\
\hline
0 
    &
    \begin{array}[t]{ll}
(0, 0)
\end{array}
    &
    \begin{array}[t]{ll}
(-3, -1)
\end{array}
    &
    \begin{array}[t]{ll}
(-3, -1)\\(-3, 1)\\(-1, -3)\\(-1, 3)\\(1, -3)\\(1, 3)\\(3, -1)\\(3, 1)
\end{array}
    \\
1 
    &
    \begin{array}[t]{ll}
(1, 0)
\end{array}
    &
    \begin{array}[t]{ll}
(7, -1)
\end{array}
    &
    \begin{array}[t]{ll}
(-7, -1)\\(-7, 1)\\(-5, -5)\\(-5, 5)\\(-1, -7)\\(-1, 7)\\(1, -7)\\(1, 7)\\(5, -5)\\(5, 5)\\(7, -1)\\(7, 1)
\end{array}
    \\
2 
    &
    \begin{array}[t]{ll}
(0, 1)
\end{array}
    &
    \begin{array}[t]{ll}
(-3, 9)
\end{array}
    &
    \begin{array}[t]{ll}
(-9, -3)\\(-9, 3)\\(-3, -9)\\(-3, 9)\\(3, -9)\\(3, 9)\\(9, -3)\\(9, 3)
\end{array}
    \\
3 
    &
    \begin{array}[t]{ll}
(0, -1)\\(1, 1)
\end{array}
    &
    \begin{array}[t]{ll}
(-3, -11)\\(7, 9)
\end{array}
    &
    \begin{array}[t]{ll}
(-11, -3)\\(-11, 3)\\(-9, -7)\\(-9, 7)\\(-7, -9)\\(-7, 9)\\(-3, -11)\\(-3, 11)\\(3, -11)\\(3, 11)\\(7, -9)\\(7, 9)\\(9, -7)\\(9, 7)\\(11, -3)\\(11, 3)
\end{array}
    \\
4 
    &
    \begin{array}[t]{ll}
(-1, 0)\\(1, -1)
\end{array}
    &
    \begin{array}[t]{ll}
(-13, -1)\\(7, -11)
\end{array}
    &
    \begin{array}[t]{ll}
(-13, -1)\\(-13, 1)\\(-11, -7)\\(-11, 7)\\(-7, -11)\\(-7, 11)\\(-1, -13)\\(-1, 13)\\(1, -13)\\(1, 13)\\(7, -11)\\(7, 11)\\(11, -7)\\(11, 7)\\(13, -1)\\(13, 1)
\end{array}
    \\
5 
    &
    \begin{array}[t]{ll}

\end{array}
    &
    \begin{array}[t]{ll}

\end{array}
    &
    \begin{array}[t]{ll}

\end{array}
    \\
6 
    &
    \begin{array}[t]{ll}
(-1, 1)
\end{array}
    &
    \begin{array}[t]{ll}
(-13, 9)
\end{array}
    &
    \begin{array}[t]{ll}
(-15, -5)\\(-15, 5)\\(-13, -9)\\(-13, 9)\\(-9, -13)\\(-9, 13)\\(-5, -15)\\(-5, 15)\\(5, -15)\\(5, 15)\\(9, -13)\\(9, 13)\\(13, -9)\\(13, 9)\\(15, -5)\\(15, 5)
\end{array}
    \\
\hline
\end{array}
$$
\caption{Solutions of $x^2+y^2 = 40N+10$ versus affine Grassmannian elements of type $A_4^{(2)}$.}
\label{table_40N+10}
\end{figure}
\end{Exa}

\subsection{Type $G_2^{(1)}$}
\label{sec_Pell_G21}

Let $W$ be the Weyl group of type $G_2^{(1)}$.
By \Cref{table}, elements in ${{M}}$
write $q=q_1e_1 + q_2e_2+q_3e_3$ with $(q_1,q_2,q_3)\in\Z^3$
such that $q_3=-q_1-q_2$.
One computes
$$|\be|^2 = 2q_1^2+2q_2^2 + 2q_1q_2
\mand
\h(q) = -4q_1 -5q_2.$$
This yields
\begin{align*}
\sL_{\La_0}(t_q) &  =  6q_1^2 + 6q_2^2 + 6q_1q_2 + 4q_1+5q_2
\\
& =  6(q_1 + \frac{1}{2}q_2+ \frac{1}{3})^2 + 
\frac{9}{2}(q_2+\frac{1}{3})^2 - \frac{7}{6}
\end{align*}

which is equivalent to the following identity.

\begin{Lem}
\label{lem_G2} For all $q\in {{M}}$, we have
$$ 6\sL_{\La_0}(t_q) +7 = (6q_1+3q_2+2)^2 + 3(3q_2+1)^2.$$
\end{Lem}

Therefore, we now consider Equation \Cref{eq_2} for $a=6$ and $b=7$, that is

\begin{equation}
\label{eq_G2}
x^2+3y^2 = 6N+7
\end{equation}
and the corresponding solution set $\cU(6N+7)$.
\Cref{lem_G2} ensures that the following map is well-defined
\begin{equation}
\begin{array}{cccc}
\varphi: & \cB(N) & \longrightarrow & \cU(6N+7)
\\
& (q_1,q_2) & \longmapsto & (6q_1+3q_2+2, 3q_2+1)
\end{array}
\end{equation}

Here again, we have no direct analogue of Theorems \ref{thm_A2}, \ref{thm_C2}, \ref{thm_C2_1} and \ref{thm_D3t}, as illustrated in the example below.
However, we are able to prove the following.
Recall the $V_4$ action on $\cU(6N+7)$ given in \Cref{actionV4}.

\begin{Prop}
\label{prop_G21}
Let $q\in {{M}}$. Then the $V_4$-orbit of $\varphi(q)$ has size $4$.
\end{Prop}

\begin{proof}
Let $q\in {{M}}$, denote $(x,y)=\varphi(q)\in\cU(6N+7)$.
Since $x=2\mod 3$ and $y=1\mod 3$, 
\Cref{lem_V4} applies, which proves the claim.
\end{proof}

\begin{Exa}
We determine in \Cref{table_6N+7}
the affine Grassmannian elements of type $G_2^{(1)}$ of atomic length $N$, give the corresponding solution
of $x^2+3y^2 = 6N+7$, and compare this to the remaining solutions of the equation.
There is an action of $V_4$ (by allowing signs in each component) which 
is not free. For instance, for $N=3$, the orbit of $(5,0)$ has size only $2$.
As in \Cref{sec_Pell_A42}, it is unclear how to get the remaining solutions.

\begin{figure}[h!]
$$
\footnotesize
\rowcolors{2}{gray!25}{white}
\renewcommand\arraystretch{0.9}
\begin{array}{@{}l@{\hskip 20pt} @{}l@{\hskip 20pt} @{}l@{\hskip 20pt} @{}l@{}}
\hline
N & q\in\cB(N) & \varphi(q)
& \text{solutions of } x^2+3y^2 =  6N+7
\\
\hline
0 
    &
    \begin{array}[t]{ll}
(0, 0)
\end{array}
    &
    \begin{array}[t]{ll}
(2, 1)
\end{array}
    &
    \begin{array}[t]{ll}
(-2, -1)\\(-2, 1)\\(2, -1)\\(2, 1)
\end{array}
    \\
1 
    &
    \begin{array}[t]{ll}
(0, -1)
\end{array}
    &
    \begin{array}[t]{ll}
(-1, -2)
\end{array}
    &
    \begin{array}[t]{ll}
(-1, -2)\\(-1, 2)\\(1, -2)\\(1, 2)
\end{array}
    \\
2 
    &
    \begin{array}[t]{ll}
(-1, 0)
\end{array}
    &
    \begin{array}[t]{ll}
(-4, 1)
\end{array}
    &
    \begin{array}[t]{ll}
(-4, -1)\\(-4, 1)\\(4, -1)\\(4, 1)
\end{array}
    \\
3 
    &
    \begin{array}[t]{ll}

\end{array}
    &
    \begin{array}[t]{ll}

\end{array}
    &
    \begin{array}[t]{ll}
(-5, 0)\\(5, 0)
\end{array}
    \\
4 
    &
    \begin{array}[t]{ll}

\end{array}
    &
    \begin{array}[t]{ll}

\end{array}
    &
    \begin{array}[t]{ll}
(-2, -3)\\(-2, 3)\\(2, -3)\\(2, 3)
\end{array}
    \\
5 
    &
    \begin{array}[t]{ll}
(1, -1)
\end{array}
    &
    \begin{array}[t]{ll}
(5, -2)
\end{array}
    &
    \begin{array}[t]{ll}
(-5, -2)\\(-5, 2)\\(5, -2)\\(5, 2)
\end{array}
    \\
\hline
\end{array}
$$
\caption{Solutions of $x^2+3y^2 = 6N+7$ 
versus affine Grassmannian elements of type $G_2^{(1)}$.}
\label{table_6N+7}
\end{figure}
\end{Exa}

\subsection{Type $D_4^{(3)}$}
\label{sec_Pell_D43}

Let $W$ be the Weyl group of type $D_4^{(3)}$.
By \Cref{table}, elements in ${{M}}$
write $q=q_1e_1 + q_2e_2+q_3e_3$ with $(q_1,q_2,q_3)\in\Z^3$
such that $q_3=-q_1-q_2$.
One computes
$$|\be|^2 = 2q_1^2+2q_2^2 + 2q_1q_2
\mand
\h(q) = -2q_1 -3q_2.$$
This yields
\begin{align*}
\sL_{\La_0}(t_q) &  =  4q_1^2 + 4q_2^2 + 4q_1q_2 + 2q_1+3q_2
\\
& =  4(q_1 + \frac{1}{2}q_2+ \frac{1}{4})^2 + 
3(q_2+\frac{1}{3})^2 - \frac{7}{12}
\end{align*}
which is equivalent to the following identity.

\begin{Lem}
\label{lem_D43} For all $q\in {{M}}$, we have
$$ 12\sL_{\La_0}(t_q) +7 = (6q_2+2)^2 + 3(4q_1+2q_2+1)^2.$$
\end{Lem}

Therefore, we consider Equation \Cref{eq_2} 
for $a=12$ and $b=7$, that is
\begin{equation}
\label{eq_D43}
x^2+3y^2 = 12N+7
\end{equation}
and the corresponding solution set $\cU(12N+7)$.
\Cref{lem_D43} ensures that the following map is well-defined
\begin{equation}
\begin{array}{cccc}
\varphi: & \cB(N) & \longrightarrow & \cU(12N+7)
\\
& (q_1,q_2) & \longmapsto & (6q_2+2, 4q_1 +2q_2+1)
\end{array}
\end{equation}

\begin{Th}
\label{thm_D43}
Let $N\in\N$.
\begin{enumerate}
\item The action of the Klein four-group $V_4$ on $\cU(12N+7)$ given in \Cref{actionV4} is free.
\item The set $\Im(\varphi)$ is a complete set of representatives of the $V_4$-orbits of $\cU(12N+7)$.
\end{enumerate}
\end{Th}

\begin{proof} Let $N\in\N$.
\begin{enumerate} 
\item 
Let $(x,y)\in \cU(12N+7)$. This implies that
$x^2=1\mod 3$, so $x=\pm 1\mod 3$.
Moreover, 
$x^2-y^2=3\mod 4$,
so by examining the quadratic residues modulo $4$,
one must have $x^2=0\mod 4$ and $y^2=1\mod 4$,
therefore $x=0\mod 2$.
By the Chinese remainder theorem (or by a direct check), these conditions on $x$
imply that $x=\pm 2\mod 6$.
Moreover, since $y^2=1\mod 4$, $y=\pm1\mod 4$.
We obtain $x\in\{\pm 1, \pm 4\} \mod 12$ and
 $y\in\{\pm 1, \pm 3, \pm 5, \pm 7\} \mod 12$.
Therefore, \Cref{lem_V4} applies and the action is free.
\item Let $\cO$ be a $V_4$-orbit in $\cU(12N+7)$.
We have seen in (1) that there exists a unique $(x,y)\in \cO$ such that $x=2\mod 6$ and $y=\pm 1\mod 4$.
Write $x=6q_2+2$ for some $q_2\in\Z$.
\begin{itemize}
\item Suppose that $q_2$ is even. 
Then we choose the element $y$ such that $y=1\mod 4$, that is, $y=4k+1$ for some $k\in \Z$.
Thus, we have
$$(x,y) = (6q_2+2, 4(k-\frac{q_2}{2}) +2q_2+1)=\varphi(q_1,q_2)\text{\quad with\quad} q_1=k-\frac{q_2}{2}.$$
\item Suppose that $q_2$ is odd. 
Then we choose the element $y$ such that $y=-1\mod 4$, that is, $y=4k-1$ for some $k\in \Z$.
Thus, we have
$$(x,y) = (6q_2+2, 4(k-\frac{q_2+1}{2}) +2q_2+1)=\varphi(q_1,q_2)\text{\quad with\quad} q_1=k-\frac{q_2+1}{2}.$$
\end{itemize}
\end{enumerate}
\end{proof}

Recall that we had a bijection 
$$\cB(N)\overset{\sim}{\longrightarrow} \cD_4^\flat(N),$$
where $\cD_4^\flat$ is the combinatorial model of \Cref{sec_cores_D43}.

\begin{Cor}
\label{cor_D43}
For all $N\in\N$, we have
$$|\cD_4^\flat(N) | = \frac{1}{4}|\cU(12N+7)|.$$
\end{Cor}

\begin{Exa}
\label{exa_D43}
For each $0\leq N\leq 7$,
we have recorded in \Cref{table_12N+7} which elements $\be\in {{M}}$ have 
atomic length $N$ and their image under $\varphi$.
In the last column, we have computed $\cU(12N+7)$.
Comparing the last two columns illustrates \Cref{thm_D43}.
For completeness, we illustrate \Cref{cor_D43} by listing the elements of $\cD_4^\flat$
in the second column, which can be computed from the third column 
using Formulas \Cref{formules_david} of \Cref{sec_comb_models} (see also \Cref{exa_cores_D43}).
\begin{figure}
\footnotesize
$$
\rowcolors{2}{gray!25}{white}
\renewcommand\arraystretch{1}
\begin{array}{@{}l@{\hskip 20pt} @{}l@{\hskip 20pt} @{}l@{\hskip 20pt} @{}l@{\hskip 20pt} @{}l@{}}
\hline
N &  \text{elements of }\cD_4^\flat &\text{elements } \be\in \cB(N) & \text{elements of } \varphi(\cB(N)) & \text{solutions of } x^2+3y^2 = 12 N +7
\\
\hline
0 
    &
\emptyset
    &
    \begin{array}[t]{ll}
(0, 0)
\end{array}
    &
    \begin{array}[t]{ll}
(2, 1)
\end{array}
    &
    \begin{array}[t]{ll}
(-2, -1)\\(-2, 1)\\(2, -1)\\(2, 1)
\end{array}
    \\
1 
&
\Yboxdim{5pt}
\yng(1)
    &
    \begin{array}[t]{ll}
(0, -1)
\end{array}
    &
    \begin{array}[t]{ll}
(-4, -1)
\end{array}
    &
    \begin{array}[t]{ll}
(-4, -1)\\(-4, 1)\\(4, -1)\\(4, 1)
\end{array}
    \\
2 
&
 \Yboxdim{5pt} \yng(2)
    &
    \begin{array}[t]{ll}
(-1, 0)
\end{array}
    &
    \begin{array}[t]{ll}
(2, -3)
\end{array}
    &
    \begin{array}[t]{ll}
(-2, -3)\\(-2, 3)\\(2, -3)\\(2, 3)
\end{array}
    \\
3 
&
 \Yboxdim{5pt} \yng(2,1)
    &
    \begin{array}[t]{ll}
(1, -1)
\end{array}
    &
    \begin{array}[t]{ll}
(-4, 3)
\end{array}
    &
    \begin{array}[t]{ll}
(-4, -3)\\(-4, 3)\\(4, -3)\\(4, 3)
\end{array}
    \\
4 
&
    &
    \begin{array}[t]{ll}

\end{array}
    &
    \begin{array}[t]{ll}

\end{array}
    &
    \begin{array}[t]{ll}

\end{array}
    \\
5 
&
 \Yboxdim{5pt} \yng(3,2)
    &
    \begin{array}[t]{ll}
(-1, 1)
\end{array}
    &
    \begin{array}[t]{ll}
(8, -1)
\end{array}
    &
    \begin{array}[t]{ll}
(-8, -1)\\(-8, 1)\\(8, -1)\\(8, 1)
\end{array}
    \\
6 
&
 \Yboxdim{5pt} \yng(4,2)
    &
    \begin{array}[t]{ll}
(1, 0)
\end{array}
    &
    \begin{array}[t]{ll}
(2, 5)
\end{array}
    &
    \begin{array}[t]{ll}
(-2, -5)\\(-2, 5)\\(2, -5)\\(2, 5)
\end{array}
    \\
7 
&
    \begin{array}[t]{ll}
 \Yboxdim{5pt} \yng(4,2,1) \\ \Yboxdim{5pt} \yng(4,3) 
\end{array}
    &
    \begin{array}[t]{ll}
(-1, -1)\\(0, 1)
\end{array}
    &
    \begin{array}[t]{ll}
(-4, -5)\\(8, 3)
\end{array}
    &
    \begin{array}[t]{ll}
(-8, -3)\\(-8, 3)\\(-4, -5)\\(-4, 5)\\(4, -5)\\(4, 5)\\(8, -3)\\(8, 3)
\end{array}
    \\
8 
&
 \Yboxdim{5pt} \yng(5,2,1)
    &
    \begin{array}[t]{ll}
(1, -2)
\end{array}
    &
    \begin{array}[t]{ll}
(-10, 1)
\end{array}
    &
    \begin{array}[t]{ll}
(-10, -1)\\(-10, 1)\\(10, -1)\\(10, 1)
\end{array}
    \\
9 
&
    &
    \begin{array}[t]{ll}

\end{array}
    &
    \begin{array}[t]{ll}

\end{array}
    &
    \begin{array}[t]{ll}

\end{array}
    \\
10 
&
 \Yboxdim{5pt} \yng(5,4,1)
    &
    \begin{array}[t]{ll}
(0, -2)
\end{array}
    &
    \begin{array}[t]{ll}
(-10, -3)
\end{array}
    &
    \begin{array}[t]{ll}
(-10, -3)\\(-10, 3)\\(10, -3)\\(10, 3)
\end{array}
    \\
11 
&
 \Yboxdim{5pt} \yng(6,3,2)
    &
    \begin{array}[t]{ll}
(-2, 1)
\end{array}
    &
    \begin{array}[t]{ll}
(8, -5)
\end{array}
    &
    \begin{array}[t]{ll}
(-8, -5)\\(-8, 5)\\(8, -5)\\(8, 5)
\end{array}
    \\
12 
&
 \Yboxdim{5pt} \yng(6,4,2)
    &
    \begin{array}[t]{ll}
(-2, 0)
\end{array}
    &
    \begin{array}[t]{ll}
(2, -7)
\end{array}
    &
    \begin{array}[t]{ll}
(-2, -7)\\(-2, 7)\\(2, -7)\\(2, 7)
\end{array}
    \\
13 
&
 \Yboxdim{5pt} \yng(6,4,2,1)
    &
    \begin{array}[t]{ll}
(2, -1)
\end{array}
    &
    \begin{array}[t]{ll}
(-4, 7)
\end{array}
    &
    \begin{array}[t]{ll}
(-4, -7)\\(-4, 7)\\(4, -7)\\(4, 7)
\end{array}
    \\
14 
&
 \Yboxdim{5pt} \yng(6,5,2,1)
    &
    \begin{array}[t]{ll}
(2, -2)
\end{array}
    &
    \begin{array}[t]{ll}
(-10, 5)
\end{array}
    &
    \begin{array}[t]{ll}
(-10, -5)\\(-10, 5)\\(10, -5)\\(10, 5)
\end{array}
    \\
15 
&
    &
    \begin{array}[t]{ll}

\end{array}
    &
    \begin{array}[t]{ll}

\end{array}
    &
    \begin{array}[t]{ll}

\end{array}
    \\
16 
&
 \Yboxdim{5pt} \yng(7,4,3,2)
    &
    \begin{array}[t]{ll}
(-1, 2)
\end{array}
    &
    \begin{array}[t]{ll}
(14, 1)
\end{array}
    &
    \begin{array}[t]{ll}
(-14, -1)\\(-14, 1)\\(14, -1)\\(14, 1)
\end{array}
    \\
17 
&
 \Yboxdim{5pt} \yng(8,4,3,2)
    &
    \begin{array}[t]{ll}
(1, 1)
\end{array}
    &
    \begin{array}[t]{ll}
(8, 7)
\end{array}
    &
    \begin{array}[t]{ll}
(-8, -7)\\(-8, 7)\\(8, -7)\\(8, 7)
\end{array}
    \\
18 
&
 \Yboxdim{5pt} \yng(7,6,3,2)
    &
    \begin{array}[t]{ll}
(-2, 2)
\end{array}
    &
    \begin{array}[t]{ll}
(14, -3)
\end{array}
    &
    \begin{array}[t]{ll}
(-14, -3)\\(-14, 3)\\(14, -3)\\(14, 3)
\end{array}
    \\
19 
&
    &
    \begin{array}[t]{ll}

\end{array}
    &
    \begin{array}[t]{ll}

\end{array}
    &
    \begin{array}[t]{ll}

\end{array}
    \\
\hline
\end{array}
$$
\caption{Solutions of $x^2+3y^2 = 12N+7$,
corresponding affine Grassmannian elements of type $D_4^{(3)}$,
and corresponding bar partitions in $\cD_4^\flat$.}
\label{table_12N+7}
\end{figure}

\end{Exa}

\section{Solution sets via cores: type $A_3^{(1)}$}
\label{sec_pell_A3}

In this section, we focus on type $A_3^{(1)}$
and proceed to a detailed study of the relationship between (extended) affine Grassmannian elements 
and solutions of the naturally attached Diophantine equation. 
In particular, we establish analogues of \Cref{thm_A2} and \Cref{same orbit A2} in rank $3$, this will be 
\Cref{thm size orbits} and \Cref{pig A3} respectively.
The main difference here is that the action of 
the group $G$ on the corresponding solution set will be significantly more complicated.

All along the section, $M$ is the root lattice of type $A_3^{(1)}$ and $L$ is the corresponding weight lattice.

\subsection{From $4$-cores to integer solutions}
\label{sec_A3_phi}

For $N\in\N$, recall that ${\mathcal{B}}(N) = \left\{ q\in {{M}} \mid \sL_{\La_0}(t_q) = N \right\}$ and $\widehat{\mathcal{B}}(N) = \{ q\in L \mid \sL_{\Lambda_0}(t_q) = N\}$.
Recall also from \Cref{exa_Pell_A2} (2) that we established the following identity 
for all $q = q_1\eps_1 + q_2\eps_2 + q_3 \eps_3 \in {{M}}$.
\begin{Lem}\label{lem_A3}
For all $q\in {{M}}$, we have
$$
48 \sL_{\La_0}(t_q) +30 =  (12q_2 + 4q_3 - 1)^2 + 2(8q_3+1)^2 + 3(8q_1 + 4q_2 + 4q_3 - 3)^2.
$$
\end{Lem}
Therefore, we consider the Diophantine equation 
\begin{equation}\label{eq_A3}
x^2+2y^2+3z^2 = 48N+30
\end{equation}
where $N\in\N$ is fixed, and simply denote by $\cU(48N+30)$ 
the set of integer solutions to \Cref{eq_A3}.
We have also explained in \Cref{exa_Pell_A2} (2) why the map
\begin{equation}\label{map p A3}
\begin{array}{ccccc}
\varphi & : & \widehat \cB(N) & \longrightarrow & \cU(48N+30) \\
   & & q &\longmapsto &(12q_2 + 4q_3 -1, 8q_3+1,8q_1 + 4q_2+4q_3-3)
   \end{array}
\end{equation}
is well-defined,
that is, we can attach an integral solution of \Cref{eq_A3} to any extended affine Grassmannian element.

\medskip

Let us express \Cref{map p A3} in matrix form, with respect to
the basis $\cC$ of $V_0$ and the standard basis of $\mathbb R^3$.
\begin{equation}
\label{eq:matricephi_rk3}
\varphi(q)=Q_0q+u_0,
\end{equation}
where
\begin{equation}\label{matrix B def rk3}
Q_0 = 
 \begin{pmatrix}
0 & 12 & 4  \\
0 & 0 & 8   \\
8 & 4 &4 
\end{pmatrix}
\quad \text{and}\quad 
u_0 = 
 \begin{pmatrix}
 -1\\1\\-3
\end{pmatrix}.
\end{equation}

\medskip

\subsection{Extended affine Grassmannian elements}

Similarly to \Cref{sec_Pell_A2}, we have by \Cref{L in terms of M}
\begin{equation}\label{dec L A3}
L  = {{M}} \sqcup (\varpi_1+ w_{0,1}w_0({{M}}) ) \sqcup (\varpi_2+ w_{0,2}w_0({{M}}) ) \sqcup (\varpi_3+ w_{0,3}w_0({{M}}) ).
\end{equation}
Again, we will use matrix representation for elements of the extended affine Weyl group $\widehat W$,
and we will always work with the basis $\cC$ of ${{M}}$ introduced in \Cref{basis_eps}.
From \Cref{fund_weights_A}, we see that the fundamental weights decompose in the basis $\mathcal{C}$ by
\begin{equation}\label{eq:fund_weight_A3}
\varpi_1 =
\frac{1}{4}
 \begin{pmatrix}
3 \\
-1 \\
-1 
\end{pmatrix},
 \quad
\varpi_2 =
\frac{1}{2}
 \begin{pmatrix}
1 \\
1 \\
-1 
\end{pmatrix}
\quad \text{ and } \quad
\varpi_3 =
\frac{1}{4}
 \begin{pmatrix}
1 \\
1 \\
1 
\end{pmatrix}.
\end{equation}
From \Cref{w01w0_can_A},  we deduce that the matrix of $w_{0,1}w_0$ (viewed as a linear map on $V_0$) is
\begin{equation}\label{matrix w01w0 def A3}
    B= \text{Mat}_{\cC}(w_{0,1}w_{0}) =
 \begin{pmatrix}
-1 & -1 & -1 \\
1 & 0 & 0  \\
0 & 1 & 0 
\end{pmatrix}.
\end{equation}
Moreover, as we already explained in \Cref{section Sigma type A}, 
$\mathrm{Mat}_{\cC} (w_{0,i}w_0) = B^i$ for all $i=1,2,3$.
If
we define 
$$Q_i =
Q_{0}B^i\quad\text{and}\quad u_i = Q_0\varpi_i  + u_0,$$
then for all $q\in {{M}}$, we have
  \begin{equation}\label{matrices Q_i}
  \varphi(\varpi_i+w_{0,i}w_0(q)) = Q_iq + u_i.
  \end{equation}
Writing $q=(q_1,q_2,q_3)$ one computes $\varpi_1 + w_{0,1}w_0(q)$, $\varpi_2 + w_{0,2}w_0(q)$, $\varpi_3 + w_{0,3}w_0(q)$, respectively given by
\begin{equation}
\label{eq:sigmatq_A3}
\frac{1}{4}
\begin{pmatrix}
-4(q_1+q_2+q_3)+3\\
4q_1-1\\
4q_2-1
\end{pmatrix},
\quad 
\frac{1}{2}
\begin{pmatrix}
2q_3+1\\
-2(q_1+q_2+q_3)\\
2q_1-1
\end{pmatrix},
\quad
\frac{1}{4}
\begin{pmatrix}
4q_2+1\\
4q_3+1\\
-4(q_1+q_2+q_3)+1
\end{pmatrix},
\end{equation}
which constitute, together with $(q_1,q_2,q_3)$, the typical elements of $L$ by \Cref{dec L A3}.

\subsection{Main theorem}

We have just seen how to attach a solution of  \Cref{eq_A3}
to any extended affine Grassmannian element.
We now turn to the converse: which solutions arise from extended
affine Grassmannian elements? \Cref{pig A3} below asserts that this is
always the case, up to the action of a suitable finite group.

\medskip

More precisely, in line with the spirit of the matrix $R$ introduced in
\cite[Section 2]{BrunatNath2022} and recalled in \Cref{matrix R BN},  we
define $R$ to be the following matrix
\begin{equation}
\label{eq:matR}
R :=  
 \begin{pmatrix}
1 & 0 & 0\\
0 & \frac{1}{\sqrt{2}} & 0 \\
0 & 0 & \frac{1}{\sqrt{3}}
\end{pmatrix}
 \begin{pmatrix}
\text{cos}(\frac{\pi}{3}) & 0 & -\text{sin}(\frac{\pi}{3})\\
0 & 1 & 0 \\
\text{sin}(\frac{\pi}{3}) & 0 & \text{cos}(\frac{\pi}{3})
\end{pmatrix}
 \begin{pmatrix}
1 & 0 & 0\\
0 & \sqrt{2}  & 0 \\
0 & 0 & \sqrt{3} 
\end{pmatrix}
=  \frac{1}{2}\begin{pmatrix}
1 & 0 & -3\\
0 & 2 & 0 \\
1 & 0 & 1 
\end{pmatrix}.
\end{equation}

By construction, the matrix $R$ is 
an element of the orthogonal group $\mathcal{O}_f(\mathbb{R})$ of the quadratic form $f(x,y,z) = x^2 + 2y^2 + 3z^2$, 
implying in particular that $f(RX) = f(X)$ for any $X \in \mathbb{R}^3$.
We introduce now the following group
\begin{align}\label{def group G rank 3}
    G := \langle R, s \rangle \leq \mathcal{O}_f(\mathbb{R})
\end{align}
where $s$ is the matrix, representing the reflection through the plane spanned by $e_1$ and
$e_2$. In particular, $G$ is isomorphic to the dihedral group of order
$12$, and we have the relations 
\begin{equation}\label{relation G rank 3}
   sR^ks = R^{-k} \quad \text{for all } k. 
\end{equation}

\begin{Prop}\label{prop xyz properties and action of G on the solutions}
Let $X =(x,y,z) \in \mathcal{U}(48N + 30)$.  Then
$y$ is odd,
 $x$ and $z$ have the same parity and the group $G$ acts on $\mathcal{U}(48N + 30)$.
\end{Prop}

\begin{proof}
Reducing $x^2 + 2y^2 + 3z^2 = 48N+30$ modulo $2$ we get $x^2 = z^2 ~\text{mod} ~2$, which is equivalent to $x = z ~\text{mod} ~2$. Hence $x$ and $z$ have the same parity. Therefore, $x^2 = z^2 ~\text{mod}~ 4$ and reducing $x^2 + 2y^2 + 3z^2 = 48N+30$ modulo $4$ we get $2y^2 = 2 ~\text{mod}~ 4$, implying that $y^2$ is odd, and then $y$ as well.

Since $G$ is a subgroup of $\mathcal{O}_f(\mathbb{R})$, $G$ stabilizes the surface defined by the equation $x^2 + 2y^2 + 3z^2 = 48N+30$. But now $G$  also stabilizes the integral points of this surface, namely $\mathcal{U}(48N + 30)$. To prove it itsuffices to show that $RX$ and $sX$ belong to $\mathcal{U}(48N + 30)$ whenever $X$ does. For $sX$ this is obvious. For $RX$ this comes from the computation $RX = \frac{1}{2}(x-3z,2y,x+z)$ and from the fact that $x$ and $z$ have the same parity.
\end{proof}

We are ready to state the main theorem of this section.

\begin{Th}\label{pig A3}
Let $X = (x,y,z) \in \cU(48N+30)$. 
There exists  $r \in L$ and $g\in G$
such that $g\varphi(r) = X$. 
\end{Th}

\bigskip

This result can be viewed as a rank $3$ analogue of \cite{BrunatNath2022}.
It establishes a connection between the $G$-orbits of $\mathcal U(48N+30)$
and the set of extended Grassmannian elements of atomic length $N$.

\bigskip

\begin{Rem}
\label{rem_thm_pig}
We postpone the proof of \Cref{pig A3} to  \Cref{proof pig} as it is rather technical.
However, we now outline the key ideas underlying our approach,
and illustrate it in \Cref{fig:process_table}.
\begin{enumerate}[(1)]
\item 
Our objective is to find $r\in L$ and $g\in G$ such that
$g\varphi(r)=X$. By Part (2) of \Cref{L in terms of M},
this is equivalent to finding $1\leq i\leq 3$ and $q\in {{M}}$ such that 
\begin{equation}
\label{eq:loc1}
g\varphi(\varpi_i+w_{0,i}w_0(q))=X.
\end{equation}
Now, note that, for fixed $g\in G$ and $1\leq i\leq 3$, equation
\Cref{eq:loc1} (with unknown $q$) can be solved in $V_0$. To proceed, we
express this equation in matrix form. More precisely, using \Cref{matrices Q_i}, 
we obtain
\begin{align}\label{eq:eq}
    g \varphi(\varpi_i + w_{0,i}w_0(q)) = X \quad & 
    \Longleftrightarrow \quad Q_i q + u_i = g^{-1}X \nonumber\\
    &\Longleftrightarrow \quad q = Q_i^{-1}(g^{-1}X - u_i) \nonumber\\
    &\Longleftrightarrow \quad q = B^{-i}\left(Q_0^{-1}(g^{-1}X - u_0) -
    \varpi_i\right).
\end{align}

\item 
Define a bijection
$\tau:\{0,\,1,\,2,\,3\}\rightarrow \{1,\,3,\,5,\,7\}$ by setting
\begin{equation}\label{tau}
\tau(0)=1,\quad\tau(1)=7,\quad \tau(2)=5\quad\text{and}\quad \tau(3)=3.
\end{equation}
Write $q=(q_1,q_2,q_3)\in {{M}}$. First, we observe that, for all $g\in G$,
the second coordinate of $g\varphi(\varpi_i+w_{0,i}w_0(q))$ and
$\varphi(\varpi_i+w_{0,i}w_0(q))$ are equal. Moreover, using
\Cref{matrices Q_i}, we can compute the value of this coordinate for $i=0,1,2,3$, respectively given by 
$$8q_3+1,\quad 8q_2-1,\quad 8q_1-3\quad
\text{and}\quad-8q_1-8q_2-q_3+3.$$
Reducing these expressions modulo $8$, we conclude that the second
coordinate of $g\varphi(\varpi_i+w_{0,i}w_0(q))$ has residue $\tau(i)$. 
\item Let $\beta$ be the residue of $y$ modulo $8$. By
\Cref{prop xyz properties and action of G on the solutions}, we know that
$\beta\in \{1,3,5,7\}=\operatorname{Im}(\tau)$. Then we define
\begin{equation}\label{i_beta}
i_{\beta}=\tau^{-1}(\beta).
\end{equation}

Therefore, for any $q \in {{M}}$ the second coordinate of $\varphi\big(\varpi_{i_{\beta}}+w_{0,i_{\beta}}w_0(q)\big)$ modulo 8 is equal to $\beta$. In particular, if $\varphi(\varpi_i + w_{0,i}w_0(q)) = g\varphi\big(\varpi_j+ w_{0,j}w_0(q')\big)$ for some  $q,q' \in {{M}}$, $g \in G$ and $i,j \in \{0,1,2,3\}$ then they have the same second coordinate, forcing $i = j$.

For this specific value $i=i_{\beta}$, the equation \Cref{eq:eq} admits a
unique solution in $V_0$,
which depends only on $X$ and $g$. We denote this solution by $r(X,g)$. 
\item Our goal is then to construct an explicit element $g\in G$
such that $r(X,g)\in {{M}}$. This final step is highly technical and requires
a case-by-case analysis. It highlights the significantly greater
complexity of the rank $3$ case compared to rank~$2$. 
\end{enumerate}

\begin{figure}[h!]
\label{fig:process_table}
    \centering
    \includegraphics[scale = 0.42]{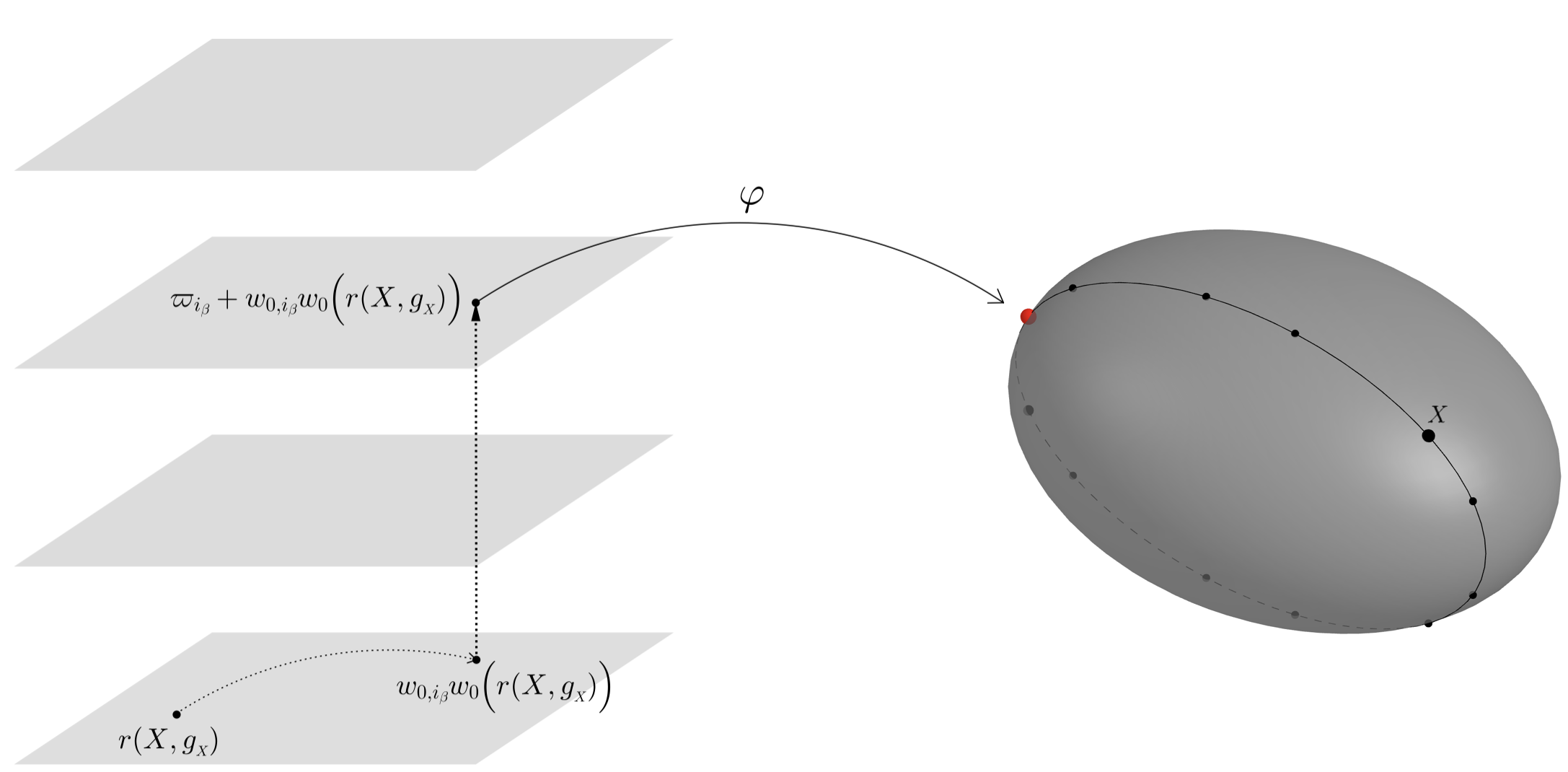}
    \caption{Let $X = (x,y,z)$ be an integral solution of the equation $x^2 + 2y^2 + 3z^2 = 48N+30$ for $N$ fixed. This picture shows the general idea of building the socle $r:=r(X,g_{\raisebox{-0.5ex}{\scalebox{0.55}{$X$}}})$ so that it belongs to the lattice ${{M}}$ and such that $g_{\raisebox{-0.5ex}{\scalebox{0.55}{$X$}}}\varphi\Big(\varpi_{i_\beta}+w_{0,i_{\beta}}w_0(r)\Big) = X$. The 4-core corresponding to $r$ via \Cref{bij_lascoux} is of size $N$. The floor on the left-hand side, where $\varpi_{i_\beta} + w_{0, i_\beta} w_0(r)$ lies, is determined by the value $i_\beta \in {0, 1, 2, 3}$, as introduced in \Cref{i_beta}. For example in this picture, $i_{\beta} = 2$.}
\end{figure}
\end{Rem}

\subsection{An application: construction of 4-cores of a given size}\label{application 4-cores}
Recall from \Cref{sec_cores} that type $A_{n}^{(1)}$ affine Grassmannian elements (or equivalently, elements $q\in {{M}}$) are in bijection with $(n+1)$-cores.
In this section, we relate  the 4-cores of size $N$ to the integral solutions of the quadratic equation 
$x^2+2y^2+3z^2 = 48N+30$.
We will use the $G$-action studied in the previous section, and  in particular \Cref{pig A3}, to prove the following theorem,
which can be seen as a constructive proof of the Granville-Ono theorem \cite{GO1996} in rank $3$.
Therefore, our result naturally falls within the framework of the already thorough study of the $4$-cores,
see \cite{Ono1994, Ono1995, HirschhornSellers1996a, HirschhornSellers1996b, GO1996, OnoSze1997, HirschhornSellers1999}.

\begin{Th}\label{thm 4-cores}
For any $N \in \mathbb{N}$ we can construct a 4-core of size $N$ from any integral solution of the equation $x^2 + 2y^2 + 3z^2 = 48N+30$.
\end{Th}

\begin{proof}
We split the proof into two steps.
\begin{itemize}
    \item 
Firstly, in \cite{bhargava2000conway} Bhargava characterises the integers not represented by the quadratic form $x^2+2y^2+3z^2$ (see \cite[Table 1]{bhargava2000conway}), namely
\begin{equation}
\{n \in \mathbb{N}~|~ \text{for all } (x,y,z) \in \mathbb{Z}^3, x^2+2y^2+3z^2 \neq n \} 
 = \{2^d(8k+5)~|~d~\text{odd},~ k \in \mathbb{N}\}.
\end{equation}
But for any $N \in \mathbb{N}$ the integer $48N+30$ is not in the above set. Indeed, assume otherwise, that is suppose that there exists an odd integer $d$ and $k \in \mathbb{N}$ such that $48N+30 = 2^d(8k+5)$.  Then $d = 1$ and it follows that $24N + 15 = 8k+5$, which is impossible after reducing modulo 8.
Therefore,  for any choice of $N \in \mathbb{N}$  there exists $(x,y,z) \in \mathbb{Z}^3$ such that $x^2 + 2y^2 + 3z^2 = 48N+30$.

\item
Secondly, we construct a $4$-core of size $N$ in terms of an integral solution of the equation $x^2 + 2y^2 + 3z^2 = 48N+30$. Let $X \in \cU(48N+30)$ be such a solution. By \Cref{pig A3} there exists $g \in G$ and 
$r \in L$, that we can explicitly construct from $X$ (see \Cref{Table_parametrisation}), such that $g\varphi(r) = X$, that is $\varphi(r) = g^{-1}X$, which also belongs to $\mathcal{U}(48N+30)$ by \Cref{prop xyz properties and action of G on the solutions}.

Setting $q = w_{0,i}w_0^{-1} (r-\varpi_i)$,
where $i \in \{0,1,2,3\}$ is determined by \Cref{i_beta},
we have $q\in {{M}}$ by \Cref{L in terms of M}, 
and $\varphi(q)\in \cU(48N+30)$ by \Cref{new_sols_A} since $\varphi(r)\in \cU(48N+30)$.
We conclude by using the explicit bijection between ${{M}}$ and the $4$-cores explained in \Cref{sec_cores_alcoves}.
\end{itemize}
\end{proof}

\begin{Exa}[\Cref{thm 4-cores} in action]\label{example 4-cores}
    Let us take $N=57$ and let us consider the set of integral solutions of the equation $x^2 + 2y^2 + 3z^2 = 48\cdot 57 + 30 = 2766$, namely $\mathcal{U}(2766)$. An easy computational verification with SageMath shows that $|\mathcal{U}(2766)| = 590$. We treat below two solutions: $X = (-37, 25, -7)$ and $Y = (-5, 37, 1)$, where from each one of them we build a 4-core of size 57.

    \begin{itemize}
        \item[(1)] \underline{From $X$}. We have $-37 = 12\cdot(-4) + 11$, $25 = 8\cdot 3 + 1$ and $-7 = 4\cdot(-2)+1$, so that $k = -4$, $n=3$, $m = -2$ and $\al = 11$, $\beta = 1$, $\ga = 1$. Hence $X \mod H = (11,1,1)$ and we need to look at the third "large" row of \Cref{Table_parametrisation} for $\beta = 1$. Finally, by \Cref{i_beta} we have $i_{\beta} = 0$. We have then all the ingredients to compute the 4-th column of \Cref{Table_parametrisation} corresponding to $X$.

    Since $k,m$ have same parity and $n \mod 3 = 0$, we need to take the first row of the small two rows in this big row, and then the corresponding $g$ is 1. Thus $r_0(X,1) \in {{M}}$ and is expressed in the basis $\mathcal{C}$ by
    \begin{align*}
        r_0(X,1) &= \frac{1}{24}\Big(-(-37)-25+3\cdot(-7)+9,~ 2\cdot(-37)-25+3,~ 3\cdot 25-3\Big) \\
        & = \frac{1}{24}\big(0,~-96,~72\big) \\
        & = 0\varepsilon_1 -4\varepsilon_2 + 3\varepsilon_3.
    \end{align*}

    Now, to find the corresponding 4-core it is more convenient to express this vector in the standard basis: $r_0(X,1) = (0,-4,3,1)$. Thus, the corresponding 4-core of size 57 (via \Cref{bij_lascoux}) is given by $(11,8,6,6,4,4,3,3,3,2,2,2,1,1,1)$.

    \item[(2)] \underline{From $Y$}. We have $-5 = 12\cdot(-1) + 7$, $37 = 8\cdot 4 + 5$ and $1 = 4\cdot 0 +1$, so that $k = -1$, $n=4$, $m = 0$ and $\al = 7$, $\beta = 5$, $\ga = 1$. Hence $Y \mod H =(7,5,1)$ and we need to look at the 9-th "large" row of \Cref{Table_parametrisation} for $\beta = 5$, giving $i_{\beta} = 2$ by \Cref{i_beta}. Since $n \mod 3 = 1$, by \Cref{abc} the right $g\in G_{(7,5,1)}$ is $g = Rs$. Thus $r(5,Rs)(Y) \in {{M}}$ and 
    \begin{align*}
        r(5,Rs)(Y) &= \frac{1}{24} B^{-2}\Big((-5-37-3\cdot 1 + 9,~-5 - 37 + 3\cdot 1 + 3,~3\cdot 37 - 2)   -24\varpi_2 \Big) \\
        &= \frac{1}{24} B^{-2}(-48,~-48,~120) \\
        & = 5\varepsilon_1 -\varepsilon_2 -2\varepsilon_3.
    \end{align*} 

Similarly, $r(5,Rs)(Y) = (5,-1,-2,-2)$ in the standard basis and its corresponding 4-core (via \Cref{bij_lascoux}) is given by $(17,14,11,8,5,2)$.

    \end{itemize}

\end{Exa}

\begin{Rem}
Observe that the quadratic form $x^2+2y^2+3z^2$ just fails to be \textit{universal}, that is, to represent all integers.
Indeed, the so-called \say{Fifteen Theorem}, announced by Conway and Schneeberger in 1993 (a sketch of proof can be found in \cite{Schneeberger1995} and a proof in \cite{bhargava2000conway}),
states that if an integral positive definite quadratic form represents every positive integer up to 15, then it represents every positive integer. 
Our quadratic form  $x^2+2y^2+3z^2$ represents all integers up to $15$ but $10$, so the theorem does not apply.
Interestingly, in 1995 Kaplansky \cite{kaplansky1995ternary} showed that the quadratic form $x^2+2y^2+3z^2$ is \textit{$O$-universal}, that is, represents every odd positive integer (obviously though, $48N+30$ is never odd).
\end{Rem}

\subsection{Proof of  \Cref{pig A3}}\label{proof pig}

We are now concerned with the proof of the main theorem of this section.
This will be achieved in two steps:
\begin{enumerate}
    \item First, we consider only the residue classes of solutions $(x,y,z)$ modulo $12,8,4$ respectively.
    \item Second, we  completely characterise how solutions are obtained from extended affine Grassmannian elements and the $G$-action, see \Cref{lemma fonda table} and \Cref{Table_parametrisation}.
\end{enumerate}

\subsubsection{Residue classes of solutions}\label{section_arbre}

We are now interested in studying the integral solutions of \Cref{eq_A3}.
First of all, we 
give general necessary conditions for an integer triplet to be in  $\cU(48N+30)$ .

We now consider
\begin{equation}\label{group H}
    H= 12\mathbb{Z} \times 8\mathbb{Z} \times 4\mathbb{Z}.
\end{equation}

For $(x,y,z) \in \mathcal{U}(48N+30)$, we write $x = 12k+\alpha$, $y = 8n+\beta$ and $z = 4m + \gamma$ where $\alpha \in \llbracket 0,11 \rrbracket$ ,
 $\beta \in \llbracket 0,7 \rrbracket$ and $\gamma \in \llbracket 0,3\rrbracket$, so that $(x,y,z) \mod H = (\alpha,\beta,\gamma)$. 
By \Cref{prop xyz properties and action of G on the solutions}, $\alpha$ and $\gamma$ have same parity, and by \Cref{rem_thm_pig} (2), $\beta = 1,3,5,7$. 
We will now give a natural classification of 
the classes of $(x,y,z)$ modulo $H$, see \Cref{prop neat decomposition} below.
This will be crucial for establishing the main results of this chapter.

\medskip

\newcommand{\fB}{\mathfrak{B}}

Starting with a solution $(x,y,z)\in\cU(48N+30)$, we will consider the different cases represented in \Cref{fig:cases} and we will write $(x,y,z) \in \mathfrak{B}_i$ if we are in Case $\mathfrak{B}_i$.

\medskip

\begin{figure}[h!]
    \centering
\begin{tikzpicture}[
  level distance=1.5cm,
  level 1/.style={sibling distance=4cm},
  level 2/.style={sibling distance=2cm},
  level 3/.style={sibling distance=1cm}
]
\node {} 
  child {node {$\beta=1,7$} 
    child {node {$C_1$} 
      child {node {$
      \begin{array}[t]{c}
      0,2
      \\
      \fB_1
      \end{array}
      $}}
      child {node {$
      \begin{array}[t]{c}
      1
      \\
      \fB_2
      \end{array}
      $}}
    }
    child {node {$C_2$}
      child {node {$
      \begin{array}[t]{c}
      0,2
      \\
      \fB_3
      \end{array}
      $}}
      child {node {$
      \begin{array}[t]{c}
      1
      \\
      \fB_4
      \end{array}
      $}}
    }
  }
  child {node {$\beta=3$}
    child {node {$C_1$}
      child {node {$
      \begin{array}[t]{c}
      2,1
      \\
      \fB_5
      \end{array}
      $}}
      child {node {$
      \begin{array}[t]{c}
      0
      \\
      \fB_6
      \end{array}
      $}}
    }
    child {node {$C_2$}
      child {node {$
      \begin{array}[t]{c}
      2,1
      \\
      \fB_7
      \end{array}
      $}}
      child {node {$
      \begin{array}[t]{c}
      0
      \\
      \fB_8
      \end{array}
      $}}
    }
  }
  child {node {$\beta=5$}
    child {node {$C_1$}
      child {node {$
      \begin{array}[t]{c}
      1,0
      \\
      \fB_{9}
      \end{array}
      $}}
      child {node {$
      \begin{array}[t]{c}
      2
      \\
      \fB_{10}
      \end{array}
      $}}
    }
    child {node {$C_2$}
      child {node {$
      \begin{array}[t]{c}
      1,0
      \\
      \fB_{11}
      \end{array}
      $}}
      child {node {$
      \begin{array}[t]{c}
      2
      \\
      \fB_{12}
      \end{array}
      $}}
    }
  };
\end{tikzpicture}
\caption{The different cases, 
where  \\
$C_1$ : \;($\alpha, \gamma$ are even) or ($\alpha$, $\gamma$ are odd and $k,m$ have same parity)  \\
$C_2$ : \; ($\alpha$, $\gamma$ are odd and $k,m$ have different parity). \\
The bottom values indicate the residue of $n$ modulo $3$. 
}

\label{fig:cases}
\end{figure}

\begin{Lem}\label{lemma 12x8x4}
Let $(x,y,z) \in \cU(48N+30)$.  Then
$(x^2 + 2y^2 + 3z^2) - (\alpha^2 + 2\beta^2 + 3\gamma^2)$ equals
\begin{itemize}
\item $0 \mod 48$ if we are in Cases $\fB_1$, $\fB_6$, $\fB_9$,
\item $8 \mod 48$ if we are in Case $\fB_7$,
\item $16 \mod 48$ if we are in Cases $\fB_2$, $\fB_{10}$,
\item $24 \mod 48$ if we are in Cases $\fB_3$, $\fB_8$, $\fB_{11}$,   
  \item $32 \mod 48$ if we are in Case $\fB_5$,
\item $40 \mod 48$ if we are in Cases $\fB_4$, $\fB_{12}$.
\end{itemize}
\end{Lem}

\begin{proof}
One has
\begin{align*}
x^2 + 2y^2 + 3z^2 & = (12k+\alpha)^2 + 2(8n+\beta)^2 + 3(4m+\gamma)^2 \\
& = 144k^2 + 24k\alpha + \alpha^2 + 128n^2 + 32\beta n + 2\beta^2 + 48m^2 + 24m\gamma + 3\gamma^2 \\
& = 48(3k^2 + 2n^2+ m^2) + (\alpha^2 + 2\beta^2 + 3\gamma^2) + 32n(n+\beta) + 24(k\alpha+ m\gamma) \\
& =  (\alpha^2 + 2\beta^2 + 3\gamma^2) + 24(k\alpha + m\gamma) + 32n(n+\beta) ~\mod 48.
\end{align*}

We need now to analyse the two terms $24(k\alpha+ m\gamma)$ and $32n(n+\beta)$ modulo 48.  
\begin{itemize}
    \item  By \Cref{prop xyz properties and action of G on the solutions} we know that $\alpha$ and $\gamma$ have same parity. 
If $\alpha$ and $\gamma$ are even then $24(k\alpha+ m\gamma) = 0 ~\mod 48$.  If $\alpha$ and $\gamma$ are odd, then $24(k\alpha+ m\gamma) ~\mod 48 = 24(k+m) ~\mod 48$,   which gives $24(k+m) = 24 ~\mod 48$ if $k$ and $m$ do not have same parity and $24(k+m) = 0 ~\mod 48$ otherwise.
\item As $(x,y,z) \in \mathcal{U}(48N+30)$ we know that  $\beta$ belongs to $\{1,3,5,7\}$.  We have then 4 cases to consider but in fact,  the cases $\beta = 1$ and $\beta = 7$ will give the same answer.  If $\beta = 1,7$ then a direct computation gives for $n = 0 \mod 3$ that $32n(n+1) = 0 \mod 48$, for $n = 1 \mod 3$ that $32n(n+1) = 16 \mod 48$, for $n = 2 \mod 3$ that $32n(n+1) = 0 \mod 48$. If $\beta = 3$ then again a direct computation gives for $n = 0 \mod 3$ that $32n(n+3) = 0 \mod 48$, for $n = 1 \mod 3$ that $32n(n+3) = 32 \mod 48$, for $n = 2 \mod 3$ that $32n(n+3) = 32 \mod 48$. Finally, if $\beta = 5$, for $n = 0 \mod 3$ that $32n(n+5) = 0 \mod 48$, for $n = 1 \mod 3$ that $32n(n+5) = 0 \mod 48$, for $n = 2 \mod 3$ that $32n(n+5) = 16 \mod 48$.
\end{itemize}
Combining the two points above we obtain the announced statements.
\end{proof}

\bigskip

For $\beta=1,3,5,7$, we now introduce the following subsets of $\Z^3$:
    \begin{align*}
        \Omega_{hex}(\beta) & := \Bigl\{(1, \beta, 3),(4, \beta, 2),(7, \beta, 3),(11,\beta,1),(8, \beta, 2),(5,\beta,1)\Bigr\} \\
        \Omega_{squ}(\beta) & := \Bigl\{(1,\beta,1),(5,\beta,3),(7,\beta,1),(11,\beta,3)\Bigr\} \\
        \Omega_{tri}(\beta) & := \Bigl\{(0,\beta,2),(3,\beta,1),(9,\beta,3)\Bigr\} \\
        \Omega_{seg}(\beta) & := \Bigl\{(3,\beta,3),(9,\beta,1)\Bigr\}
    \end{align*}
and we set 
\begin{align}\label{OMEGA}
    \Omega(\beta) := \Omega_{hex}(\beta)\sqcup  \Omega_{squ}(\beta)\sqcup \Omega_{tri}(\beta) \sqcup  \Omega_{seg}(\beta)
\text{\quad along with \quad }
\Omega=\bigsqcup_{\beta=1,3,5,7}\Omega(\beta).
\end{align}

\bigskip

\begin{Prop}\label{prop neat decomposition}
Let $(x,y,z)\in\cU(48N+30)$. Write $ (\alpha,\beta,\gamma) = (x,y,z) \mod H$. We have 
\begin{itemize}
    \item[(i)] $(x,y,z) \in \mathfrak{B}_1 \sqcup \mathfrak{B}_5 \sqcup \mathfrak{B}_9$ if and only if $(\alpha,\beta,\gamma) \in \Omega_{hex}(\beta)$.
    \item[(ii)] $(x,y,z) \in \mathfrak{B}_3 \sqcup \mathfrak{B}_7 \sqcup \mathfrak{B}_{11}$ if and only if $ (\alpha,\beta,\gamma) \in \Omega_{squ}(\beta)$.
    \item[(iii)] $(x,y,z) \in \mathfrak{B}_2 \sqcup \mathfrak{B}_6 \sqcup \mathfrak{B}_{10}$ if and only if $ (\alpha,\beta,\gamma) \in \Omega_{tri}(\beta)$.
    \item[(iv)] $(x,y,z) \in \mathfrak{B}_4 \sqcup \mathfrak{B}_8 \sqcup \mathfrak{B}_{12}$ if and only if $ (\alpha,\beta,\gamma) \in \Omega_{seg}(\beta)$.
\end{itemize}

\end{Prop}

\begin{proof}
Using \Cref{lemma 12x8x4},
we observe that if $(x,y,z)\in\fB_1\sqcup\fB_5\sqcup\fB_9$, then $\alpha^2+3\gamma^2 = 30 -2\beta^2-R$
where $R=0$ if $(x,y,z)\in\fB_1\sqcup\fB_9$ and $R=32$ if $(x,y,z)\in\fB_5$.
This yields $\alpha^2+3\gamma^2=28 \mod 48$ in each case.
By a straightforward computer calculation, the triples satisfying the previous identity are
$$\Bigl\{(1, \beta, 3),(4, \beta, 2),(7, \beta, 3),(11,\beta,1),(8, \beta, 2),(5,\beta,1)\Bigr\}$$
and we recognize $\Omega_{hex}(\beta)$.
By a similar argument, we find the 
sets
$$
\begin{array}{ll}
        \Bigl\{ (2,\beta,  0), (10,\beta,  0) (1,\beta,1),(5,\beta,3),(7,\beta,1),(11,\beta,3)\Bigr\}  
        & \text{\quad if } (x,y,z)\in\fB_3\sqcup\fB_7\sqcup\fB_{11},
        \\
        \Bigl\{(0,\beta,2),(3,\beta,1),(9,\beta,3)\Bigr\}  
        & \text{\quad if } (x,y,z)\in\fB_2\sqcup\fB_6\sqcup\fB_{10},
        \\
        \Bigl\{ (6,\beta,  0), (3,\beta,3),(9,\beta,1)\Bigr\} 
        & \text{\quad if } (x,y,z)\in\fB_4\sqcup\fB_8\sqcup\fB_{12}
    \end{array}
$$
and we recognize the sets $ \Omega_{squ}(\beta), \Omega_{tri}(\beta), \Omega_{seg}(\beta)$
up to the presence of $(2,\beta,  0), (6,\beta,  0), (10,\beta,  0)$ in the first and last sets.
We remark that these are exactly those triples with $\gamma=0$.
Let us show that this condition on $\gamma$ leads to a contradiction.
Note first that $\gamma=0$ means that $z=0\mod 4$, and thus $3z^2=0\mod 48$. 
Our equation therefore simplifies to $x^2+2y^2=30 \mod 48.$
Computing the quadratic residues modulo $48$ (there are $8$ possible values), 
we check by computer that the above equation has no solution, whence the contradiction.
At this point, we have proved the four direct implications. Each reverse implication follows directly by contraposition.
\end{proof}

\subsubsection{Points induced by Grassmannians}\label{section calcul matrices A3}

We will prove \Cref{pig A3} constructively, by establishing \Cref{Table_parametrisation} (on Page \pageref{Table_parametrisation}).
First of all, let us explain how to read \Cref{Table_parametrisation}.
Start with an integer solution $X\in \cU(48N+30)$.

\begin{itemize}
\item \textbf{Column 1:} Use \Cref{prop neat decomposition} to 
list all possible classes $(\al,\beta,\ga)$ of $X$ modulo $H$: $(\al,\beta,\ga) = X \mod H$.
These were given in \Cref{OMEGA}.

\item \textbf{Columns 2 and 3:}
Write $x = 12k+\alpha$, $y = 8n+\beta$ and $z = 4m + \gamma$ where $\alpha \in \llbracket 0,11 \rrbracket$, $\beta \in \llbracket 0,7 \rrbracket$ and $\gamma \in \llbracket 0,3\rrbracket$.
List the possible parities of $k,m$ and consider the residue of $n$ modulo $3$ using \Cref{fig:cases}. Accordingly, we set
\begin{equation}\label{abc}
(\mathfrak{a},\mathfrak{b},\mathfrak{c}) =
\left\{ 
\begin{array}{ll}
(0,2,1) & \text{ if } \beta=1,7 
\\
(2,1,0) & \text{ if } \beta=3 
\\
(1,0,2) & \text{ if } \beta=5
\end{array}
\right.
\end{equation}   
    
\item \textbf{Column 4:}
We give elements $g_{\raisebox{-0.5ex}{\scalebox{0.55}{$X$}}} \in G$.
At first sight, these seem to come out of nowhere, but  
we will prove in \Cref{lemma fonda table} that these will be precisely the elements $g$ appearing in \cref{pig A3}.
Note that there will always be either one or two such $g$'s.
We denote $G_{(\alpha,\beta,\gamma)} =  \{g_{\raisebox{-0.5ex}{\scalebox{0.55}{$X$}}} ~|~ (\al,\beta,\gamma) = X \mod H\}$.

\item \textbf{Column 5:} 
We compute the socle $r(X,g_{\raisebox{-0.5ex}{\scalebox{0.55}{$X$}}})$ 
as introduced in \Cref{rem_thm_pig} (1), and we write down $24 r(X,g_{\raisebox{-0.5ex}{\scalebox{0.55}{$X$}}})$.
\end{itemize}

\begin{table}[h!]
{
\centering
{\footnotesize
\begin{longtable}{lllll}
\label{Table_parametrisation}
\\
$(\alpha,\beta,\gamma) \in \Omega$~ & Parity $k,m$ ~& $n \mod 3$ & $g_{\raisebox{-0.5ex}{\scalebox{0.55}{$X$}}}$ &  Vector \textrm{$24 r(X,g_{\raisebox{-0.5ex}{\scalebox{0.55}{$X$}}}) = (r_1,r_2,r_3)$} in the basis $\mathcal{C}$ \\
\hline
$(4,\beta,2)$  &  diff.  &   $\mathfrak{a}$&  $R$     &   $B^{-i_{\beta}}\big( (-2x-y+9,~x-y+3z+3,~3y-3) -24\varpi_{i_{\beta}}\big)$             \\
               &  diff.  &  $\mathfrak{b}$ &  $R^2s$  &   $B^{-i_{\beta}}\big((2x-y+9,~ -x-y+3z+3,~ 3y-3)-24\varpi_{i_{\beta}}\big)$                   \\
               &  same   &  $\mathfrak{a}$ &  $R^5s$  &   $B^{-i_{\beta}}\big((-2x-y+9,~ x-y-3z+3,~3y-3)-24\varpi_{i_{\beta}}\big)$     \\
               &   same   & $\mathfrak{b}$  &  $R^4$   &   $B^{-i_{\beta}}\big((2x-y+9,~-x-y-3z + 3,~3y-3)-24\varpi_{i_{\beta}}\big)$                   \\
                         
\hline
$(7,\beta,3)$  &  same   & $\mathfrak{a}$  & $R^5$ &  $B^{-i_{\beta}}\big((x-y+3z+9, ~x-y-3z+3,~3y-3)-24\varpi_{i_{\beta}}\big)$                      \\
               &   same  &  $\mathfrak{b}$ & $s$   &   $B^{ -i_{\beta}}\big((-x-y-3z+9, ~2x-y+3,~3y-3)-24\varpi_{i_{\beta}}\big)$             \\
\hline

$(11,\beta,1)$ &  same   & $\mathfrak{a}$     & $1$  &   $B^{-i_{\beta}}\big((-x-y+3z+9, ~2x-y+3,~3y-3)-24\varpi_{i_{\beta}}\big)$             \\
               & same    &   $\mathfrak{b}$   & $Rs$ &   $B^{-i_{\beta}}\big((x-y-3z+9, ~x-y+3z+3,~3y-3)-24\varpi_{i_{\beta}}\big)$           \\
\hline

$(8,\beta,2)$  &  diff.   &  $\mathfrak{b}$    & $R$    &   $B^{-i_{\beta}}\big((-2x-y+9,~x-y+3z+3,~3y-3)-24\varpi_{i_{\beta}}\big)$                           \\
               &  diff.   &   $\mathfrak{a}$   & $R^2s$ &   $B^{-i_{\beta}}\big((2x-y+9,~ -x-y+3z+3,~3y-3)-24\varpi_{i_{\beta}}\big)$                       \\
               &  same    &  $\mathfrak{b}$    & $R^5s$ &   $B^{-i_{\beta}}\big((-2x-y+9,~ x-y-3z+3,~3y-3)-24\varpi_{i_{\beta}}\big)$                   \\
               &  same    & $\mathfrak{a}$     & $R^4$  &   $B^{-i_{\beta}}\big((2x-y+9,~-x-y-3z + 3,~3y-3)-24\varpi_{i_{\beta}}\big)$                       \\
\hline

$(5,\beta,1)$  &  same  &$\mathfrak{a}$      & $R^2$  & $B^{-i_{\beta}}\big((-x-y-3z+9,~-x-y+3z+3,~3y-3)-24\varpi_{i_{\beta}}\big)$                          \\
               &  same  & $\mathfrak{b}$     & $R^3s$ & $B^{-i_{\beta}}\big((x-y+3z+9,~-2x-y+3,~3y-3)-24\varpi_{i_{\beta}}\big)$                          \\
               
\hline

$(1,\beta,3)$  &  same  & $\mathfrak{a}$      & $R^3$   & $B^{-i_{\beta}}\big((x-y-3z+9, ~-2x-y+3,~3y-3)-24\varpi_{i_{\beta}}\big)$                            \\
               &  same  & $\mathfrak{b}$      & $R^4s$  & $B^{-i_{\beta}}\big((-x-y+3z+9, ~-x-y-3z+3,~3y-3)-24\varpi_{i_{\beta}}\big)$                          \\
               
\hline
$(1,\beta,1)$    &  diff.  & $\mathfrak{b}$     & $R^2$  & $B^{-i_{\beta}}\big((-x-y-3z+9,~-x-y+3z+3,~3y-3)-24\varpi_{i_{\beta}}\big)$                  \\
                 &  diff.  &  $\mathfrak{a}$    & $R^3s$ & $B^{-i_{\beta}}\big((x-y+3z+9,~-2x-y+3,~3y-3)-24\varpi_{i_{\beta}}\big)$                    \\
\hline

$(5,\beta,3)$    & diff.   & $\mathfrak{b}$     & $R^3$  & $B^{-i_{\beta}}\big((x-y-3z+9, ~-2x-y+3,~3y-3)-24\varpi_{i_{\beta}}\big)$                         \\
                 & diff.   & $\mathfrak{a}$     & $R^4s$ & $B^{-i_{\beta}}\big((-x-y+3z+9, ~-x-y-3z+3,~3y-3)-24\varpi_{i_{\beta}}\big)$                        \\
                 
\hline
$(7,\beta,1)$  & diff.   & $\mathfrak{b}$      &  $1$  &   $B^{-i_{\beta}}\big((-x-y+3z+9, ~2x-y+3,3y-3)-24\varpi_{i_{\beta}}\big)$                  \\
               & diff.   & $\mathfrak{a}$      & $Rs$  &   $B^{-i_{\beta}}\big((x-y-3z+9, ~x-y+3z+3,~3y-3)-24\varpi_{i_{\beta}}\big)$                     \\
               
\hline
$(11,\beta,3)$            &  diff.   &  $\mathfrak{b}$     & $R^5$ & $B^{-i_{\beta}}\big((x-y+3z+9, ~x-y-3z+3,~3y-3)-24\varpi_{i_{\beta}}\big)$                        \\
                          &  diff.   &  $\mathfrak{a}$     & $s$   & $B^{-i_{\beta}}\big((-x-y-3z+9, ~2x-y+3,~3y-3)-24\varpi_{i_{\beta}}\big)$                        \\
                          
\hline
$(0,\beta,2)$  & same    &  $\mathfrak{c}$ & $R$    & $B^{-i_{\beta}}\big((-2x-y+9,~x-y+3z+3,~3y-3)-24\varpi_{i_{\beta}}\big)$                       \\
               &         &                 & $R^2s$ & $B^{-i_{\beta}}\big((2x-y+9,~ -x-y+3z+3,~3y-3)-24\varpi_{i_{\beta}}\big)$                    \\
               &  diff.  &  $\mathfrak{c}$ & $R^5s$ & $B^{-i_{\beta}}\big((-2x-y+9,~ x-y-3z+3,~3y-3)-24\varpi_{i_{\beta}}\big)$                     \\
               &         &                 & $R^4$  & $B^{-i_{\beta}}\big((2x-y+9,~-x-y-3z + 3,~3y-3)-24\varpi_{i_{\beta}}\big)$                     \\
               
\hline

$(3,\beta,1)$  & same   &  $\mathfrak{c}$ & $1$  &  $B^{-i_{\beta}}\big((-x-y+3z+9, ~2x-y+3,~3y-3)-24\varpi_{i_{\beta}}\big)$                         \\
               &        &                 & $Rs$ &  $B^{-i_{\beta}}\big((x-y-3z+9, ~x-y+3z+3,~3y-3)-24\varpi_{i_{\beta}}\big)$                          \\
               
\hline
$(9,\beta,3)$  & same  &  $\mathfrak{c}$  & $R^3$  & $B^{-i_{\beta}}\big((x-y-3z+9, ~-2x-y+3,~ 3y-3)-24\varpi_{i_{\beta}}\big)$                    \\
               &       &                  & $R^4s$ & $B^{-i_{\beta}}\big((-x-y+3z+9, ~-x-y-3z+3,~3y-3)-24\varpi_{i_{\beta}}\big)$                    \\
               
\hline
$(3,\beta,3)$  & diff.   & $\mathfrak{c}$  & $R^5$ & $B^{-i_{\beta}}\big((x-y+3z+9, ~x-y-3z+3,~3y-3)-24\varpi_{i_{\beta}}\big)$                      \\
               &         &                 &  $s$  & $B^{-i_{\beta}}\big((-x-y-3z+9, ~2x-y+3,~3y-3)-24\varpi_{i_{\beta}}\big)$                    \\
\hline
$(9,\beta,1)$  & diff.   & $\mathfrak{c}$ & $R^2$    & $B^{-i_{\beta}}\big((-x-y-3z+9,~-x-y+3z+3,~3y-3)-24\varpi_{i_{\beta}}\big)$                   \\
               &         &                & $R^3s$   & $B^{-i_{\beta}}\big((x-y+3z+9,~-2x-y+3,~3y-3)-24\varpi_{i_{\beta}}\big)$          \\   
\hline
               &         &                &          & 
\end{longtable}
\addtocounter{table}{-1}
}
}
\caption{Attaching elements of $G$ and affine Grassmannian elements to a given solution $X = (x,y,z)$.}

\end{table}

\newpage

We want to point out that the neat decomposition of $(x,y,z) \mod H$ (\Cref{prop neat decomposition}) is essential in the proof of the following lemma. It used to show that the $r_i$'s belong to $\mathbb{Z}$.

\bigskip

\begin{Lem} \label{lemma fonda table}
    Let $X=(x,y,z)\in\cU(48N+30)$ and $(\al,\beta,\gamma) = X \mod H$.
    \begin{enumerate}
        \item If $g \notin G_{(\alpha,\beta,\gamma)}$ then $g\varphi(r) \neq X$ for any $r \in L$.
        \item Let $g_{\raisebox{-0.5ex}{\scalebox{0.55}{$X$}}} \in G_{(\alpha,\beta,\gamma)}$. Then $r(X,g_{\raisebox{-0.5ex}{\scalebox{0.55}{$X$}}}) \in {{M}}$ if and only if the residue of $n$ modulo $3$ and the parity of $k,m$ are the same as those given in the row of $g_{\raisebox{-0.5ex}{\scalebox{0.55}{$X$}}}$ in \Cref{Table_parametrisation}.
    \end{enumerate}
\end{Lem}

\begin{proof}
\begin{enumerate}
    \item 
    The proof is a long list of computations and we only treat one case: $(\al, \beta, \ga) = (4,1,2)$, that is the first row of \Cref{Table_parametrisation} with $\beta = 1$, and leave the rest to the reader where we first need to change $\beta \in \{3,5,7\}$ and then go to the next row.
    Let $r \in L$ and assume that $g\varphi(r) = X$. Easy computations show that for the eight elements $ g\in G \setminus G_{(4,1,2)} = \{1, R^2,R^3,R^5,s,Rs,R^3s, R^4s\}$, the last coordinate of $g\varphi(r)$ modulo 4 is always 1, preventing $g\varphi(r) \mod H \neq (4,1,2)$, a contradiction.

    \item Let us prove the direct implication first. The proof is a long verification. 
    We treat here only one case in detail (out of $144$), the remaining cases being analogous. Our case is $(\al,\beta,\ga)= (4,1,2)$ and $g_{\raisebox{-0.5ex}{\scalebox{0.55}{$X$}}} = R$ (first row of \Cref{Table_parametrisation} for $\beta = 1$, which gives $i_{\beta}=0$ and $\varpi_{i_{\beta}}=0$). Therefore one has
    \begin{align*}
        24r_0(X,R) & = r_1\varepsilon_1 + r_2\varepsilon_2 + r_3\varepsilon_3 \\
                   & =(-2x-y+9)\varepsilon_1 + (x-y+3z+3)\varepsilon_2 + (3y-3)\varepsilon_3.
    \end{align*}

    By the first point above we know that $G_{(4,1,2)} = \{R,R^2s,R^5s, R^4\}$. We have then three situations to consider and the claim is that for the second, third and forth case, taking the corresponding data of the row (that is the residue of $n$ modulo 3 and the parity of $k,m$), one has $r_0(X,R) \notin M$. To do so we must show that at least one of the $r_i$'s does not belong to $24\mathbb{Z}$. As $\beta = 1$, one has by \Cref{abc} that $\mathfrak{a} = 0$ and $\mathfrak{b}=2$. Therefore
    \begin{itemize}
        \item if ~$n \mod 3 = 2$ and $k,m$ have different parity then the first coordinate of $24r_0(X,R)$ yields the problem. Write $k = 2a+\varepsilon$ with $\varepsilon \in \{0,1\}$ and $n = 3b+2$. Then $-2x-y+3 =  24(-2a-b-\varepsilon) + 8$, which does not belong to $24\mathbb{Z}$.
        
        \item If ~$n \mod 3 = 0$ and $k,m$ have same parity then the second coordinate of $24r_0(X,R)$ yields the problem. Write $k = 2a+\epsilon$ and $m = 2c + \epsilon$ with $\epsilon =0,1$. Then $x-y+3z+3 = 24(-a-n+c+\varepsilon) + 12$, which does not belong to $24\mathbb{Z}$.
        
        \item if ~$n \mod 3 = 2$ and $k,m$ have same parity then again the first coordinate of $24r_0(X,R)$ yields the problem. Write $k = 2a$, $n = 3b+2$ and $m = 2c$ (the case $k$ and $m$ odd is similar). Then $-2x-y+9 = 24(-2a-b-1) + 8$, which does not belong to $24\mathbb{Z}$.

    \end{itemize}
    
    Let us now prove the converse. Therefore, we follow the setting of each row of \Cref{Table_parametrisation} and we show under these conditions that
        $r(X,g_{\raisebox{-0.5ex}{\scalebox{0.55}{$X$}}}) \in {{M}}$. 
        Again, we only treat one case and leave the others (that are similar) to the reader. 
        Our case is $(\al,\beta,\ga) = (4,7,2)$, with $n = 0 \mod 3$, $k$ and $m$ having same parity and $g_{\raisebox{-0.5ex}{\scalebox{0.55}{$X$}}} = R^5s$ (that is the second row of \Cref{Table_parametrisation} for $\beta = 7$, giving then $i_{\beta} = 1$). Therefore we have 
        \begin{align*}
            24r(X,R^5s)     & = r_1\varepsilon_1 + r_2\varepsilon_2 + r_3\varepsilon_3 \\
                            & =B^{-1}\big((- 2x - y +9,~ x - y - 3z + 3,~3y-3)-24\varpi_{1}\big) \\
                            & = (x-y-3z+9)\varepsilon_1 + (3y+3)\varepsilon_2 + (x-y+3z-3)\varepsilon_3.
        \end{align*}

        We show now that $r(X,R^5s) \in {{M}}$, i.e., the $r_i$'s all belong to $24\mathbb{Z}$. Our context is the following: $x = 12k + 4$, $y = 8n + 7$, $z = 4m+2$ with $k$ and $m$ having same parity and $n = 3b$ for some $b \in \mathbb{Z}$, implying $y = 24b+7$. Write $k = 2a+\epsilon$ and $m = 2c + \epsilon$ with $\epsilon =0,1$.
        \begin{itemize}
            \item $r_1 \in 24\mathbb{Z}$: one has $x-y-3z+9 = 24(a-b-c)$. 
            
            \item $r_2 \in 24\mathbb{Z}$: one has $3y+3 = 24(3b+1)$.

            \item $r_3 \in 24\mathbb{Z}$: one has $x-y+3z-3 = 24(a-b-c+\varepsilon)$.
        \end{itemize}
\end{enumerate}
\end{proof}

We are finally ready to prove our main theorem.
\medskip

\textit{Proof of \Cref{pig A3}.}
Write $(\al,\beta,\ga) = X \mod H$.
Using the appropriate row of \Cref{Table_parametrisation} associated to $X$ that provides the element(s) $g_{\raisebox{-0.5ex}{\scalebox{0.55}{$X$}}} \in G_{(\al,\beta,\ga) }$, we know by \Cref{lemma fonda table} that $q:=r(X,g_{\raisebox{-0.5ex}{\scalebox{0.55}{$X$}}}) \in {{M}}$. Moreover, by \Cref{rem_thm_pig} (1) we have $g_{\raisebox{-0.5ex}{\scalebox{0.55}{$X$}}}\varphi(\sigma_{i_\beta} t_q)=X$, that is 
 $g_{\raisebox{-0.5ex}{\scalebox{0.55}{$X$}}}\varphi\big(\varpi_{i_{\beta}}+ w_{0,i_{\beta}}w_0(q)\big) = X$.  But then, setting $r = \varpi_{i_{\beta}}+  w_{0,i_{\beta}}w_0(q)$, since $q \in {{M}}$, it follows by \Cref{L in terms of M} (2) that $r \in L$. Hence $g_{\raisebox{-0.5ex}{\scalebox{0.55}{$X$}}}\varphi(r) = X$.
\hfill $\square$

    \begin{figure}[h!]
    \begin{center}
    \includegraphics[scale=0.6]{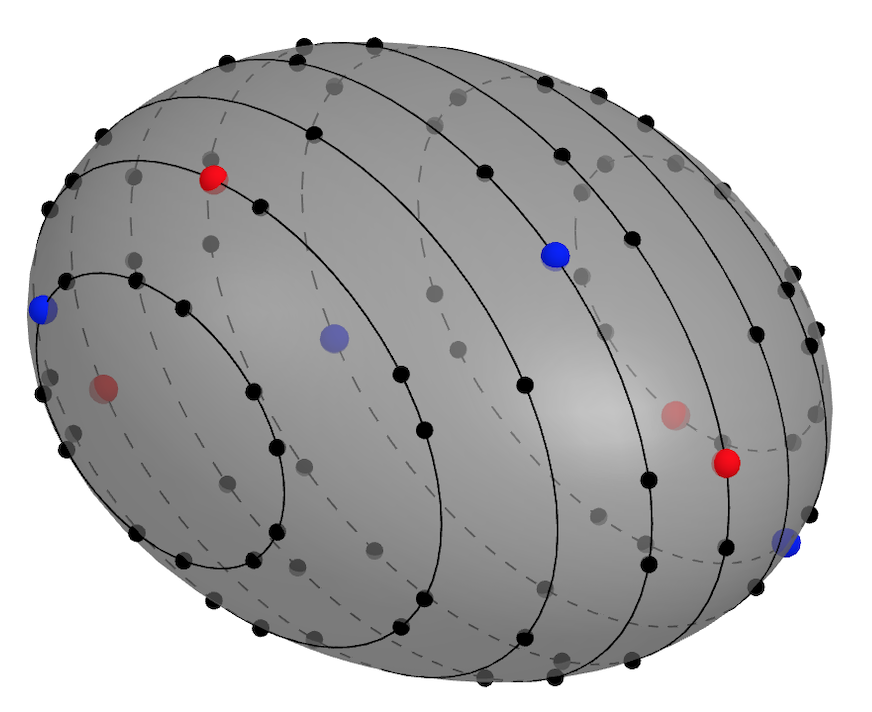}
    \caption{Integral solutions of the equation $x^2+2y^2 + 3z^2 = 30+48\cdot 2$, that is when the $\Lambda_0$-atomic length is equal to 2.  There are two affine Grassmannian elements of $\Lambda_0$-atomic length 2: $q:=(1, 0, -1, 0)$ and $q':=(0, 1, 0, -1)$ in the standard basis.  Each of them generates 3 other elements of same $\Lambda_0$-atomic length, namely $\sigma_it_{q}$ and $\sigma_i t_{q'}$ for $i =1,2,3$.  The blue points are the solutions $\varphi(q)$ and $\varphi(\varpi_i + w_{0,i}w_0(q))$ while the red points are the solutions $q'$ and $\varphi(\varpi_i + w_{0,i}w_0(q))$.}\label{image 3d N=2}
    \end{center}
    \end{figure}

\subsection{$G$-set structure on solutions}
\label{sec_variety}

In this section, we describe the $G$-orbits of solutions.
We give in \Cref{thm size orbits} a characterisation of the orbits of a given size (which is either $6$ or $12$)
by certain arithmetic conditions.
\Cref{intersection vide orbites A3} is an analogue of \Cref{same orbit A2}

\begin{Prop}\label{prop orbites}
Let $N \in \mathbb{N}$.
\begin{enumerate}
\item For all $X\in \cU(48N+30)$, we have 
$|\cO_{\langle R\rangle}(X)| = 6$.
\item For all $X\in \cU(48N+30)$, we have 
$|\cO_G(X)|\in\{6,12\}$. 
\item For all $X\in \cU(48N+30)$,
the $G$-orbit of $X$ has size $6$ if and only if
$\mathrm{Stab}_G(X)=\langle sR^k\rangle$ for some odd $k$.
\end{enumerate}
\end{Prop} 

\textit{Proof.}\
\begin{enumerate}
\item 
The proof is exactly the same as that of \cite[Theorem 2.1]{BrunatNath2022}.
\item
It is enough to consider the action of $\langle R \rangle$ and to show that for any $X \in \mathcal{U}(48N+30)$, the stabiliser of $X$ is trivial. Therefore, by (1) as $\langle R \rangle$ is a subgroup (of order 6) of $G$ (of order 12), we deduce the announced statement. 
\item
We have seen in the previous point that the group $\langle R \rangle $ acts regularly on $\cU(48N+30)$. 
Therefore, for $G = \langle R,s \rangle$, the only situation that prevents the $G$-orbit of $X$ to have 
size 12 is if $sX$ belongs to the $\langle R \rangle$-orbit of $X$, that is if there exists $k \in \llbracket 1,6\rrbracket $ such that $sX = R^kX$.
Therefore, it remains to discard the situations $sX = R^kX$ for $k = 0,2,4$. If $sX = X$ then it follows that $z = 0$, but by \Cref{prop neat decomposition} we know that $\gamma$ is either $1,2,3,$ implying that $z = 4m+\ga$ can never be equal to $0$. If now $sX = R^2X$ then one has $(x,y,-z)=1/2(-x-3z,2y,x-z)$, that is $x = 0$ and $z = 0$, which for the same reason as above is impossible (notice that an other interesting argument is the following: $2y^2 = 48N+30$ implies that $y^2 = 28N+15$ and then $y^2 = 7 \mod 8$, however the residues of a square modulo $8$ are $0,1,4$). If $sX = R^4X$ then $(x,y,-z)=1/2(-x+3z,2y,-x-z)$, which is equivalent to $x = 0$ and $z = 0$ and then is not possible. 
\hfill
$\square$
\end{enumerate}

We now characterise which points $X$ yield an orbit of size 6 or 12.
Let $X=(x,y,z)\in\cU(48N+30)$. As before, let us write $x = 12k+\alpha$, $y = 8n+\beta$ and $z = 4m + \gamma$ where $\alpha \in \llbracket 0,11 \rrbracket$, $\beta \in \llbracket 0,7 \rrbracket$ and $\gamma \in \llbracket 0,3\rrbracket$. 

\bigskip

\begin{Th}\label{thm size orbits}
    The $G$-orbit of $X$ has size $6$ if and only one of the three following situations occurs
    \begin{enumerate}
        \item $k = \al = 0$ (in which case $X \mod H = (0,\beta,2)$).
        \item $k+m = -1$ and $\al + 3\ga = 12$ (in which case $X \mod H \in \Omega_{seg}(\beta)$).
        \item $k = m$ and $\al - 3\ga = 0$ (in which case $X \mod H \in \Omega_{tri}(\beta)$).
    \end{enumerate}
\end{Th}

\begin{proof}
    By \Cref{prop orbites} (3), we only need to look at $sX = R^kX$ for $k = 1,3,5$.
    \begin{itemize}
        \item If $sX = R^3X$ then $(x,y,-z)=(-x,y,-z)$, that is $x = 0$ and then $12k+\al = 0$. Hence $k = \al = 0$. By \Cref{prop neat decomposition} there is only one possible class modulo $H$ allowing $\al = 0$ and this is $(0,\beta,2)$, whence $(1)$.
        \item If $sX = RX$ then $(x,y,-z)=1/2(x-3z,2y,x+z)$, which is equivalent to $x = -3z$. Therefore $12k+\al = -3(4m+\ga)$, that is $12(k+m) + (\al + 3\ga) = 0$. But now, in order to have $12(k+m) + (\al + 3\ga) = 0$, the integer $\al + 3\ga$ must be a  multiple of 12 and an easy verification on $\Omega$, via \Cref{section_arbre}, shows there are only two situations allowing it: $(3,\beta,3)$ and $(9,\beta,1)$, that is the elements of $\Omega_{seg}(\beta)$. However, for these two classes one has $\al+3\ga = 12$, which shows that in order to have $12(k+m) + (\al + 3\ga) = 0$ we must have $k+m = -1$.
        \item If $sX = R^5X$ then $(x,y,-z)=1/2(x+3z,2y,-x+z)$, i.e., $x=3z$. Then $12k+\al = 3(4m+\ga)$, that is $12(k-m) + (\al - 3\ga) = 0$. As above, $\al - 3\ga$ must be a  multiple of 12 and there are only two situations allowing it: $(3,\beta,1)$ and $(9,\beta,3)$. For this two classes one has $\al-3\ga = 0$, implying then $k=m$.
    \end{itemize}
\end{proof}

\bigskip

\begin{center}
\begin{figure}[h!]
    \centering
  \includegraphics[scale=0.55]{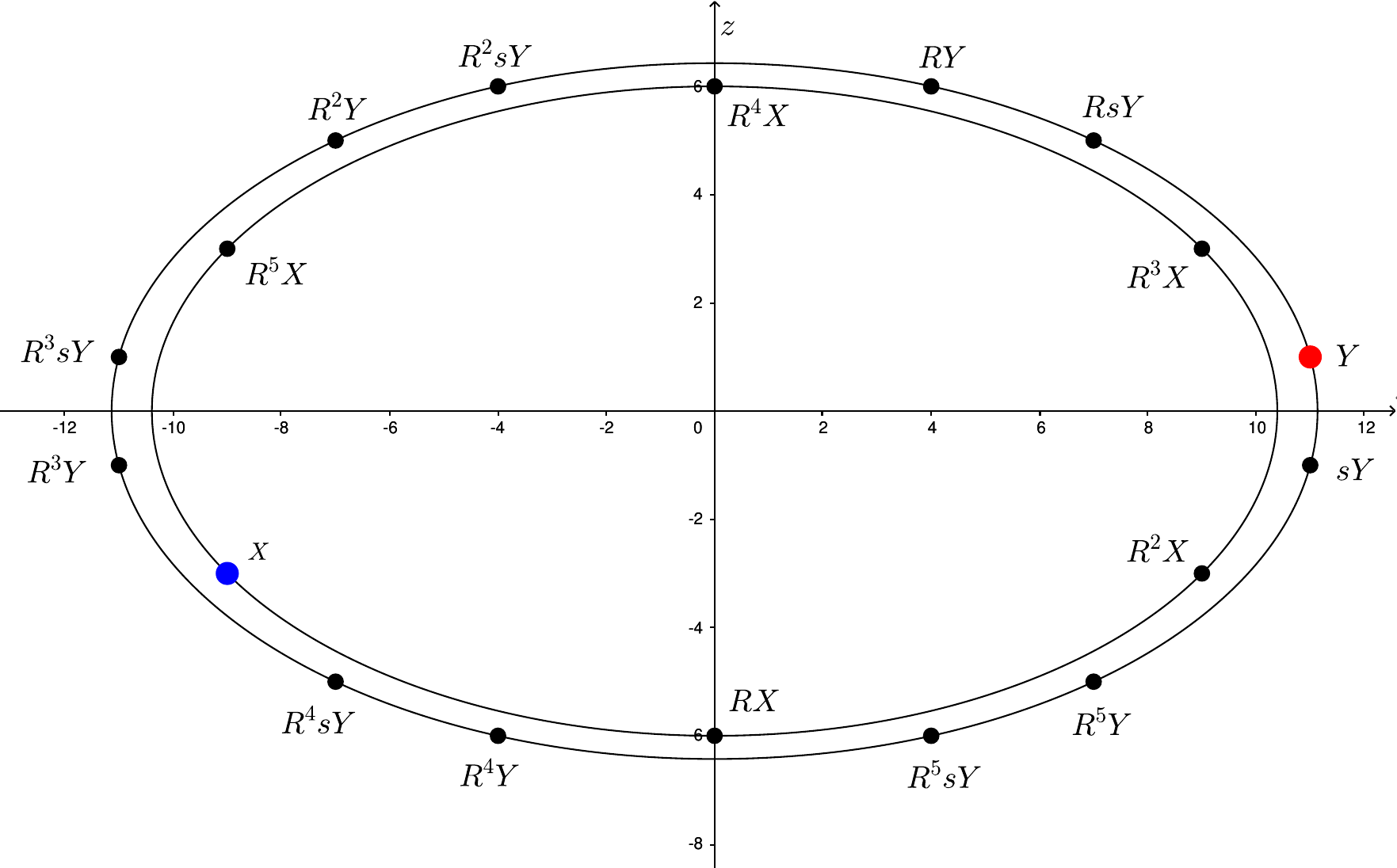}  
    \caption{Two slices of solutions for $N = 2$.}
    \label{Two slices of solutions}
\end{figure}

\end{center}

\begin{Exa}
We continue the example in \Cref{image 3d N=2} where we take two slices: $y = 1$ and $y = 3$.
We represent them in \Cref{Two slices of solutions} and we detect the $G$-orbit of $X = (-9,3,-3)$ and $Y = (11,1,1)$. The $G$-orbit of $X$ has $6$ elements since $-9 = 12\cdot(-1) + 3$ and $-3 = 4\cdot (-1) + 1$, that is $k = -1, m = -1$ and $\al = 3,\ga = 1$. Hence $k=m$ and $\al - 3\ga = 0$ (\Cref{thm size orbits} (3)).
On the contrary, the $G$-orbit of $Y$ has size $12$ because none of the conditions of \Cref{thm size orbits} are satisfied.
\end{Exa}

\bigskip

The following proposition explores the case where the image via $\varphi$ of multiple extended affine Grassmannian elements are contained within the same $G$-orbit. In other words, this approach provides insight about the structure of 
$G$-orbits when several \say{same-colored points} are found in the same orbit (see \Cref{image 3d N=2} and imagine two blue points in the same orbit). It will be particularly helpful to understand \Cref{Table_parametrisation}.

\begin{Th}\label{intersection vide orbites A3}
Let $r, r' \in L$ with $r \neq r'$.
Assume that $g\varphi(r) = \varphi(r')$ with $g \in G$. Then
\begin{enumerate}
    \item $g = Rs$,
    \item $|\mathcal{O}_G(\varphi(r))| =  12$,
    \item there is no other $r'' \in L$ such that $ \varphi(r'') \in \mathcal{O}_G(\varphi(r))$.
\end{enumerate}
\end{Th}

\begin{proof} 
By \Cref{dec L A3} write $r = \varpi_i + w_{0,i}w_0(q)$ and $r' = \varpi_j + w_{0,j}w_0(q')$ with $q,q' \in {{M}}$ and $i,j \in \{0,1,2,3\}$.
\begin{enumerate}
    \item 
By \Cref{rem_thm_pig} (3), since $g\varphi(r) = \varphi(r')$, 
it follows that $i=j$, from which we obtain the following equality
\begin{equation}\label{equ proof}
    g\varphi(\varpi_i + w_{0,i}w_0(q)) = \varphi(\varpi_i + w_{0,i}w_0(q')).
\end{equation}
Write $q=q_1\eps_1+q_2\eps_2 +q_3\eps_3$ and 
$q'=q_1'\eps_1 + q_2'\eps_2 +q_3'\eps_3$.  
We are going to show that \Cref{equ proof} is impossible unless $g = sR^5 = Rs$, see \Cref{relation G rank 3}.
\medskip

Recall that in matrix form, 
$\varphi(\varpi_i + w_{0,i}w_0(q')) = Q_i q' + u_i$ with the notations of \Cref{sec_A3_phi}. 
The first coordinate of 
$Q_i q' + u_i$, for $i=0,1,2,3$, is
$12q_2' + 4q_3' - 1$, ~
$12q_1' + 4q_2' - 5$,~
$-8q_1' - 12q_2' - 12q_3' + 3$
and
$-4q_1' - 4q_2' + 8q_3' + 3$,
respectively.
In particular their residue modulo $4$ is always $3$. 
Moreover, for $m=0,1$, the first coordinate of $s^m R^k\varphi(\varpi_i + w_{0,i}w_0(q))$ (for $k=0,1,2,3,4,5$) is respectively:

\begin{itemize}
    \item If $i = 0$: 
    $12q_2+4q_3-1,~ 
    -12q_1 - 4q_3 + 4,~
    -12q_1 - 12q_2 - 8q_3 + 5,~
    -12q_2 - 4q_3 + 1,~
    12q_1 + 4q_3 - 4,~
    12q_1 + 12q_2 + 8q_3 - 5
    $. 
    \item[]
    \item If $i = 1$: 
    $
    12q_1 + 4q_2 - 5, 8q_2 - 1,~
    12q_1 + 8q_2 + 12q_3 - 4,~
    4q_2 + 12q_3 + 1,~
    -12q_1 - 4q_2 + 5,~
    -12q_1 - 8q_2 - 12q_3 + 4,~
    -4q_2 - 12q_3 - 1
    $.   
    \item[] 
    \item If $i = 2$: 
    $
    -8q_1 - 12q_2 - 12q_3 + 3,~
 -4q_1 - 12q_3,~
 4q_1 + 12q_2 - 3,~
 8q_1 + 12q_2 + 12q_3 - 3,~
 4q_1 + 12q_3,~
 -4q_1 - 12q_2 + 3
    $.  
    \item[]
    \item If $i = 3$: 
    $-4q_1 - 4q_2 + 8q_3 + 3,~
 4q_1 - 8q_2 + 4q_3,~
 8q_1 - 4q_2 - 4q_3 - 3,~
 4q_1 + 4q_2 - 8q_3 - 3,~
 -4q_1 + 8q_2 - 4q_3,~
 -8q_1 + 4q_2 + 4q_3 + 3
    $.
\end{itemize}

We see that for each $i$, the only case having a residue modulo 4 equal to 3 is always the last one, that is for $k = 5$. We thus shrunk our candidates to only two: $ g= R^5$ or $g = sR^5$. We discard now $R^5$.

On one side, the third coordinate of $\varphi(\varpi_i + w_{0,i}w_0(q'))$ (for $i=0,1,2,3$) is respectively 
$8q_1 + 4q_2 + 4q_3 - 3,
-4q_1 - 4q_2 - 8q_3 + 1,
-4q_2 + 4q_3 + 1,
-4q_1 + 4q_2 + 1$. 
Therefore, their residue modulo $4$ is $1$ for each one of them. On the other side, the third coordinate of $R^5\varphi(\varpi_i +  w_{0,i}w_0(q))$ (for $i=0,1,2,3$) 
is respectively 
$4q_1 - 4q_2 - 1,
-8q_1 - 4q_2 - 4q_3 + 3,
4q_1 + 4q_2 + 8q_3 - 1,
4q_2 - 4q_3 - 1$. 
Their residues modulo 4 are all equal to 3. This prevents the equality $R^5\varphi(\varpi_i + w_{0,i}w_0(q)) = \varphi(\varpi_i + w_{0,i}w_0(q'))$.

\item
By \Cref{prop orbites} (1), we know that the $\langle R \rangle$-orbit of $\varphi(r)$ is of size 6. Moreover, by what we have done in the proof of (1) for $m = 0$, we know that $R^k\varphi(r) \neq \varphi(r')$ for any $k=0,1,2,3,4,5$ ($0$ coming from our assumption $r \neq r'$). However, since $\varphi(r)$ and $\varphi(r')$ belong to the same $G$-orbit by assumption, it follows that there are at least 7 elements in this orbit and by \Cref{prop orbites} (2) the result follows.

\item 
Assume that we have a third element $r''$ belonging to the $G$-orbit of $r$ and $r'$. Therefore, by (1) it follows that $\varphi(r'') = Rs \varphi(r') =  RsRs\varphi(r) \overset{\Cref{relation G rank 3}}{=} \varphi(r)$. This completes the proof since $\varphi$ is bijective.
\end{enumerate}
\end{proof}


\section{Underlying type $B_n/C_n$ and action of the hyperoctahedral group}
\label{sec_Pell_hyperoct}

In this section we generalise \Cref{prop_A42} to Dynkin types
that have underlying finite type $B_n$ or $C_n$. According to \Cref{table},
these are type
$B_n^{(1)}$, $C_n^{(1)}$,  $A_{2n-1}^{(2)}$, $D_{n+1}^{(2)}$, $A_{2n}^{(2)}$, by considering the action of the hyperoctohedral group $H_n$ (which is nothing but the Weyl group of type $BC_n$).
For $a,b\in \Z$, denote again $\cU(aN+b)$ the set of integer solutions of the Diophantine equation
\begin{equation}
\label{eq_hyperoct_ab}
x_1^2+\ldots + x_n^2 = a N + b.
\end{equation}

Using \Cref{thm_LA}, it is easy to establish the general formula for $\sL_{\La_0}(t_q)$, and
to compute the corresponding Gaussian reduction.
We obtain the following results.

\subsection{Type $B_n^{(1)}$}

For all $q\in {{M}}$, we have
\begin{align*}
\sL_{\La_0}(t_q) 
& =
n\sum_{i=1}^nq_i^2 - \sum_{i=1}^n(n-i+1)q_i
\\
& =
 n \sum_{i=1}^n \left( q_i - \frac{n-i+1}{2n}  \right)^2 -  \frac{n(n+1)(2n+1)/6}{4n}
\end{align*}
where the numerator of the last fraction is obtained by the formula for the sum of the first $n$ squares.
This gives the following equivalent identity.

\begin{Lem}
\label{lem_Bn1}
$$4n \sL_{\La_0}(t_q) + \frac{n(n+1)(2n+1)}{6} = \sum_{i=1}^n \left( 2n q_i - (n-i+1) \right)^2.$$
\end{Lem}

Therefore, if we take $a=4n$ and $b=\frac{n(n+1)(2n+1)}{6}$,
we obtain a map
\begin{equation}
\label{phi_Bn1}
\begin{array}{cccc}
\varphi: & \cB(N)  & \longrightarrow & \cU(aN+b)
\\
& (q_1,\ldots,q_i, \ldots, q_n) & \longmapsto & (\ldots, 2n q_i - (n-i+1), \ldots).
\end{array}
\end{equation}
In other terms, each $q\in {{M}}$ with produces a solution of \Cref{eq_hyperoct_ab}
with $N=\sL_{\La_0}(t_q)$.

\subsection{Type $C_n^{(1)}$}

For all $q\in {{M}}$, we have
\begin{align*}
\sL_{\La_0}(t_q) 
& =
2n\sum_{i=1}^nq_i^2 - \sum_{i=1}^n(2(n-i)+1)q_i
\\
& =
 \sum_{i=1}^n 2n\left( q_i - \frac{2(n-i)+1}{4n}  \right)^2 -  \frac{n(2n+1)(2n-1)/3}{8n}
\end{align*}
where the numerator of the last fraction is obtained by the formula for the sum of the first $n$ odd squares.
This gives the following equivalent identity.

\begin{Lem}
\label{lem_Cn1}
$$8n \sL_{\La_0}(t_q) + \frac{n(2n+1)(2n-1)}{3} = \sum_{i=1}^n \left( 4n q_i - (2(n-i)+1) \right)^2.$$
\end{Lem}

Therefore, if we take $a=8n$ and $b=\frac{n(2n+1)(2n-1)}{3}$,
we obtain a map
\begin{equation}
\label{phi_Cn1}
\begin{array}{cccc}
\varphi: & \cB(N)  & \longrightarrow & \cU(aN+b)
\\
& (q_1,\ldots,q_i, \ldots, q_n) & \longmapsto & (\ldots, 4n q_i - (2(n-i)+1), \ldots).
\end{array}
\end{equation}
In other terms, 
by the results of \Cref{sec_cores_C},
each self-conjugate $2n$-core of size $N$ produces a solution of \Cref{eq_hyperoct_ab}.

\subsection{Type $A_{2n-1}^{(2)}$}

For all $q\in {{M}}$, we have
\begin{align*}
\sL_{\La_0}(t_q) 
& =
\sum_{i=1}^nq_i^2 - \sum_{i=1}^n\frac{(2(n-i)+1)}{2}q_i
\\
& =
\frac{2n-1}{2} \sum_{i=1}^n \left( q_i - \frac{2(n-i)+1}{2(2n-1)}  \right)^2 -  \frac{2n-1}{2}\frac{n(2n+1)(2n-1)/3}{4(2n-1)^2}
\end{align*}
This gives the following equivalent identity.

\begin{Lem}
\label{lem_A2n-12}
$$(16n-8)\sL_{\La_0}(t_q) + \frac{n(2n+1)(2n-1)}{3} = \sum_{i=1}^n \left( (4n-2)q_i - (2(n-i)+1) \right)^2.$$
\end{Lem}

Therefore, if we take $a=(16n-8)$ and $b=\frac{n(2n+1)(2n-1)}{3}$,
we obtain a map
\begin{equation}
\label{phi_A2n-12}
\begin{array}{cccc}
\varphi: & \cB(N)  & \longrightarrow & \cU(aN+b)
\\
& (q_1,\ldots,q_i, \ldots, q_n) & \longmapsto & (\ldots, (4n-2) q_i - (2(n-i)+1), \ldots).
\end{array}
\end{equation}
In other terms, each $q\in {{M}}$ with produces a solution of \Cref{eq_hyperoct_ab}
with $N=\sL_{\La_0}(t_q)$.

\subsection{Type $D_{n+1}^{(2)}$}

For all $q\in {{M}}$, we have
\begin{align*}
\sL_{\La_0}(t_q) 
& =
(n+1)\sum_{i=1}^nq_i^2 - \sum_{i=1}^n(n-i+1)q_i
\\
& =
 \sum_{i=1}^n (n+1)\left( q_i - \frac{n-i+1}{2(n+1)}  \right)^2 -  \frac{n(n+1)(2n+1)/6}{4(n+1)}.
\end{align*}.
This gives the following equivalent identity.

\begin{Lem}
\label{lem_Dn+12}
$$4(n+1) \sL_{\La_0}(t_q) + \frac{n(n+1)(2n+1)}{6} = \sum_{i=1}^n \left( 2(n+1) q_i - (n-i+1) \right)^2.$$
\end{Lem}

Therefore, if we take $a=4(n+1)$ and $b=\frac{n(n+1)(2n+1)}{6}$,
we obtain a map
\begin{equation}
\label{phi_Dn+12}
\begin{array}{cccc}
\varphi: & \cB(N)  & \longrightarrow & \cU(aN+b)
\\
& (q_1,\ldots,q_i, \ldots, q_n) & \longmapsto & (\ldots, 2(n+1) q_i - (n-i+1), \ldots).
\end{array}
\end{equation}
In other terms,  by the results of \Cref{sec_cores_Dt},
each bar core partition in $\cD_{2n+2}(N)$ produces a solution of \Cref{eq_hyperoct_ab}.

\subsection{Type $A_{2n}^{(2)}$}

For all $q\in {{M}}$, we have
\begin{align*}
\sL_{\La_0}(t_q) 
& =
\frac{2n+1}{2}\sum_{i=1}^nq_i^2 - \sum_{i=1}^n(n-i+\frac{1}{2})q_i
\\
& =
 \sum_{i=1}^n \frac{2n+1}{2}\left( q_i - \frac{2(n-i)+1 }{4n+2}  \right)^2 -  \frac{(2n+1)}{2}.\frac{n(2n+1)(2n-1)/3}{(4n+2)^2}
\\
& =
 \sum_{i=1}^n \frac{2n+1}{2}    \left( q_i - \frac{2(n-i)+1 }{4n+2}  \right)^2 -  \frac{n(2n+1)(2n-1)/3}{8(2n+1)}.
\end{align*}
By multiplying both sides by $8(2n+1)$, This gives the following equivalent identity.

\begin{Lem}
\label{lem_A2n2}
$$(16n+8) \sL_{\La_0}(t_q) + \frac{n(2n+1)(2n-1)}{3} = \sum_{i=1}^n \left( (4n+2) q_i - (2(n-i)+1) \right)^2.$$
\end{Lem}

Therefore, if we take $a=(16n+8)$ and $b=\frac{n(2n+1)(2n-1)}{3}$,
we obtain a map
\begin{equation}
\label{phi_A2n2}
\begin{array}{cccc}
\varphi: & \cB(N)  & \longrightarrow & \cU(aN+b)
\\
& (q_1,\ldots,q_i, \ldots, q_n) & \longmapsto & (\ldots, (4n+2) q_i - (2(n-i)+1) , \ldots).
\end{array}
\end{equation}
In other terms, each $q\in {{M}}$ with produces a solution of \Cref{eq_hyperoct_ab}
with $N=\sL_{\La_0}(t_q)$.

\subsection{Action of the hyperoctahedral group}

Now, let $H_n=C_2^n\rtimes S_n$ be the hyperoctahedral group in dimension $n$, of order $2^nn!$.
Similarly to \Cref{sec_pell},
$H_n$ acts on the solution set of 
\begin{equation}
\label{eq_hyperoct}
x_1^2+\cdots +x_n^2 = k,
\end{equation}
where the generator of $C_2$ acts by $x_1\mapsto -x_1$, and $S_n$ permutes the $x_i$'s.
Denote as usual $H_n(x)$ the $H_n$-orbit of a solution $x=(x_1,\ldots, x_n)$ of \Cref{eq_hyperoct}. 
We are able to generalise \Cref{prop_A42}
to all of the above types.

\begin{Prop}
\label{prop_hyperoct}
Let $W$ be the one of the above Weyl groups,
and consider the corresponding map $\varphi$.
For all $q\in {{M}}$, $$\left|H_n\varphi(q)\right| = 2^nn!.$$
\end{Prop}

In other words, each $q\in {{M}}$ with $\sL_{\La_0}(t_q)=N$ produces $2^nn!$ solutions of \Cref{eq_hyperoct_ab}.

\begin{proof}
The proof is completely analogous to that of \Cref{prop_A42}.
Let $x_1^2+\cdots +x_n^2 = aN+b$ be the associated Diophantine equation.
Let $q\in {{M}}$ and write $x=\varphi(q)$. 
The $H_n$-orbit of $x$ is $$\{ (\pm x_{\sigma(1)}, \ldots, \pm x_{\sigma(n)}) \,;\, \sigma\in S_n\}.$$
Let $c$ be the coefficient of each $q_i$ in $\varphi(q)$.
Now we have in all cases
\begin{itemize}
\item $x_i\neq \pm x_j\mod c$ and hence $x_i\neq \pm x_j$ for all $1\leq i\neq j\leq n$,
\item $x_i\neq 0$ for all $1\leq i\leq n$.
\end{itemize}
By the direct analogue of \Cref{lem_free_1} for $H_n$, the element in the orbit above are all distinct.
\end{proof}

\begin{Rem} Let us focus on type $C_n^{(1)}$, where self-conjugate $2n$-cores appear naturally
because $\{q\in {{M}}\mid \sL_{\La_0}(t_q)=N\} \overset{\sim}{\longleftrightarrow}\SC_{2n}(N)$ (see \Cref{sec_comb_models}).
\begin{enumerate}
\item 
The case $n=3$ yields, by \Cref{lem_Cn1}
$$(12q_1-5)^2+(12q_2-3)^2+(12q_3-1)^2 = 24\sL_{\La_0}(t_q) +35.$$
By \Cref{prop_hyperoct} each self-conjugate $6$-core of size $N$ produces $48$ distinct solutions of the equation
$$x^2+y^2+z^2=24N+35.$$
In \cite[Theorem 2]{HansonJameson2022}, it is proved that
$$\left|\SC_6(N)\right| = \frac{1}{12}\left|\left\{ (x,y,z)\in\Z^3 \mid 3x^2+32y^2+32yz+32z^2 = 24N+35  \right\}\right|,$$
correcting an error in \cite[Theorem 6]{Alpoge2014}.
It is interesting to notice that the quantity $24N+35$ appears naturally in both contexts.
\item 
The case of $n=4$ yields, by \Cref{lem_Cn1}
$$32 \sL(w) + 84 = (16q_1-7)^2+(16q_2-5)^2+(16q_3-3)^2+(16q_4-1)^2.$$
By \Cref{prop_hyperoct} each self-conjugate $8$-core of size $N$ produces 
$384$ distinct solutions of the equation
$x^2+y^2+z^2+t^2 = 32N+21$.
This can also be compared to \cite[Theorem 4]{Alpoge2014}, where 
$\left|\SC_8(N)\right|$ is claimed to be in relationship with the number of integer solutions of $x^2+y^2+2z^2+2t^2 = 8N+21$, and we remark that
$8N+21 = \frac{1}{4}(32N+84).$
However, Alpoge's equation is not quite correct, and it would be interesting to tweak his statement using an appropriate quadratic expression.
\end{enumerate}
\end{Rem}

\section*{Acknowledgements}

We thank Cédric Lecouvey and David Wahiche for sharing and discussing the results of \cite{LW2024}. We also thank James Parkinson and Arthur Garnier for useful communications, and Michael Hanson for sharing and discussing some notes. The authors were supported by the Agence Nationale de la Recherche funding ANR CORTIPOM 21-CE40-001.

\bibliographystyle{alphaurl}
\bibliography{references}

\end{document}

%% file: table.tex
\afterpage{%
  \clearpage
\begin{landscape}
\begin{figure}
\rowcolors{2}{gray!25}{white}
\footnotesize
$$
\begin{array}{@{}l@{\hskip 20pt} @{}l@{\hskip 20pt} @{}l@{\hskip 20pt} @{}l@{\hskip 20pt} @{}l@{\hskip 20pt} @{}l@{\hskip 20pt}  @{}l@{\hskip 20pt} @{}l@{}}
\hline 
\text{Type}
&
a_0, a_1,\ldots, a_n
&
a_0^\vee, a_1^\vee ,\ldots, a_n^\vee
&
W_0
&
\text{Simple roots realisation}
&
\text{Lattice $M^{\dagger}$ and its realisation}
&
W
&
h
\\
\hline
A_n^{(1)} \ n\geq 1
&
1, 1,\ldots, 1
&
1, 1,\ldots, 1
&
A_n
&
\begin{array}[t]{@{}ll@{}}
\textcolor{blue!80}{\al_i=e_i-e_{i+1} ,\, 1\leq i\leq n}
\end{array}
&
\begin{array}[t]{@{}ll@{}}
\Z\al_1+\cdots+\Z\al_{n-1}+\Z\al_n =
\\
\left\{ 
\be \in \Z e_1+\cdots+\Z e_{n+1}
\mid 
\sum_{1\leq i\leq n+1} \be_i =0 \right\} 
\end{array}
&
\widetilde{A}_n
&
n+1
\\
B_n^{(1)}  \ n\geq 3
&
1, 1,2,\ldots, 2,2
&
1, 1,2,\ldots, 2,1
&
B_n
&
\begin{array}[t]{@{}ll@{}}
\textcolor{blue!80}{\al_i=e_i-e_{i+1} ,\, 1\leq i\leq n-1}
\\
\al_n=e_n
\end{array}
&
\begin{array}[t]{@{}ll@{}}
\Z\al_1+\cdots+\Z\al_{n-1}+2\Z\al_n =
\\
\left\{ 
\be \in \Z e_1+\cdots+\Z e_n
\mid 
\sum_{1\leq i\leq n} \be_i \text{ even } \right\} 
\end{array}
&
\widetilde{B}_n
&
2n
\\
C_{n}^{(1)}  \ n\geq 2
&
1, 2,\ldots, 2,1
&
1, 1,\ldots, 1,1
&
C_n
&
\begin{array}[t]{@{}ll@{}}
\al_i=\frac{1}{\sqrt{2}}\left(e_i-e_{i+1}\right) ,\, 1\leq i\leq n-1
\\
\textcolor{blue!80}{\al_n=\sqrt{2}e_n}
\end{array}
&
\begin{array}[t]{@{}ll@{}}
2\Z\al_1+\cdots+2\Z\al_{n-1}+\Z\al_n =
\sqrt{2}\Z e_1+\cdots+\sqrt{2}\Z e_n
\\
\end{array}
&
\widetilde{C}_n
&
2n
\\
D_{n}^{(1)}   \ n\geq 4
&
1, 1,2,\ldots, 2,1,1
&
1, 1,2,\ldots, 2,1,1
&
D_n
&
\begin{array}[t]{@{}ll@{}}
\textcolor{blue!80}{\al_i=e_i-e_{i+1} ,\, 1\leq i\leq n-1}
\\
\textcolor{blue!80}{\al_n=e_{n-1}+e_n}
\end{array}
&
\begin{array}[t]{@{}ll@{}}
\Z\al_1+\cdots+\Z\al_{n-1}+\Z\al_n =
\\
\left\{ 
\be \in \Z e_1+\cdots+\Z e_n
\mid 
\sum_{1\leq i\leq n} \be_i \text{ even } \right\} 
\end{array}
&
\widetilde{D}_n
&
2n-2
\\ 
\begin{array}[t]{@{}ll@{}}
E_{n}^{(1)}\  n=6
\\
\phantom{E_{n}^{(1)}} \  n=7
\\
\phantom{E_{n}^{(1)}}  \  n=8
\end{array}
&
\begin{array}[t]{@{}ll@{}}
1,1,2,3,2,2,1 
\\
1,1,2,3,4,2,3,2
\\
1,2,3,4,5,6,3,4,2
\end{array}
&
\begin{array}[t]{@{}ll@{}}
1,1,2,3,2,2,1 
\\
1,1,2,3,4,2,3,2
\\
1,2,3,4,5,6,3,4,2
\end{array}
&
E_n
&
\begin{array}[t]{@{}ll@{}}
\textcolor{blue!80}{\al_i=e_i-e_{i+1} ,\, 1\leq i\leq n-2}
\\
\textcolor{blue!80}{\al_{n-1}=e_{n-2}+e_{n-1}}
\\
\textcolor{blue!80}{\al_n=-\frac{1}{2}\sum_{1\leq i \leq 8} e_i}
\end{array}
&
\begin{array}[t]{@{}ll@{}}
\Z\al_1+\cdots+\Z\al_n =
\\
\left\{ 
\begin{array}{r}
\be \in \frac{1}{2}\Z e_1+\cdots+\frac{1}{2}\Z e_8
\mid 
\be_i-\be_j \in \Z, 1\leq i,j\leq 8, 
\\
\hfill \sum_{1\leq i\leq n} \be_i \text{ even } , \text{ and }\be_n=\cdots=\be_8
\end{array}\right\} 
\end{array}
&
\widetilde{E}_n
&
\begin{array}[t]{@{}ll@{}}
12
\\
18
\\
30
\end{array}
\\
F_{4}^{(1)}
&
1, 2,3,4,2
&
1, 2,3,2,1	
&
F_4
&\begin{array}[t]{@{}ll@{}}
\textcolor{blue!80}{\al_i=e_i-e_{i+1} ,\, 1\leq i\leq 2}
\\
\al_3=e_3
\\
\al_4=\frac{1}{2}(-e_1-e_2-e_3+e_4)
\end{array}
&
\begin{array}[t]{@{}ll@{}}
\Z\al_1+\Z\al_2+2\Z\al_3+2\Z\al_4 =
\\
\left\{ 
\be \in \Z e_1+\cdots+\Z e_4
\mid 
\sum_{1\leq i\leq 4} \be_i \text{ even}
\right\} 
\end{array}
&
\widetilde{F}_4
&
12
\\ 
G_{2}^{(1)}
&
1,2,3
&
1, 2,1
&
G_2
&
\begin{array}[t]{@{}ll@{}}
\textcolor{blue!80}{\al_1=e_1-e_2}
\\
\al_2=\frac{1}{3}(-2e_1 + e_2 +e_3)
\end{array}
&
\begin{array}[t]{@{}ll@{}}
\Z\al_1+3\Z\al_2 =
\left\{ 
\be \in \Z e_1+\Z e_2+\Z e_3
\mid 
\sum_{1\leq i\leq 3} \be_i =0 \right\} 
\end{array}
&
\widetilde{G}_2
&
6
\\
A_{2n-1}^{(2)}  \ n\geq 3
&
1, 1,2,\ldots, 2,1
&
1, 1,2,\ldots, 2,2
&
C_n
&
\begin{array}[t]{@{}ll@{}}
\al_i=e_i-e_{i+1} ,\, 1\leq i\leq n-1
\\
\textcolor{blue!80}{\al_n=2e_n}
\end{array}
&
\begin{array}[t]{@{}ll@{}}
\Z\al_1+\cdots+\Z\al_{n-1}+\Z\al_n =
\\
\left\{ 
\be \in \Z e_1+\cdots+\Z e_n
\mid 
\sum_{1\leq i\leq n} \be_i \text{ even } \right\} 
\end{array}
&
\widetilde{B}_n
&
2n-1
\\
D_{n+1}^{(2)}  \ n\geq 2
&
1, 1,\ldots, 1,1
&
1, 2,\ldots, 2,1
&
B_n
&
\begin{array}[t]{@{}ll@{}}
\textcolor{blue!80}{\al_i=\sqrt{2}(e_i-e_{i+1}) ,\, 1\leq i\leq n-1}
\\
\al_n=\sqrt{2}e_n
\end{array}
&
\Z\al_1+\cdots+\Z\al_{n-1}+\Z\al_n =
\sqrt{2}\Z e_1+\cdots+\sqrt{2}\Z e_n
&
\widetilde{C}_n
&
n+1
\\
E_{6}^{(2)}
&
1, 2,3,2,1
&
1,2,3,4,2
&
F_4
&
\begin{array}[t]{@{}ll@{}}
\al_i=e_i-e_{i+1} ,\, 1\leq i\leq 2
\\
\textcolor{blue!80}{\al_3=2e_3}
\\
\textcolor{blue!80}{\al_4=-e_1-e_2-e_3+e_4}
\end{array}
&
\begin{array}[t]{@{}ll@{}}
\Z\al_1+\Z\al_2+\Z\al_3+\Z\al_4 =
\\
\left\{ 
\be \in \Z e_1+\cdots+\Z e_4
\mid 
\sum_{1\leq i\leq 4} \be_i \text{ even}
\right\} 
\end{array}
&
\widetilde{F}_4
&
9
\\
D_{4}^{(3)}
&
1, 2,1
&
1, 2,3
&
G_2
&
\begin{array}[t]{@{}ll@{}}
\al_1=e_1-e_2 
\\
\textcolor{blue!80}{\al_2=-2e_1 + e_2 +e_3}
\end{array}
&
\begin{array}[t]{@{}ll@{}}
\Z\al_1+\Z\al_2 = 
\left\{ 
\be \in \Z e_1+\Z e_2+\Z e_3
\mid 
\sum_{1\leq i\leq 3} \be_i =0 \right\} 
\end{array}
&
\widetilde{C}_n
&
4
\\
A_{2n}^{(2)}  \ n\geq 2
&
2, 2,\ldots, 2, 1
&
1, 2,\ldots, 2,2
&
C_n
&
\begin{array}[t]{@{}ll@{}}
\al_i=e_i-e_{i+1} ,\, 1\leq i\leq n-1
\\
\textcolor{blue!80}{\al_n=2e_n}
\end{array}
&
\Z\al_1+\cdots+\Z\al_{n-1}+\frac{1}{2}\Z\al_n =
\Z e_1+\cdots+\Z e_n
&
\widetilde{C}_n
&
2n+1
\\
A_2^{(2)}
&
2,1
&
1,2
&
A_1
&
\textcolor{blue!80}{\al_1=\sqrt{2}(e_1-e_2)}
&
\frac{1}{2}\Z\al_1  =
\left\{ 
\be \in \frac{\sqrt{2}}{2}\Z e_1+\frac{\sqrt{2}}{2}\Z e_2
\mid 
\be_1+\be_2 =0 \right\} 
&
\widetilde{A}_1
&
3
\\ 
\hline
\end{array}
$$
\normalsize
\caption{The different affine Kac-Moody types and associated important data. 
We have indicated in blue the long (real) simple roots.
We have written $\be=\be_1e_1+\be_2e_2+\cdots$ for $\be\in M^{\dagger}$ in the sixth column.}
\label{table}
 \end{figure}
 \end{landscape}
 \clearpage
}